\documentclass[11pt,letterpaper]{article}

\usepackage{graphicx}

\usepackage{ifthen,latexsym,amssymb,amsmath,bbm,fixmath}
\usepackage{color}
\usepackage[procnames]{listings}
\definecolor{keywords}{RGB}{255,0,90}
\definecolor{comments}{RGB}{0,0,113}
\definecolor{myred}{RGB}{160,0,0}
\definecolor{green}{RGB}{0,150,0}
\usepackage[hyphens]{url}

\renewcommand{\baselinestretch}{1.0}
\newcommand{\C}[1]{{\protect\cal #1}}
\newcommand{\B}[1]{{\bf #1}}
\newcommand{\I}[1]{{\mathbb #1}}
\newcommand{\V}[1]{\mathbold{#1}}

\newcommand{\ceil}[1]{\lceil #1\rceil}
\newcommand{\e}{\varepsilon}
\newcommand{\floor}[1]{\lfloor #1\rfloor}

\renewcommand{\mid}{:}

\newif\ifnotesw\noteswtrue
\newcommand{\comment}[1]{\ifnotesw $\blacktriangleright$\ {\sf #1}\
  $\blacktriangleleft$ \fi}
\noteswfalse	
\newcommand{\hide}[1]{}

\newcommand{\beq}[1]{\begin{equation}\label{eq:#1}}
\newcommand{\eeq}{\end{equation}}
\newcommand{\req}[1]{\textrm{(\ref{eq:#1})}}

\newtheorem{theorem}{Theorem}[section]
\newcommand{\bth}[2][nothing]{\ifthenelse{\equal{#1}{nothing}}
 {\begin{theorem}} {\begin{theorem}[#1]}\label{th:#2}}

\newtheorem{lemma}[theorem]{Lemma}
\newcommand{\blm}[2][nothing]{\ifthenelse{\equal{#1}{nothing}}
 {\begin{lemma}} {\begin{lemma}[#1]}\label{lm:#2}}

\newtheorem{problem}[theorem]{Problem}
\newcommand{\bpr}[2][nothing]{\ifthenelse{\equal{#1}{nothing}}
 {\begin{problem}} {\begin{problem}[#1]}\label{pr:#2}}

\newtheorem{corollary}[theorem]{Corollary}
\newtheorem{definition}[theorem]{Definition}

\newtheorem{claim}[theorem]{Claim}

\newcommand{\Claim}[2]{\begin{claim}\label{cl:#1} #2\end{claim}}

\newcommand{\bpf}[1][Proof.]{\smallskip\noindent{\it #1} }
\newcommand{\qed}{\nolinebreak\mbox{\hspace{5 true pt}%
  \rule[-0.85 true pt]{3.9 true pt}{8.1 true pt}}}
\newcommand{\cqed}{\nolinebreak\mbox{\hspace{5 true pt}%
  \rule[-0.85 true pt]{2.0 true pt}{8.1 true pt}}}
\newcommand{\epf}{\qed \medskip}
\newcommand{\bcpf}{\bpf[Proof of Claim.]}
\newcommand{\ecpf}{\cqed \medskip}

\newcommand{\brm}{\smallskip\noindent{\bf Remark.} }

\newcommand{\blow}[2]{#1(#2)}

\newcommand{\pp}{P}

\newcommand{\Forb}[1]{\mathrm{Forb}(#1)}
\newcommand{\dedit}{\delta_{\mathrm{edit}}}
\newcommand{\Dedit}{\Delta_{\mathrm{edit}}}
\newcommand{\E}{\B E}
\newtheorem{assumption}[theorem]{Assumption}
\newtheorem{example}[theorem]{Example}
\newcommand{\g}[1]{\gamma(#1)}
\newcommand{\gmax}{\|\gamma\|_\infty}
\newcommand{\ex}{\mathrm{ex}}

\renewcommand{\dots}{\hspace{0.9pt}.\hspace{0.3pt}.\hspace{0.3pt}.\hspace{1.5pt}
}
\renewcommand{\geq}{\geqslant}
\renewcommand{\leq}{\leqslant}
\renewcommand{\ge}{\geqslant}
\renewcommand{\le}{\leqslant}

\begin{document}

\title{Strong Forms of Stability from Flag Algebra Calculations}
\author{Oleg Pikhurko\mbox{$^{1,}$}\thanks{Supported by ERC
grant~306493 and EPSRC grant~EP/K012045/1.}\and Jakub
Slia\v{c}an\mbox{$^{2,}$}\footnote{Part of this research was done when this author
was supported by ERC
grant~306493.}\and
Konstantinos Tyros\mbox{$^{3,}$}\footnote{Part of this research was done when this author
was supported by ERC
grant~306493.}}
\footnotetext[1]{Mathematics Institute and DIMAP,
University of Warwick,
Coventry CV4 7AL, UK}
\footnotetext[2]{Department of Mathematics and Statistic, The Open University,
Milton Keynes,
MK7 6AA, UK}
\footnotetext[3]{Department of Mathematics, Ko\c c University, Saryer, Istanbul, 34450, Turkey}
\date{}
\maketitle

\begin{abstract}
 Given a hereditary family $\C G$ of \emph{admissible} graphs and a function $\lambda(G)$
that linearly depends on the statistics of order-$\kappa$ subgraphs
in a graph $G$, we consider the extremal problem of determining $\lambda(n,\C G)$,
the maximum of $\lambda(G)$ over all admissible graphs $G$ of order~$n$.
We call the problem
\emph{perfectly $B$-stable} for a graph $B$ if there is a constant $C$ such that \textbf{every}
admissible graph $G$ of order $n\ge C$
can be made into a blow-up of $B$ by changing at most $C(\lambda(n,\C G)-\lambda(G)){n\choose2}$ adjacencies. As special cases,
this property describes all almost extremal graphs of order $n$  within $o(n^2)$ edges 
and shows that
every extremal graph of order $n\ge n_0$ is a blow-up of~$B$.

We develop general methods for establishing stability-type results from
flag algebra computations and apply them to concrete examples.
In fact, one of our sufficient
conditions for perfect stability is stated in a
way
that allows automatic verification by a computer. This gives
a unifying way to obtain computer-assisted proofs of many new results.
\end{abstract}


\section{Introduction}\label{Intro}

By the term graph, we mean finite \emph{simple} graph,
that is, without loops and multiple edges. For a graph $G$, we refer to the
cardinality of its vertex set as the \emph{order} of $G$ and we denote it by~$v(G)$.
Moreover, for  a subset $X$ of the vertex set $V(G)$ of $G$, we denote by $G[X]$ the
subgraph \emph{induced} by $X$ in  $G$, that is, the graph having $X$ as the vertex
set and two nodes $x,y\in{}X$ are connected in $G[X]$ if and only if they are
connected in $G$.

Let $F$ and $G$ be graphs of orders $v(F)\le v(G)$. Call $F$ a \emph{subgraph} of
$G$ if there is a subset $X$ of $V(G)$ such that $G[X]$ is isomorphic to $F$. (Thus a subgraph means
an induced subgraph.)
Let $P(F,G)$ be the number of $v(F)$-subsets of $V(G)$ that induce
a subgraph isomorphic to $F$.
Also, let $$
 p(F,G):=P(F,G)\Big/{v(G)\choose v(F)}$$ be the
\emph{(induced) density} of $F$ in $G$; equivalently, $p(F,G)$ is the
probability that
a random $v(F)$-subset $X$ of $V(G)$ induces a subgraph isomorphic to~$F$.

Suppose that we have a (possibly infinite) family $\C F$ of \emph{forbidden} graphs. Call a graph $G$ \emph{admissible} or \emph{$\C F$-free}
(and denote this by $G\in\Forb{\C F}$)
if no $F\in{}\C F$ is a subgraph of $G$.
Let $\C G:=\Forb{\C F}$ be the family of all admissible graphs; clearly,
$\C G$ is a hereditary graph property, that is, every subgraph of some member of $\C G$ belongs to $\C G$, too.

Let $\kappa$ be a positive integer. We denote by $\C G_{\kappa}^0$ the set obtained by taking one representative from
each isomorphism class of
graphs in $\C G$ of order $\kappa$. Clearly $\C G_{\kappa}^0$ is a finite set. Let $\gamma$ be a function from $\C G_\kappa^0$ to the reals.
It gives rise to two other functions defined on graphs $G$ with $v(G)\ge \kappa$:
\begin{equation}
  \label{eq:lambdadefn}
  \begin{split}
    \Lambda(G)&:=\sum_{H\in\C G_\kappa^0} \g H P(H,G),\\
 \lambda(G)&:=\sum_{H\in\C G_\kappa^0} \g H p(H,G)={v(G)\choose \kappa}^{-1}\cdot\Lambda(G).
  \end{split}
\end{equation}
  One can view $\lambda(G)$
as the expected value of $\g{G[X]}$, where $X$ is a random $\kappa$-subset
of $V(G)$.
 \hide{It may be convenient to extend $\gamma$ on all
graphs of order $\kappa$ by agreeing that $\g H=0$ for
every non-admissible $H$.}

Under the above conventions, consider the  problem of maximising
$\Lambda(G)$ over admissible graphs $G$ of given order $n$. Namely, we define
the extremal function
 \begin{equation}\label{eq:GenProblem}
 \Lambda(n,\C G):=\max\{\Lambda(G)\mid G\in\C G,\ v(G)=n\}
 \end{equation}
 and its density version $\lambda(n,\C G):=\Lambda(n,\C G)/{n\choose \kappa}$.
  It is not hard to show
  (see Lemma~\ref{lm:lim}) that the sequence $(\lambda(n,\C G))_{n=\kappa}^\infty$ is
  non-increasing and therefore the following limit exists:
   \begin{equation}\label{eq:GenProblemAsympt}
   \lambda(\C G):=\lim_{n\to\infty} \lambda(n,\C G).
   \end{equation}

This is a rather general setting. As an illustration, here is one example (and
the reader is encouraged to consult
other concrete examples that can be found in Section~\ref{se:examples}).

\begin{example}[Tur\'{a}n function]\label{ex:TuranFunction}
 Let $\kappa=2$, $\g{\overline K_2}=0$
and $\g{K_2}=1$, where by $K_m$ we denote the complete graph of order $m$
and by $\overline G$ the complement of a graph $G$.
(Thus $\Lambda(G)=e(G)$ is the number of edges in $G$.)
If $\C H$ is any family of graphs and $\C H^\uparrow\supseteq \C H$ consists
of all graphs that can be obtained from $H\in\C H$ by adding some
missing edges, then
$\Lambda(n,\Forb{\C H^\uparrow})$ is the well-known Tur\'an function
$\ex(n,\C H)$.\end{example}

\hide{
Since the case of small $n$ is often exceptional and difficult to deal with,
we are interested in understanding the problem $n\to\infty$.
Informally speaking, we are going to present some general results that directly imply or help
with proving that the extremal function $\Lambda(n,\C G)$ is given by a blow-up of
fixed graph $B$ and, furthermore, all almost extremal graphs $G$ of order $n\to\infty$
(those with $\Lambda(G)=\Lambda(n,\C G)+o(n^{\kappa})$) are $o(n^2)$ close in edit distance to some blow-up of~$B$.
}%

Fix a graph $B$ with vertex set $[m]:=\{1,\dots,m\}$. For pairwise disjoint sets $V_1,\dots,V_m$
(some of which may be empty), let
the \emph{blow-up} $\blow{B}{V_1,\dots,V_m}$ be obtained from the empty graph on $V_1\cup\dots\cup V_m$ by adding for every edge $\{i,j\}\in E(B)$ the complete
bipartite graph with parts $V_i$ and $V_j$. Note that no part $V_i$ spans an edge. Let $\blow{B}{}$ be the family of
all possible blow-ups of $B$. It consists of graphs
that can be obtained from $B$ by a sequence of vertex duplications and vertex deletions.

Suppose that $\blow{B}{}\subseteq \C G$.
Trivially, we get the lower bound $\Lambda(n,\C G)\ge \Lambda(n,\blow{B}{})$ valid
for every integer $n\ge \kappa$, where, in accordance with our general notation, $\Lambda(n,\blow{B}{})$ is the maximum value of $\Lambda$
over all blow-ups of  $B$ of order $n$. For a vector
$\V a=(a_1,\dots,a_m)$ in
the \emph{$m$-simplex}
 $$
 \I S_{m}:=\{\V x \in \I R^{m}:\V x\ge 0,\ x_1+\dots+x_m=1\},
 $$
 define
 \begin{equation}\label{eq:lambdaBA}
 \lambda(\blow{B}{\V a}):=\lim_{n\to\infty}
 \lambda(\blow{B}{V_{1,n},\dots,V_{m,n}}),
 \end{equation}
 where $|V_{i,n}|=a_in+O(1)$ for each $i\in[m]$. In other words,
we look at the limiting value of the function
$\lambda$ evaluated on a blow-up $G$ of $B$ of order $n\to\infty$ where
the $i$-th part occupies
$(a_i+O(1/n))$-fraction
of all vertices. It is easy to see that the limit in~\eqref{eq:lambdaBA} exists
(and  does not depend on the choice of the sizes $|V_{i,n}|$).
 In fact, $\lambda(\blow{B}{\V
a})$ is a polynomial in $\V a$ of degree
at most~$\kappa$, so the rate of convergence in~\eqref{eq:lambdaBA} is
$O(1/n)$. An easy argument based on the compactness of $\I S^m$ and
the continuity of $\lambda(\blow{B}{\V a})$ as a function of $\V a\in \I S_m$ shows that
 \begin{equation}\label{eq:MaxLambdaBlow}
  \lambda(\blow{B}{})=\max\{\lambda(\blow{B}{\V x})\mid \V x\in\I S_m\}.
 \end{equation}
 By our assumption, $\blow{B}{}\subseteq \C G$ so $\lambda(\blow{B}{})$ is a
 lower bound on
$\lambda(\C G)$.

Here is an illustration.  In
Example~\ref{ex:TuranFunction}, if $\C
H=\C H^\uparrow=\{K_t\}$ consists of the $t$-clique,
then a good choice is to take $B=K_{t-1}$. Then $\blow{B}{}$ is a subset of
$\C G=\Forb{\C F}$ and
$\lambda(\blow{B}{\V a})=2\prod_{1\le i<j\le t-1} a_ia_j$. This is clearly
maximised if all entries $a_i$ are equal to each other, giving 
$\lambda(\blow{B}{})= (t-2)/(t-1)$, which is a lower bound on the Tur\'an density of $K_t$. The
classical
results
of Mantel~\cite{mantel:07} (for $t=3$) and Tur\'an~\cite{turan:41} (for any
$t$) imply that this is an equality. More strongly, they showed that
$\Lambda(n,\C G)=\Lambda(n,\blow{B}{})$ for all $n$, while an easy optimisation
shows that $\Lambda(n,\blow{B}{})$ is attained
by the (unique) blow-up of $B$ of order $n$ with parts as equal as possible.

In general, there is no hope for a theory that allows to determine
$\lambda(\C G)$ for every $\lambda$ and $\C G$. Namely, as it was shown by
Hatami and Norin~\cite{hatami+norine:11},
the question if $\lambda(\C G)\le c$ on input $(\C G,\Lambda,c)$ is undecidable in
general (even if $\C G$ consists of all graphs).
However, one can determine $\lambda(\C G)$ for various concrete examples
of interest.
Many of these proofs utilise the powerful flag algebra approach of Razborov~\cite{razborov:07,razborov:10}, where a computer can be used to generate
a \emph{certificate} $\C C$ proving $\lambda(\C G)\le c$.

We will \textsl{}discuss certificates
in more detail in Section~\ref{certificates}. For the time being, let us just
remark that the desired
inequality $\lambda(\C G)\le c$ follows by symbolically representing
$c-\lambda(G)$, for $G\in \C G_n^0$
as a sum of squares within error term of order $O(1/n)$ as $n\to\infty$.
One illustrative example of such a sum is $s(G):={n\choose 3}^{-1}\sum_{x\in
V(G)}
(d_G(x)-d_{\overline G}(x))^2$, where $d_H(x)$ is the degree of $x$ in a graph~$H$.
Clearly, $s(G)$ is non-negative while it is routine to see that
$s(G)=6p(K_3,G)+6p(\overline{K_3},G)-2p(P_3,G)-2p(\overline{P_3},G)+O(1/n)$, where $P_i$ is the path with $i$ vertices. This gives
an asymptotic inequality that always holds between 3-vertex subgraph densities.
One advantage of the flag algebra approach is that it allows us to
generate and manipulate
such equalities automatically; here finding optimal coefficients amounts
to solving a certain semi-definite programme (which is independent of $n$).
We refer the reader to Section~\ref{certificates} for details and formal
definitions.

Thus, if a flag algebra calculation proves $\lambda(\C G)\le c$ while we can find an order-$m$ graph $B$ with
$\blow{B}{}\subseteq \C G$ and
$\V a\in \I S_m$ with $\lambda(\blow{B}{\V a})=c$, then we know
$\lambda(n,\C G)$ within an error term of~$O(1/n)$:
 $$
  c-O(1/n)\le \lambda(n,\blow{B}{\V a})\le \lambda(n,\C
 G)\le c+O(1/n).
  $$

\hide{When
the flag algebras and blow-ups gave the value of $\lambda(\C G)$ as above, the authors
were typically able to close the gap; however, previous approaches were rather
ad-hoc and dependent on the concrete choice of $\Lambda$ and $\C G$. Here we
present some results that help to prove that $\Lambda(n,\C G)=\Lambda(n,\blow{B}{})$
for all large $n$. Our methods also helps with establishing the following very strong
stability properties.}

In addition to determining $\lambda(\C G)$, it is often desirable to obtain information on the structure
of all large admissible graphs $G$ for which the value $\lambda(G)$ is close to the
maximal possible. In particular, we look for sufficient conditions establishing that
every such $G$ is necessarily close to a blow-up of $B$, in which case we regard
the problem as stable. This paper will consider a few non-equivalent
versions of stability, with the corresponding definitions following shortly. The
stability is a very useful property in extremal graph theory as it is often
indispensable
in determining the exact value of $\lambda(n, \C G)$ as well as the set of all extremal graphs of
large order~$n$.
Besides being an important property on its own, stability also
helps in solving the randomised or counting versions of extremal results.

We will use only one notion of distance on graphs. Namely, the \emph{(edit) distance} between two graphs $G$ and $H$ of the same order $n$ is
 $$
 \Dedit(G,H):= \min_\sigma \Big|\,E(G)\bigtriangleup \big\{\{\sigma(u),\sigma(v)\}\mid uv\in E(H)\big\}\,\Big|,
 $$
  where the minimum is taken over all bijections $\sigma:V(H)\to V(G)$. In other words, $\Dedit(G,H)$ is the minimum number of adjacencies that one needs to change in $G$ in order to obtain a graph isomorphic to $H$. We define the \emph{normalised (edit) distance} to be $\dedit(G,H):=\Dedit(G,H)/{n\choose 2}$. For a family $\C H$ of graphs we define $\Dedit(G,\C H):=\min\{\Dedit(G,H)\mid H\in \C H_n^0\}$ and $\dedit(G,\C H):=\min\{\dedit(G,H)\mid H\in \C H_n^0\}$.

We say that our problem~\eqref{eq:GenProblem} is \emph{robustly $B$-stable}
(resp.\ \emph{perfectly $B$-stable})
if there is $C>0$ such that for every graph $G\in\C G$ of order $n\ge C$ we have $$
 \dedit(G,\blow{B}{})\le  C \max\left(1/{n},\lambda(n,\C G)-\lambda(G)\right),
 $$
 (resp.\ $\dedit(G,\blow{B}{})\le  C (\lambda(n,\C G)-\lambda(G))$).
 For comparison,
the \emph{classical $B$-stability} states that for every $\e>0$ there is $\delta>0$ such that $\dedit(G,\blow{B}{})\le \e$
for
every $G\in\C G$ with $n\ge 1/\delta$ vertices and $\lambda(G)\ge
\lambda(\C G)-\delta$.
Clearly, the perfect stability implies the
robust stability which in turn implies the classical stability. (Our
Theorem~\ref{th:turan} will in particular
show
that these notions of stability are not equivalent, already for
such a natural problem as the Tur\'an function.)
Also,
if the problem is perfectly stable, then  for all $n\ge C$
we have $\Lambda(n,\C G)=\Lambda(n,\blow{B}{})$ and every extremal graph is
a blow-up of $B$ (which one may call
an \emph{exact result}).

For example, for the Tur\'an function $\ex(n,\C F)$, the classical stability
was established independently by Erd\H os~\cite{erdos:67} and
Simonovits~\cite{simonovits:68}. The perfect (and thus also robust)
stability in the case when $\C F$ consists of a clique $K_t$ follows from results of F{\"u}redi~\cite{furedi:15} (who considered the distance to being $(t-1)$-partite instead of complete $(t-1)$-partite as we do in this paper). Very recently, Roberts and Scott~\cite{roberts+scott} improved on F\"uredi's result by extending it to
all colour critical graphs and giving a sharper bound on the distance.

As far as we know, the above  results in~\cite{furedi:15,roberts+scott} and some recent work of  Norin and
Yepremyan~\cite{norin+yepremyan:extensions,NorinYepremyan17}
(who considered the Tur\'an problem for hypergraphs) are the only known
examples where perfect stability was established for a non-trivial problem.
Furthermore, almost all proofs where the classical stability and the exact
result were derived from a flag algebra computation were rather ad-hoc and
taylored to a particular problem.

The purpose of this paper is to present
some general sufficient conditions that imply some version of stability.
This allows us to give a unified proof of many previous stability and exactness
results. Also, we can establish perfect
stability (a new result) for a number of problems.

More specifically,
Theorem~\ref{th:stab1} gives a sufficient condition for robust stability
and Theorems~\ref{th:exact}, \ref{th:exact2}, and~\ref{thm:pc} give various
sufficient conditions for perfect stability. Furthermore, all assumptions
of Theorem~\ref{thm:pc} are stated in a way that allows automatic verification
by a computer. We also present an openly available computer code (written in
\texttt{sage} by adopting the \texttt{flagmatic} package of Emil
Vaughan~\cite{flagmatic:1.5}) that allows us to both generate and verify
certificates for
general problems based on Theorem~\ref{thm:pc}. In all the cases where we could
verify assumptions of
Theorem~\ref{thm:pc} and derive perfect stability, the procedure
was essentially mechanical and the final high-level code
is very short (having 6--10 lines, each invoking some function).

Even if one knows that $\lambda(\C G)=\lambda(\blow{B}{})$ for a concrete $B$,
the determination of asymptotically optimal part ratios (namely, finding all
$\V a\in \I S_m$ with $\lambda(\blow{B}{})=\lambda(\blow{B}{\V a})$) may
still be a non-trivial task that amounts to polynomial maximisation. While the combination of
Lagrange multipliers and Gr\"obner bases provides a general computational framework,
in an extremal problem one often has a candidate $\V a\in\I S_m$ and
wishes to prove that this is the only vector (up to a symmetry of $B$)
that achieves $\lambda(\blow{B}{})$. Clearly, if a flag algebra certificate $\C
C$
proves that $\lambda(\C G)\le \lambda(\blow{B}{\V a})$ then this automatically
implies that $\V a$ is a maximiser and it is possible that the information
in $\C C$ is enough to imply the uniqueness of~$\V a$. We present such
a sufficient condition on $\C C$ in Lemma~\ref{lm:unique} which can
be automatically verified by a computer and seems to work very well in practice.

The exact statements of the above sufficient conditions rely on understanding
flag
algebra  certificates, so we postpone them until later. Here, let us list the
extremal problems for which our method gives perfect
stability. In almost all cases, Theorem~\ref{thm:pc} and
Lemma~\ref{lm:unique}
apply directly, immediately giving perfect
stability and implying the uniqueness of asymptotically optimal part ratios.

However, there are a few natural problems where the assumptions of
Theorem~\ref{thm:pc} are not satisfied. As the test case that our method
can still give perfect stability, we chose the inducibility function for
the paw graph, see Theorem~\ref{th:paw}. The asymptotic value of this
function was determined by Hirst~\cite{hirst-inducibility} but the classical stability
and the exact result were not known. By utilising our other results
(Theorem \ref{th:exact}) we derive perfect stability for the inducibility
problem for the paw graph. Since
the proof is rather long and was meant mainly as an illustration of
the flexibility of our method, we decided to include only one such example in this paper.

\subsection{Examples of extremal problems for which we can prove perfect
stability}\label{se:examples}

\subsubsection{Minimising the number of independent sets in
triangle-free graphs}\label{se:fkl}

Erd\H os~\cite{erdos:62} asked for the value of $f(n,k,l)$, the minimum
number of independent sets of size $k$ in a $K_l$-free graph of order $n$.

Consider first the case $l=3$, when we forbid a triangle.
Goodman~\cite{goodman:59} determined $f(2n,3,3)$; his bounds
also give the asymptotic value of $f(2n+1,3,3)$.
Lorden~\cite{lorden:62} showed that, for $n\ge 12$, the value of $f(n,3,3)$ is attained by taking $K_{\floor{n/2},\ceil{n/2}}$ and removing any (possibly empty) matching from it. Some partial results were
obtained by Nikoforov~\cite{nikiforov:05:arxiv,nikiforov:01}.
The problem for $k\ge 4$ remained
open until recently when the classical stability and the exact result were
established
independently by Das, Huang, Ma, Naves and
Sudakov~\cite{das+huang+ma+naves+sudakov:13} (for $k\in\{4,5\}$)
and Pikhurko and Vaughan~\cite{pikhurko+vaughan:13} (for
$k\in\{4,5,6,7\}$) when $n$ is sufficiently large: if $k\in\{4,5\}$,
then all extremal graphs are blow-ups of $C_5$ and if $k\in \{6,7\}$, then
all extremal graphs are blow-ups of the \emph{Clebsch
graph}~$L$. The \emph{Clebsch
graph} $L$ has, as  vertices, binary $5$-sequences of even \emph{weight} (that is, the number of ones), with
two vertices being adjacent
if the point-wise sum modulo 2 of the corresponding sequences has weight $4$.
It easily follows from this description that the graph $L$ has 16 vertices and is triangle-free and
$5$-regular.

This question of Erd\H os is a special case of our general problem.
It turns out that our computer codes can prove the perfect stability
in the following cases.
(Note that the $f(n,3,3)$-problem is not perfectly stable by the above mentioned result of Lorden~\cite{lorden:62}.)

\begin{theorem}\label{th:ErdosPr3} Let
 \begin{itemize}
 \item $k\in \{4,5\}$ and $B=C_5$, or
 \item $k\in\{6,7\}$ and $B=L$.
 \end{itemize}
 Let $\C F=\{K_3\}$ (thus $\C G$ consist of all triangle-free graphs) and let $\gamma(H)$ be 0 except
 $\gamma(\overline K_k)=-1$.
 Then the corresponding problem $\Lambda(n,\C G)$ (that is,  Erd\H os' problem
 of determining
 $f(n,k,3)$) is perfectly $B$-stable. Furthermore, for each $k\in\{4,\dots,7\}$
 the unique
probability vector $\V a\in \I S_{v(B)}$ that maximises $\lambda(\blow{B}{\V
a})$ is the uniform vector.
 \end{theorem}

If $l\ge 4$, then the asymptotic value of $f(n,k,l)$ is known only for $k=3$
and $l\le 7$,
see~\cite{das+huang+ma+naves+sudakov:13,pikhurko+vaughan:13}.
The papers~\cite{das+huang+ma+naves+sudakov:13,pikhurko+vaughan:13} also showed
that in each of these cases the problem is classically $B$-stable,
where $B=K_{l-1}$, and the value of
$f(n,k,l)$ is attained by a blow-up of $B$ for all large $n$. However,
the problem is not perfectly $B$-stable since it is possible to remove a few edges from the blow-up of $B$ (e.g.\ a matching between two parts) so that the number of copies of $\overline K_3$ does not change.

The above results and our Theorem~\ref{th:exact2} imply that, in fact, the $f(n,3,l)$-problem is not robustly stable for
$l\in\{3,4,5,6,7\}$. Alternatively, the same conclusion can be derived directly by taking the optimal blow-up $\blow{K_{l-1}}{V_1,\dots,V_{l-1}}$, fixing some sets $X\subseteq V_1$ and $Y\subseteq V_2$ each of size $\e n$, where $\e$ is a small constant, and then flipping all pairs between $X$ and~$Y$. (Such tranformation is done in the proof of Theorem~\ref{th:exact2} and is carefully analysed there.) 

\subsubsection{Maximising the number of pentagons in triangle-free graphs}

Erd\H os~\cite{erdos:84} asked if $c(5m)\le m^5$ for every natural number $m$,
where
$c(n)$ is the maximum number of 5-cycles that a triangle-free graph
of order $n$ can have. Note that this bound
is sharp for every $m\in\I N$ which is witnessed by the balanced blow-up
of $C_5$.

Some partial progress on this problem was obtained by
Gy\H{o}ri~\cite{gyori:89} who proved $c(n)\le 1.03\, (n/5)^5$ and
F\"uredi (unpublished) who proved $c(n)\le 1.01\, (n/5)^5$. Recently,
Grzesik~\cite{grzesik-pentagon} and independently
Hatami et al.~\cite{hatami-hladky-kral-pentagon} proved that, as $n\to\infty$,
there can be at most $(1/5^5+o(1))n^5$ copies of $C_5$. Furthermore, Hatami et al.~\cite{hatami-hladky-kral-pentagon} proved the exact result for all large $n$
(and the
classical stability can also be derived from their method). In fact,
the validity of Erd\H os' conjecture follows from the asymptotic result by a simple
blow-up trick (see~\cite[Corollary~3.3]{hatami-hladky-kral-pentagon}).
Interestingly, if $n=8$, there is another extremal example which is not a
blow-up of $C_5$ that was discovered by Michael~\cite{michael:11}. Very recently,
the value of $c(n)$ and the description of all extremal graphs for every $n$ was obtained 
by Lidick\'{y} and Pfender~\cite{LidickiPfender18arxiv}.

This problem fits into our general framework and we can prove (again
in a completely automated way) that it is perfectly stable.

\begin{theorem}
\label{th:pentagons}
Let $\C F=\{K_3\}$, $\kappa=5$, and $\gamma(H)$ be zero,
except $\gamma(C_5)=1$. Let $B=C_5$. Then the corresponding problem is
perfectly $B$-stable
(with the unique maximiser for $\lambda(\blow{B}{\V a})$ being
the vector $\V a\in\I S_5$ with each entry equal to $1/5$).
\end{theorem}

\subsubsection{Inducibility}

Given a graph $F$, the \emph{inducibility
problem} for $F$ asks for the maximal possible (induced) density of the
graph $F$ among all graphs of order $n$. In our general notation,
it can be expressed as follows. Let $\kappa=v(F)$, $\C F=\emptyset$ and
let $\gamma$ take
the value 0 on every graph with $\kappa$ vertices except $F$, where it
takes value 1. Thus we are interested in $i(n,F):=\Lambda(n,\C G)$. We call
$i(F):=\lambda(\C
G)$ the
\emph{inducibility} of $F$. Equivalently,
 $$
 i(F):= \lim_{n\to\infty} \max\{p(F,G)\mid v(G)=n\}.
 $$
Observe that the inducibility of each graph is equal to the inducibility of its complement.%
\hide{
In \cite{Pippenger+Golumbic:75} Pippenger and Golumbic proved that this is the only case, showing that each graph $F$,
that is neither complete nor empty, satisfies $1>i(F)\ge k!/(k^k-k)$, where $k$ denotes the cardinatily of the vertex set of $F$.
Moreover, they provide lower bounds for the inducibility of cycles and they determine the inducibility of complete bipartite graphs
which are almost balanced, that is, the cardinalities of their parts differ by at most 1.
In \cite{siran:84} \v{S}ir\'a\v{n}, improved on the lower bound for the inducibility of a graph.
The inducibility problem drawn a considerable amount of interest \cite{Balogh+Hu+Lidicky+Pfender:16,evenzohar+linial:16,Hatami+Hirst+Norine:14,hirst-inducibility},
while, the inducibility of partite graphs, in particular, triggered a big
amount of research \cite{Bollobas+Egawa+Jin:95,bollobas+nara+tachibana:86,
Brown+Sidorenko:94,exoo:86}.
Below we list inducibility problems that we show that they are perfectly stable.
}
The inducibility problem has drawn a considerable amount of interest, see for
example
\cite{Balogh+Hu+Lidicky+Pfender:16,Bollobas+Egawa+Jin:95,
bollobas+nara+tachibana:86,Brown+Sidorenko:94,evenzohar+linial:16,exoo:86,
Pippenger+Golumbic:75,
Hatami+Hirst+Norine:14,HefetzTyomkyn17arxiv,hirst-inducibility,siran:84}.

Before we look at
concrete examples, let us mention the following general result of Brown
and Sidorenko~\cite[Proposition~1]{Brown+Sidorenko:94}: if $F$ is
\emph{complete partite} (i.e.\ a blow-up of some clique), then for every $n\in\I N$ at
least one
$i(n,F)$-extremal graph is complete partite. The proof
in~\cite{Brown+Sidorenko:94} uses the
symmetrisation method and it is not clear how to extract a stability-type
result from it. Also, Even-Zohar and Linial~\cite[Table~2]{evenzohar+linial:16}
systematically looked at the inducibility of 5-vertex graphs but without
trying to convert the numerical bounds coming from flag algebra calculations
into
computer-verifiable
mathematical proofs.

\subsubsection{Inducibility of the cycle on four vertices}

The inducibility of the 4-cycle, denoted by $C_4$, follows from the
above mentioned result of  Brown
and Sidorenko~\cite{Brown+Sidorenko:94}. Previously,  Pippenger 
and Golumbic
\cite{Pippenger+Golumbic:75} determined $i(n,K_{k,l})$ for all $k,l$ with
$|k-l|\le 1$, observing that the complete balanced bipartite graph is an
extremal graph.
Here we prove perfect stability for this problem (by invoking our computer
code).

\begin{theorem}
\label{th:c4}
Let $\C F=\emptyset$, $\kappa=4$, and $\gamma(H)$ be zero,
except $\gamma(C_4)=1$. Let $B=K_2$ be a single edge. Then the corresponding problem is
perfectly $B$-stable
(with $\lambda(\C G)=i(C_4)=3/8$ and  the unique maximiser for $\lambda(\blow{B}{\V a})$ being
the vector $\V a=(1/2,1/2)$).
\end{theorem}

\subsubsection{Inducibility of $K_4$ minus an edge}
Let $K_4^-$ be the graph obtained by removing an edge from the complete
graph on four vertices. The inducibility problem for $K_4^-$ was considered
by Hirst~\cite{hirst-inducibility}, who determined
$i(K_4^-)$ using
the flag algebra method. Our new result is that this problem is perfectly
stable (and,
in particular, that $i(n,K_4^-)$ is attained by a blow-up of the complete graph on five vertices $K_5$,
for all large~$n$).

\begin{theorem}
\label{th:k4-}
Let $\C F=\emptyset$, $\kappa=4$, and $\gamma(H)$ be zero,
except $\gamma(K_4^-)=1$. Also let $B=K_5$. Then the corresponding problem is
perfectly $B$-stable
(with $\lambda(\C G)=i(K_4^-)=3/8$ and the unique maximiser for $\lambda(\blow{B}{\V a})$ being
the vector $\V a\in\I S_5$ with each entry equal to $1/5$).
\end{theorem}

\subsubsection{Inducibility of $K_{3,2}$}

The function $i(n,K_{3,2})$ was calculated by Golumbic and Pippenger in
\cite{Pippenger+Golumbic:75}, where the complete balanced bipartite graph
is an extremal graph. We show that the problem is perfectly stable.

\begin{theorem}
\label{th:k32}
Let $\C F=\emptyset$, $\kappa=5$, and $\gamma(H)$ be zero,
except $\gamma(K_{3,2})=1$. Also let $B=K_2$. Then the corresponding problem is
perfectly $B$-stable
(with $\lambda(\C G)=i(K_{3,2})=5/8$ and the unique maximiser for $\lambda(\blow{B}{\V a})$ being
the vector $\V a=(1/2,1/2)$).
\end{theorem}

\subsubsection{Inducibility of $K_{2,2,1}$}

The function $i(n,K_{2,2,1})$ can be derived from the results of Brown and
Sidorenko~\cite{Brown+Sidorenko:94}.
Here we prove perfect stability for the corresponding problem.

\begin{theorem}
\label{th:k221}
Let $\C F=\emptyset$, $\kappa=5$, and $\gamma(H)$ be zero,
except $\gamma(K_{2,2,1})=1$. Also let $B=K_3$. Then the corresponding problem is
perfectly $B$-stable
(with $\lambda(\C G)=i(K_{2,2,1})=10/27$ and the unique maximiser for $\lambda(\blow{B}{\V a})$ being
the vector $\V a=(1/3,1/3,1/3)$).
\end{theorem}

\subsubsection{Inducibility of $P_3\cup K_2$}
We consider the inducibility problem for the disjoint union of a path on 3 vertices and an edge, which we denote by $P_3\cup K_2$.
We prove the following.
\begin{theorem}
\label{th:p3k2}
Let $\C F=\emptyset$, $\kappa=5$, and $\gamma(H)$ be zero,
except $\gamma(P_3\cup K_2)=1$. Also let $B=K_3\cup K_3$, that is, the disjoint union of two triangles.
Then the corresponding problem is perfectly $B$-stable
(with $\lambda(\C G)=i(P_3\cup K_2)=5/18$ and the unique maximiser for $\lambda(\blow{B}{\V a})$ being
the vector $\V a\in\I S_6$ with each entry equal to $1/6$).
\end{theorem}

\subsubsection{Inducibility of the ``Y'' graph}
\begin{figure}[htb]
\centering \includegraphics[]{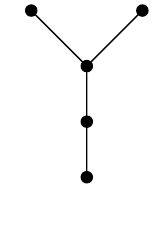}
\caption{The graph $Y$
.}\label{fg:Y}
\end{figure}

We consider the inducibility problem for the graph depicted in
Figure~\ref{fg:Y}, which we denote by $\mathrm{Y}$.
We prove the following.
\begin{theorem}
\label{th:Y}
Let $\C F=\emptyset$, $\kappa=5$, and $\gamma(H)$ be zero,
except $\gamma(\mathrm{Y})=1$. Also let $B=C_5$.
Then the corresponding problem is perfectly $B$-stable
(with $\lambda(\C G)=i(Y)=24/125$ and the unique maximiser for $\lambda(\blow{B}{\V a})$ being
the vector $\V a\in\I S_5$ with each entry equal to $1/5$).
\end{theorem}

\subsubsection{Inducibility for the paw graph}
Let us denote by $F_{\mathrm{paw}}$ the graph obtained by adding to a triangle a pendant edge.
Using flag algebras, Hirst~\cite{hirst-inducibility} determinend the value of
$i(F_{\text{paw}})$.
Here, it is more convenient to work with the (equivalent) complementary
problem. Thus we consider the inducibility
problem of the
disjoint union of a path on three vertices and a single node, that we denote by $P_3\cup K_1$.
Unfortunately, the perfect stability does not follow directly from Theorem
\ref{thm:pc}. However, it can be proved using our methods
combined with additional work, see Section \ref{paw} for the proof.

\begin{theorem}
\label{th:paw}
Let $\C F=\emptyset$, $\kappa=4$, and $\gamma(H)$ be zero,
except $\gamma(P_3\cup K_1)=1$. Also let $B=K_2\cup K_2$, that is, the
disjoint union of two edges.
Then the corresponding problem is perfectly $B$-stable
(with $\lambda(\C G)=i(P_3\cup K_1)=i(F_\mathrm{paw})=3/8$ and the unique maximiser for $\lambda(\blow{B}{\V a})$ being
the vector $\V a\in\I S_4$ with each entry equal to $1/4$).
\end{theorem}

\subsubsection{Tur\'an problem}

Recall that the Tur\'{a}n problem was introduced in Example \ref{ex:TuranFunction}.
Given a family of graphs $\C H$ we
consider
$\C H^\uparrow$, the collection of graphs obtained by adding missing edges
to the graphs in $\C H$.
While our computer code can automatically prove the perfect stability when
$\C H=\{K_t\}$ with $t\le 7$,
 this is superceded by the following result whose proof does not require a computer. For
 integer $q\ge 1$ denote by $K_m^q$ the balanced blow
 up of $K_m$ on $qm$
 vertices (i.e.\ $\blow{K_m}{V_1,\dots,V_m}$ with $|V_1|=\dots=|V_m|=q$).

\begin{theorem}\label{th:turan}
Let $\C H$ be a family of graphs and let
 $$
 m:=\min\{\chi(H)\mid
H\in\C H\}-1\ge 2,
 $$ where by $\chi(H)$ we denote the chromatic number of the
graph
$H$, that is, the minimum number of colours needed in a coloring of the vertex
set with no adjacent vertices of the same colour.  Then
the following
hold.
\begin{enumerate}
 \item The Tur\'an problem $\ex(n,\C H)$ is perfectly $K_m$-stable if and only
 if there is an integer $q$  such that $K_m^q$ plus one edge is not
 $\C
 H^\uparrow$-free.
 \item Assuming in addition that $\C H$ is finite, we have the following. The
 Tur\'an problem $\ex(n,\C H)$ is robustly $K_m$-stable if and only if
 there is an integer $q$  such that $K_m^q$ plus some forest in one of the parts of $K_m^q$ is not
  $\C
  H^\uparrow$-free.
 \end{enumerate}
 \end{theorem}

As we learned later, the non-trivial implication in Part 1 of the above theorem is apparently a folklore result. Since it  follows from~\cite[Lemma~2.3]{roberts+scott}, we omit its proof. The second part of Theorem \ref{th:turan} is proved in Section~\ref{se:turan} of this paper.

\section{Notation and preliminaries}

Some of the definitions and proofs of this paper will be more natural when stated in a more analytic way. For example, 
the definition of $\lambda(\blow{B}{\V a})$ in~\eqref{eq:lambdaBA} would not require a limit if instead we were working with vertex-weighted graphs. Such objects are quite common in extremal graph theory nowadays (appearing, for example, in the definition of the \emph{Lagrangian} of a graph that goes back to Motzkin and Straus~\cite{MotzkinStraus65}) and, of course, they are generalised in a powerful and far-reaching way by graphons (see, for example, the excellent book by Lov\'asz~\cite{Lovasz:lngl}). However, 
we believe that by staying within the universe of simple unweighted graphs, we make the paper and its ideas better accessible to a wider audience.

As usual, for each positive integer $n$, we denote by $[n]$ the set
$\{1,...,n\}$. Let $\E(X)$ denote the expected value of a random variable~$X$.
We will often abbreviate a pair $\{i,j\}$ as $ij$. For a finite set $A$ and a positive integer $k$ we denote by ${A \choose k}$ the set of all $k$-subsets of $A$.

Recall that $K_m$ denotes the complete graph of order $m$
and $\overline G:=(V(G),{V(G)\choose 2}\setminus E(G))$ denotes the complement of a
graph $G$.
Let $K_{m,n}$ be the complete bipartite graph with part sizes $m$ and $n$.
For a vertex $x\in V(G)$, let $\Gamma_G(x)=\{y\in V(G)\mid \{x,y\}\in E(G)\}$
denote the \emph{neighbourhood} of $x$ in $G$. We write $G\cong H$ if $G$ and $H$ are isomorphic graphs.

We call a graph
$B$ \emph{$\lambda$-minimal} if $\lambda(\blow{B}{})$ strictly decreases when we
remove any vertex of $B$.
By~\eqref{eq:MaxLambdaBlow},
$B$ is $\lambda$-minimal if and only if no point on the boundary of $\I S_m$
achieves the maximum.

A graph $B$ is \emph{twin-free} if it contains no two vertices $x$ and $y$ with identical neighbourhoods (i.e., for all distinct $x,y \in V(B)$ we have $\Gamma_x(B) \neq \Gamma_y(B)$).

Recall that the \emph{(edit) distance} $\Dedit(G,H)$ between two
graphs $G$ and $H$ of the
same order $n$ is the minimum of adjacencies one has to edit in $G$ to make it isomorphic to $H$.
Also, the \emph{(edit) distance} $\Dedit(G,\C H)$ from
a graph $G$ to a family $\C H$ of graphs is the minimum of $\Dedit(G,H)$ over all
$H\in\C H$ that have the same order as~$G$; this is the minimum number of adjacency edits needed
to transform $G$ into a graph in $\C H$. The respective \emph{normalised distances} are
$\dedit(G,H):=\Dedit(G,H)/{n\choose 2}$ and $\dedit(G,\C H):=\Dedit(G,\C H)/{n\choose 2}$.

Throughout this paper we will work under the following assumptions which
are collected together for future reference.

\begin{assumption}\label{as:general}
 Let $\kappa, m$ be positive integers and $\C F$ a family of graphs.
 \begin{enumerate}
 \item\label{it:problem} Set $\C G=\Forb{\C F}$.
 \item
  Let $\gamma:\C G_\kappa^0\to\I R$ be a function and
  define $\Lambda$ and $\lambda$ as in \eqref{eq:lambdadefn}.
  \comment{We had $\C G\to\I R$ but $\Lambda,\lambda$ are defined only on graph of order $\ge \kappa$}
 \item Let $B$ be
a graph on $[m]$ such that $\blow{B}{}\subseteq \C G$.
\end{enumerate}
\end{assumption}

The next lemma provides some basic information on the behaviour of the sequence $(\lambda(n,\C G))_{n=\kappa}^\infty$.


\begin{lemma}\label{lm:lim} Let $\C G$ be a graph property closed
under taking induced subgraphs.
 Then for $\kappa\le q\le n$ with $q\to\infty$ we have
 $$
 0\le \lambda(q,\C G)- \lambda(n,\C G)\le o_q(1).
 $$
 Furthermore, if $\C G$ is closed under taking blow-ups, then the error term is $O(1/q)$.
 \end{lemma}

\bpf
Take an optimal graph $G$ for $\lambda(n,\C G)$. Let $X$ be a random
$q$-subset of $V(G)$ and $G':=G[X]$. Then $G'\in\C G$. Thus $\lambda(q,\C G)\ge \E(\lambda(G'))$. Clearly, if we take a uniform $X\in {V(G)\choose q}$ and then
a uniform $Y\in {X\choose \kappa}$, then $Y$ is uniformly distributed among
all $\kappa$-subsets of $V(G)$. Thus $\E(\lambda(G'))$ equals $\lambda(G)=\lambda(n,\C G)$, giving $\lambda(q,\C G)\ge \lambda(n,\C G)$. Thus $\lambda(q,\C G)$ is non-increasing in $q$ and
 tends to a limit, implying
the other desired inequality $\lambda(n,\C G)\ge \lambda(q,\C G)+o_q(1)$.

Finally, assume that $\C G$ is also closed under taking blow-ups. To show $\lambda(q,\C G)\le \lambda(n,\C G)\le O(1/n)$, take
an optimal graph $G$ for $\lambda(q,\C G)$ on $[q]$. Consider a random
map $\phi:[n]\to [q]$ with all $q^n$ choices being equally likely and let $G'$ be the graph on $[n]$ with
$E(G')=\phi^{-1}(E(G))$. Take any $\kappa$-subset $X\subseteq [n]$. With probability $1-O(1/q)$, the map $\phi$ is injective on $X$. If we condition on this, then $\phi(X)\in {[q]\choose \kappa}$ is uniform and the average of $\gamma(G'[X])=\gamma(G[\phi(X)])$  is $\lambda(G)$. 
Thus
 $$
 \lambda(n,\C G)\ge \E_\phi(\lambda(G'))=(1-O(1/q)) \lambda(G) -O(1/q),
 $$
 giving the desired.\epf

\section{Flag algebra method}\label{lower}\label{certificates}

As we have already mentioned in the introduction of this paper, the flag algebra
method is a powerful technique pioneered by
Razborov~\cite{razborov:07,razborov:10}. In this section, we define what a
certificate is and how it implies an upper bound on $\lambda(\C G)$.
Recall that we always work under Assumption~\ref{as:general}.

\subsection{Types and flags}

A \emph{type} is a pair of the form $(H,\phi)$, where $H$ is an admissible
graph and $\phi:[v]\to V(H)$ is a bijection with $v=v(H)$.
Given a type
$\tau=(H,\phi)$ as above,
a \emph{$\tau$-flag} is a pair of the form $(G,\psi)$, where $G$ is an
admissible
graph and $\psi:[v]\to V(G)$
is an injection such that $\psi\circ \phi^{-1}: V(H)\to V(G)$
is an \emph{embedding} (that is, an injection that
preserves both edges and non-edges). Informally, a
$\tau$-flag $(G,\psi)$ is
a partially labelled graph such that the labelled vertices induce~$\tau$. The
\emph{order} $v((G,\psi))$ of the flag is~$v(G)$, the number of vertices in it.

For two $\tau$-flags $(G_1,\psi_1)$ and $(G_2,\psi_2)$
with respectively $n_1\le n_2$ vertices, let the \emph{sub-flag count}
$\pp((G_1,\psi_1),(G_2,\psi_2))$
be the number of $n_1$-subsets $X$ of $V(G_2)$ such that
$X\supseteq \psi_2([v])$ (i.e.\ $X$ contains all labelled
vertices) and the $\tau$-flags
$(G_1,\psi_1)$ and $(G_2[X],\psi_2)$ are \emph{isomorphic},
meaning that there is a graph isomorphism that preserves
the labels. Also, define the \emph{(flag) density} as
 \begin{equation}\label{eq:PForFlags}
 p((G_1,\psi_1),(G_2,\psi_2)):=\frac{\pp((G_1,\psi_1),(G_2,\psi_2))}{{
n_2-v\choose n_1-v}}.
 \end{equation}
The latter quantity can be viewed as the probability that a random $n_1$-subset
$X$ of $V(G_2)$ with
$X\supseteq \phi_2([v])$ induces a copy of the flag $(G_1,\psi_1)$ in
$(G_2,\psi_2)$.

We will also need a variation of the above notion. Let $F_1=(G_1,\psi_1)$,
$F_2=(G_2,\psi_2)$ and $(G,\psi)$ be three $\tau$-flags with respectively
$n_1,n_2$ and $n$ vertices. We define the \emph{joint sub-flag count}
$\pp(F_1,F_2,(G,\phi))$ to be the number of pairs $(X,Y)$ such that $X,Y$ are
subsets of $V(G)$ with $n_1$ and $n_2$ elements respectively, $X\cap
Y=\psi([v])$ and the $\tau$-flags $(G[X],\psi)$ and $(G[Y],\psi)$ are
isomorphic to $F_1$ and $F_2$, respectively.

The type with no vertices will be denoted by 0. Thus $0$-flags are just
unlabelled
graphs. In this case, the 0-flag density as defined by~\eqref{eq:PForFlags}
coincides
with the notion of subgraph density from the Introduction.

\subsection{Certificates}

\begin{definition} A
\emph{(flag algebra) certificate} is a triple
 \begin{equation}\label{eq:Certificate}
  \C C=(N,\C T,(Q^{\tau})_{\tau\in \C T}),
 \end{equation}
  where
   \begin{itemize}
    \item $N\ge \kappa$ is an integer;
     \item $\C T=(\tau_1,\dots,\tau_t)$ is an
  ordered
list of
some types such that
$N-v(\tau_i)$ is a positive even integer for each $i\in [t]$;
 \item $Q^{\tau_i}$
 is
 an arbitrary positive
 semi-definite $g_i\times g_i$-matrix for $i\in [t]$, where we fix
 some enumeration $(F_1^{\tau_i},\dots,F_{g_i}^{\tau_i})$ of  all
 $\tau_i$-flags  with exactly $(N+v(\tau_i))/2$ vertices, up to
 isomorphism of flags (and thus
 $g_i$ is the number of these flags).
\end{itemize} \end{definition}

Note that the third component of the certificate $\C C$ consists of exactly $t$
matrices, one per each of
the types $\tau_1,\dots,\tau_t$; one can view the rows/columns of $Q^{\tau_i}$
as indexed by the $\tau_i$-flags of order $(N+v(\tau_i))/2$.

To describe the upper bound of $\lambda(\C G)$ that a certificate $\C C$
witnesses, we need to introduce several quantities.

Let $G_1,\dots,G_g$ be the enumeration in some fixed
order of all (up to an isomorphism) admissible $N$-vertex graphs. Thus $\C
G_N^0=\{G_1,\dots,G_g\}$ with no two listed graphs being isomorphic.
For each $q\in{}[g]$ (that is, for each $G_q$), we define real numbers
 \begin{equation}\label{eq:ai}
 a_q:=\sum_{i=1}^t\sum_{h=1}^{g_i} \sum_{j=1}^{g_i} c_{h,j,q}^{\tau_i} Q_{h,j}^{\tau_i}, \quad\text{where}\quad c_{j,h,q}^{\tau_i}:=\sum_\phi\pp(F^{\tau_i}_j,F^{\tau_i}_h,(G_q,\phi))
 \end{equation}
 and the sum in the definition of $c_{h,j,q}^{\tau_i}$ is taken over all
 injective maps $\phi:[v(\tau_i)]\to V(G_q)$ that induce a copy of the flag
 $\tau_i$ in $G$.
 Also, let
 $$
 b_q:=\lambda(G_q)=\sum_{H\in\C G_\kappa^0} \gamma(H)p(H,G_q),
  $$
  and
   $$u_\lambda(\C C):=\max\{a_q+b_q:q\in[g]\}.
   $$
    A graph $G_q\in\C G_N^0$ is called \emph{$(\C C,\lambda)$-sharp} (or \emph{$\C C$-sharp}, or  just \emph{sharp}) if $a_q+b_q=u_\lambda(\C C)$. The following lemma motivates the above  definitions.


\begin{lemma}\label{lm:U(C)} Under the above notation, for every admissible graph $G$ of order $n$ we have that
\begin{eqnarray}
    u_\lambda(\C C) -\lambda(G)+O(1/n)&=& \sum_{q=1}^g (u_\lambda(\C C)-b_q)
  p(G_q,G)+O(1/n)\label{eq:ai+biEq}\\
  &\ge& \sum_{q=1}^g a_q p(G_q,G)+O(1/n)
  \ \ge\ 0.\label{eq:ai+bi_new}
\end{eqnarray}
In particular, we have that
 $\lambda(n,\C G)\le u_\lambda(\C C)+O(1/n)$ and
  $\lambda(\C G)\le u_\lambda(\C C)$.
 \end{lemma}
\bpf 
Let $G$ be an arbitrary admissible graph of order $n$. We have
 \begin{equation}
   \label{eq:bi}
   \lambda(G)=\sum_{H\in\C G_\kappa^0}\gamma(H)p(H,G)=\sum_{H\in\C G_\kappa^0}\gamma(H)\sum_{q=1}^gp(H,G_q) p(G_q,G)=\sum_{q=1}^g b_qp(G_q,G),
 \end{equation}
 proving~\eqref{eq:ai+biEq}.

Next, we define a
non-negative quantity $a$ in the following way. Initially, we set $a=0$. For each
non-negative integer $v$ such that $N-v$ is a positive even integer we work as
follows. We enumerate all
$n(n-1)\dots (n-v+1)$ injections $\psi:[v]\to V(G)$. If the induced type
$(G[\psi([v])],\psi)$ is equal to some $\tau_i\in{}\C T$, then we add the quantity $\V x Q^{\tau_i}\V x^T$ to
$a$, where
 \beq{x}
 \V x:=\big(\pp(F_1^{\tau_i},(G,\psi)),\dots,\pp(F_{g_i}^{\tau_i},(G,\psi))\big).
 \eeq
 Since each matrix $Q^{\tau_i}$ is positive semi-definite, we have that each $\V x
Q^{\tau_i}\V x\ge 0$ and that the final $a$ is non-negative.

Let $i\in{}[t]$ and set $v=v(\tau_i)$. Take any $j,h\in[g_i]$. Notice that the sum of the products
$\pp(F_j^{\tau_i},(G,\psi))\, \pp(F_h^{\tau_i},(G,\psi))$ over
all
injections $\psi:[v]\to V(G)$ such that $(G[\psi([v])],\psi)$ is
isomorphic to $\tau_i$, is equal to
 $$
 \sum_{q=1}^g c_{j,h,q}^{\tau_i} \pp(G_q,G)+O(n^{N-1}),
 $$
 where $c_{j,h,q}^{\tau_i}$ is defined in~\eqref{eq:ai}, see
e.g.~\cite[Lemma~2.3]{razborov:07}. (Informally speaking, we just count in two different
ways the number of embeddings of $F_j^{\tau_i}$ and $F_h^{\tau_i}$ into $G$ so that
the correponding labelled vertices coincide; the error term $O(n^{N-1})$ comes
from embeddings where some unlabelled vertices happen to collide.)
Thus, summing over $i\in{}[t]$, as well as, over injections $\psi$ and
expanding
each quadratic form $\V x Q^{\tau}\V x^T$, we get the
representation
 \beq{a}
 0\le \frac{a}{{n\choose N}} =O(1/n)+\sum_{q=1}^g a_q p(G_q,G),
 \eeq
 where $a_q$'s are as in~\eqref{eq:ai}. Adding this to~\eqref{eq:bi}, we obtain that
 \begin{equation}\label{eq:ai+bi}
 \lambda(G)+O(1/n)\le   \sum_{q=1}^g (a_q+b_q)p(G_q,G)\le u_\lambda(\C C)\sum_{q=1}^g p(G_q,G)=u_\lambda(\C C),
 \end{equation}
 The inequalities in \eqref{eq:ai+bi_new} follows readily by~\eqref{eq:a} and~\eqref{eq:ai+bi}. \epf

Lemmas \ref{lm:lim} and \ref{lm:U(C)} have the following immediate consequence.

\begin{corollary}
  \label{corL:lambda(n,)}
  Suppose that Assumption \ref{as:general} holds.  Let $\V a$ be a vector in $\I S_m$ and $\C C$ a certificate
 such that $\lambda(\blow{B}{\V a})\ge u_\lambda(\C C)$.
 Then
 $\lambda(\C G)=\lambda(\blow{B}{\V a})=u_\lambda(\C C)$. Moreover, if $\C G$ is closed under taking blow-ups, then $\lambda(n,\C
 G)=u_\lambda(\C C)+O(1/n)$.\qed
\end{corollary}

Finally, we close this section with the following lemma.

\begin{lemma}\label{lm:forced} Under Assumption \ref{as:general}, suppose that $\C G$ is closed under taking blow-ups and that
a certificate $\C C$ and
a vector $\V a=(a_1,\dots,a_m)\in\I S_m$ with no zero entry satisfy
$u_\lambda(\C
C)=\lambda(\blow{B}{\V a})$.
Fix $i\in [t]$ and set $v:=v(\tau_i)$.
Also, let $n$ be a large positive integer, $V_1,\dots,V_m$ a partition of $[n]$
with $|V_i|=a_in+O(1)$
and $G:=\blow{B}{V_1,\dots,V_m}$. Finally, let
$\psi:[v]\to V(G)$ be an injection such that $(G,\psi)$ is a $\tau_i$-flag. Define $\V x$ as in~\req{x}. Then
$\V xQ^{\tau_i}\V x^T=O(n^{N-v-1})$.\end{lemma}

\bpf Observe that
any modification of the injection $\psi$ such that its values stay in the same
parts $V_i$ is an embedding. These new injections give the same
vector $\V x$. There are $\Omega(n^v)$ such injections since each part $V_i$ has
size $a_in+O(1)=\Omega(n)$. Let $a$ be the quantity defined in the proof of
Lemma \ref{lm:U(C)} when applied to $G$. Note that $a$ is the sum of
non-negative quantities, some of which correspond to the above $\Omega(n^v)$
embeddings of $\tau_i$ into $G$. Thus
\begin{equation}
  \label{eq:eigenv}
   0\le \Omega(n^v) \cdot \V xQ^{\tau_i}\V x^T\le a.
\end{equation}
Also, observe that $\lambda(G)=\lambda(\blow{B}{\V a})+O(1/n)$.
Since $u_\lambda(\C C)=\lambda(\blow{B}{\V a})$, by \eqref{eq:ai+bi_new} and
\eqref{eq:a} we have that $a=O(n^{N-1})$. Invoking \eqref{eq:eigenv}, the result follows.\qed

\hide{
By \eqref{eq:ai+bi}, we have the following.

\begin{lemma}\label{lm:NonSharp} If $G_q$ is not sharp, then for any
admissible $G$ of order $n\to\infty$ with $\lambda(G)=\lambda(\C G)$ we
have $p(G_q,G)=o(1)$.\end{lemma}}

\section{Robust stability from flag algebra proofs}

The main result of this section is Theorem \ref{th:stab1} below that provides a
sufficient condition for a problem to be robustly stable.
Let $H, G$ be two graphs. We say that a map $f:V(H)\to V(G)$ is a \emph{strong
homomorphism} if it preserves both adjacency and non-adjacency.
Observe that a strong homomorphism, in contrast to an embedding, does not need to be injective, allowing pairwise non-adjacent vertices
to be mapped to the same image. Moreover, let us note that a graph $H$ admits a strong homomorphism in a graph $B$ if and only if $H$ is a blow-up of $B$.

\hide{Let $$
 \I S_{m}^*=\{\V x \in \I R^{m}: x_1+\dots+x_m=1, \mbox{ each $x_i> 0$}\}
 $$ be obtained
from the standard $(m-1)$-simplex
$\I S_m$ by removing its $m$ vertices.
 }
 \hide{Recall that for $\V a\in \I S_m$, let $\lambda(\blow{B}{\V a})$ be the limit
of $\lambda(\blow{B}{V_1,\dots,V_m})$ as $n\to\infty$, where $V_1\cup\dots\cup V_m$ is a partition
of $[n]$ with $|V_i|/n\to a_i$. It is easy to see that $\lambda(\blow{B}{\V a})$ is a polynomial in $\V a$ of
degree at most $\kappa$.
}

\begin{theorem}[Robust Stability]\label{th:stab1}
Suppose that in addition to Assumption~\ref{as:general} the following holds.
\begin{enumerate}
 \item\label{it:FAProof1} We have a vector $\V a\in\I S_m$ and a certificate
 $\C C=(N,\C T,(Q^\tau)_{\tau\in \C T})$ with $u_\lambda(\C C) \le  \lambda(\blow{B}{\V a})$.
   \item\label{it:tau1} There is a graph $\tau$ of order at most $N-2$ satisfying the following.
   \begin{enumerate}
     \item[(a)] $\lambda(\Forb{\C F})>\lambda(\Forb{\C F\cup\{\tau\}})$.
     \item[(b)] There exists a unique
(up to automorphisms of $\tau$ and $B$) strong homomorphism $f$ from $\tau$ into $B$.
     \item[(c)] For every distinct
$x_1$ and $x_2$ in $V(B)$ we have $\Gamma_B(x_1)\cap f(V(\tau))\not= \Gamma_B(x_2)\cap f(V(\tau))$.
   \end{enumerate}
 \item\label{it:sharp1} Every $\C C$-sharp graph of order $N$ admits a
strong homomorphism
into $B$.
\end{enumerate}
 Then the problem is robustly $B$-stable.
\end{theorem}

\bpf  By Corollary~\ref{corL:lambda(n,)}, we know that
$\lambda(\C G)=\lambda(\blow{B}{\V a})=u_\lambda(\C C)$.
For notational
convenience, assume that
$V(\tau)=[q]$.  Choose large constants
in the
order $C_1\ll C$. In particular,
we assume that $C>2/(\lambda(\C G)-\lambda(\Forb{\C F\cup\{\tau\}}))$. Take any
$\C F$-free graph $G$ of order $n>C$. Note that Condition~\ref{it:tau1}(c) of
Theorem~\ref{th:stab1} implies
that $B$ is twin-free.

We can assume that $\lambda(G)\ge (\lambda(\C G)+\lambda(\Forb{\C F\cup\{\tau\}}))/2$ for otherwise
 $$
 C(\lambda(\C G)-\lambda(G))\ge C(\lambda(\C G)-\lambda(\Forb{\C F\cup\{\tau\}}))/2>1\ge \dedit(G,\blow{B}{}),
 $$
 and there is nothing to do. Since $G$ is $\C F$-free
but $\lambda(G)$ is strictly larger than $\lambda(\Forb{\C F\cup\{\tau\}})$, the
supersaturation argument of Erd\H os and Simonovits~\cite{erdos+simonovits:83} or an application of the Removal Lemma
shows that
 \begin{equation}\label{eq:tauG}
  p(\tau,G)\ge 1/C_1,
   \end{equation}
 that is, $G$ has at least ${n\choose q}/C_1$ copies of $\tau$.

For every embedding $\psi:\tau\to G$, we define the following. For each binary string $\V b=(b_1,...,b_q)$ of
length $q$, let $V_{\psi,\V b}$ consist of those vertices $x\in{}V(G)$ such that the neighbourhood
of $x$ in $\psi([q])$ is given by $\V b$, that is, $\{i\in [q]\mid
\{x,\psi(i)\}\in E(G)\}=\{i\in [q]\mid b_i=1\}$.
Thus, the sets $V_{\psi,\V b}$, $\V b\in \{0,1\}^{q}$, form a partition of $V(G)$.
Observe that, if we apply the above definition to the (fixed) map $f:\tau\to B$ (instead of $\psi:\tau\to G$), then
each part in the obtained partition of $V(B)=[m]$ has at most one vertex by Condition~\ref{it:tau1}.
Let $\V b^{j}\in\{0,1\}^q$ be
the binary sequence  corresponding to the part $\{j\}$ for each $j\in [m]$; thus $\V b^{j}$ encodes the adjacencies of $j\in V(B)$ to the fixed copy of $\tau$ in $B$.  We call all other length-$q$ binary sequences
\emph{singular}. 
 Also, we call a part $V_{\psi,\V b}$ \emph{singular} if $\V b$ is singular, that is, not one of $\V b^1,\dots,\V b^m$.
Finally, we call a pair of distinct vertices $x_1,x_2\in{}V(G)$
\emph{singular} if at least one of them is in a singular part or both of them
are in non-singular parts but the adjacency relations between $x_1,x_2$ in $G$
and between $j_1,j_2$ in $B$ mismatch (that is,
one is an edge and the other is a non-edge), where $j_l$ is the unique element
of $[m]$
satisfying $x_l\in V_{\psi,\V b^{j_l}}$ for $l=1,2$. Note if we have $j_1=j_2$ above, then $\{x_1,x_2\}$ is singular if and only if $x_1$ and $x_2$ are connected in  $G$.


Observe that due to Condition \ref{it:tau1}, we have that the union of $\psi([q])$ with every singular pair $\{x_1,x_2\}$ induces a graph that does not embed into a blow-up of $B$. For example, if $x_1$ is in a singular part then already $\psi([q])\cup\{x_1\}$ spans a subgraph in $G$ that does not belong to $\blow{B}{}$. If we add an arbitrary disjoint $(N-|X|)$-set
$Y$ of vertices to $X:=\psi([q])\cup\{x_1,x_2\}$, we get a subgraph of $G$ of order $N$ that does not
belong to $\blow{B}{}$. Condition~\ref{it:sharp1} and inequality \eqref{eq:ai+bi_new} give that the total number of such subgraphs in $G$ is at most $C_1{n\choose N}\max(1/n,\lambda(\C G)-\lambda(G))$, where we assume that $1/C_1$ is smaller than $\min\{u_\lambda(\C C)-a_q-b_q:G_q\;\text{is non-sharp}\}$.
Also,  each such subgraph $H$ of $G$ can arise for at most $N!$ triples $(\psi,\{x_1,x_2\},Y)$, a rough bound
on the number of ways to embed $\tau$ into $H$,  then choose two more vertices in $H$ and
let $Y$ be the rest of $V(H)$. Thus, the number of
triples $(\psi,\{x_1,x_2\},Y)$ as above is at most $C_1{n\choose N}\max(1/n,\lambda(\C G)-\lambda(G))\times N!$.
Clearly, if we fix the first two entries, namely $(\psi,\{x_1,x_2\})$, then any choice of $Y$ will do and there
are at least ${n\choose N-q-2}$ choices of $Y$ (as $|X|$ is always at most $q+2$).
Thus the total number of possible choices of $(\psi,\{x_1,x_2\})$ as above is at most
 $$
 C_1{n\choose N}\max(1/n,\lambda(\C G)-\lambda(G))\times N! / {n\choose N-q-2}.
 $$
 Choose $\psi$ for which the number of singular pairs is at most the average. By~\eqref{eq:tauG} and Corollary~\ref{corL:lambda(n,)}, it is at most
 $$
 \frac{C_1{n\choose N}\max(1/n,\lambda(\C G)-\lambda(G))\times N! / {n\choose N-q-2}}{{n\choose q}/C_1}< Cn^2\max(1/n,\lambda(n,\C G)-\lambda(G)).
 $$
   Observe that one can convert $G$ into a blow-up of $B$ by flipping all singular pairs between non-singular parts of $G$ and merging
the singular parts into non-singular ones in an arbitrary way.  Thus, for every (and in particular this) $\psi$,
  the number of singular pairs
  is at least $\Dedit(G,\blow{B}{})$, which is by definition the minimum number of pairs that one needs to change in $G$ to make is a blow-up of~$B$. This finishes the proof of the theorem.\epf


\section{Sufficient conditions for perfect stability}
The aim of this section is to present sufficient conditions for perfect stability.
To state our results, we need the notions of strictness and flip-aversion.
Their definitions require several other concepts that we introduce in the next section.

\subsection{Notation and some preliminary results}
Throughout this section we work under the following set of assumptions.
\begin{assumption}\label{as:exact}
In addition to Assumption~\ref{as:general}, we assume the following.
\begin{enumerate}
  \item\label{as:exact.i} Each graph in $\C F$ is twin-free and
  \item\label{as:exact.ii} $\lambda(\C G)=\lambda(\blow{B}{})$.
\end{enumerate}
\end{assumption}
Observe that a trivial consequence of twin-freeness of each $F\in{}\C F$ is the following.
\begin{lemma}\label{lm:ClosedByBlowups} The set of admissible graphs $\C G$ is closed
under taking blow-ups.\qed\end{lemma}

We will also need the following pieces of notation. If $G$ is a graph and $x,y$ is a pair of distinct nodes of $G$, then by
$G\oplus xy$ we denote the graph obtained by flipping the adjacency of $x$ and $y$, while by $G - x$ we denote the graph obtained by deleting the node $x$ in $G$. Moreover, if $\kappa$ is a positive integer, for a graph $G$ of order $n\ge \kappa$ and a vertex $x$ of $G$, we define
 \begin{equation}
   \label{eq:lambdaCondDef}
   \begin{split}
     \Lambda(G,x):=&\Lambda(G)-\Lambda(G-x)\;\text{and}\\
     \lambda(G,x):=& {n-1\choose \kappa-1}^{-1}\cdot\Lambda(G,x).
   \end{split}
 \end{equation}
 The value of  $\Lambda(G,x)$ can be determined by summing $\g{G[X]}$ over all $\kappa$-subsets $X$ of $V(G)$
containing $x$. Also, $\lambda(G,x)$ is the conditional
expectation of $\g{G[X]}$ where $X$ is a random $\kappa$-subset of $V(G)$ conditioned on $X\ni x$.


Let $\V a=(a_1,...,a_m)$ in $\I S_m$ be arbitrary.
Consider a blow-up $B':=\blow{B}{V_1,\dots,V_m}$ of order $n$, where $|V_i|=a_in+O(1)$.
Let $B''$ be obtained from it by adding a new vertex $w$. Then $\Lambda(B'',w)$
is determined within additive error $O(n^{\kappa-2})$ by the vector of ratios
 \begin{equation}\label{eq:DefY}
 \V y:=\left(\,\frac{|\Gamma_{B''}(w)\cap V_1|}{|V_1|},\dots,\frac{|\Gamma_{B''}(w)\cap V_m|}{|V_m|}\,\right)\in [0,1]^m.
 \end{equation}
 In fact, we have
 \beq{y}
 \lambda(B'',w)=R_{\V a}(\V y)+O(1/n),
 \eeq
 where $R_{\V a}=R_{B,\lambda,\V a}$ is some real polynomial
in~$\V y$. One
can write $R_{\V a}$ explicitly as follows.

First, for a (not necessarily injective) map $\phi:[t]\to [m]$
and a (binary) vector $\V b=(b_1,\dots,b_t)$ in $\{0,1\}^t$, let
 $B(\phi,\V b)$ be the graph on $[t+1]$ such that two elements $i$ and $j$ of $[t]$ are adjacent if and only if
$\phi(i)$ and $\phi(j)$ are adjacent in $B$, and $\{i,t+1\}$ is an edge
if and only if $b_{i}=1$. Informally speaking, $B(\phi,\V b)$ is a graph that
we can form from a blow-up of $B$ on $[t]$ by adding a new vertex whose neighbourhood
in $[t]$ is given by the binary vector $\B b$.%
\hide{Now, for every map $\phi:[t]\to [m]$ and vector $\V b=(b_1,\dots,b_t)$ in $\{0,1\}^t$, viewing $\V b$ as a map from $[t]$ to $\{0,1\}$, we denote by $\phi\otimes\V b$ the map from $[t]$ to $[m]\times\{0,1\}$ defined by $(\phi\otimes\V b)(i):=(\phi(i),b_i)$, $i\in [t]$.}
 Then the value of the polynomial $R _{\V a}$ at $\V y=(y_1,...,y_m)$ is
\[\sum_{\phi:[\kappa-1]\to[m]}\sum_{\V b\in\{0,1\}^{\kappa-1}}(\kappa-1)!\g{B(\phi,\V b)}\prod_{p=1}^m\prod_{q=0}^1
\frac{(a_p(qy_p+(1-q)(1-y_p)))^{|\{i:\phi(i)=p,\, b_i=q\}|}}{|\{i:\phi(i)=p,\, b_i=q\}|\,!}.\]


Let us call a vector $\V y\in{}[0,1]^m$ \emph{admissible} if for every
$t\in\I N$, every map $\phi:[t]\to [m]$, and every binary vector $\B
b=(b_1,\dots,b_t)\in\{0,1\}^t$ such that
$y_{\phi(i)}=0$ implies $b_i=0$ and $y_{\phi(i)}=1$ implies $b_i=1$
(while $b_i$
can be arbitrary if $0<y_{\phi(i)}<1$), the graph $B(\phi,\V b)$ is $\C F$-free.
In other words, this condition says that if we take a blow-up $\blow{B}{V_1,\dots,V_m}$
with each $|V_i|$ large and add a vertex $w$ with $y_i|V_i|$ neighbours in $V_i$ for each
$i\in{}[m]$, then the obtained graph is still $\C F$-free. Clearly, whether $\V y = (y_1,\ldots,y_m)$ is admissible or not, depends only on the sets $\{i\in[m]:y_i=0\}$ and $\{i\in[m]:y_i=1\}$ and therefore the next claim follows easily.

\begin{claim}
 \label{cl:adm}
  The set of the admissible vectors forms a closed subset of $[0,1]^m$.\qed
\end{claim}

Let us point out that, since $\C F$ is twin-free, it suffices to check the condition in the definition of an admissible $\V y$ only for those choices of $t,\phi,\V b$ for which $B(\phi,\V b)$ is twin-free. In particular, it suffices to consider $t$ to be at most $2m$.


The following vectors will play a special role. For every $i\in{}[m]$, we define
$\V v_i=(v_{i,1},...,v_{i,m})$ in $\{0,1\}^m$ by setting
$v_{i,j}=1$ if $ij$ belongs to $E(B)$ and $v_{i,j}=0$ otherwise for all $j\in{}[m]$. Informally speaking,
the assignment $\V y=\V v_i$ corresponds to adding one extra vertex in part $V_i$.
Thus each vector $\V v_i\in[0,1]^m$ is admissible.
Under this terminology, we have, in particular, for each $i_0\in{}[m]$ and $x\in{}V_{i_0}$ that
\begin{equation}\label{eq:008}
 R _{\V a}(\V v_{i_0})=\lambda(B',x)+O(1/n).
\end{equation}
Indeed, both sides of~\eqref{eq:008} measure the (normalised) change in the objective function $\lambda$ when we remove one vertex from the $i_0$-th part of an $(\V a+O(1/n))$-blow-up of $B$ of order~$n$.

There is the following connection between $\lambda(B',x)$ and $\frac{\partial}{\partial a_{i_0}} \lambda(\blow{B}{\V a})$.
\begin{claim}
  \label{cl:partial}
  \[R _{\V a}(\V v_{i_0})=\lambda(B',x)+O(1/n)=\frac{1}{\kappa}\frac{\partial}{\partial a_{i_0}} \lambda(\blow{B}{\V a}).\]
\end{claim}
\bpf
  First, for every positive integer $t$ and map $\phi:[t]\to[m]$ we define $B(\phi)$ to be the graph having $[t]$ as the vertex set with $i$ and $j$ being adjacent if and only if $\phi(i)$ and $\phi(j)$ are adjacent. For each $H\in{}\C G_\kappa^0$ we define
  \[\Phi_H:=\{\phi:[\kappa]\to[m]: B(\phi)\cong H\}\]
  and
  \[\Phi_H^{i_0}:=\{\phi\in\Phi_H:i_0\in\phi([\kappa])\}.\]
  Then we have that
  \begin{equation}
    \label{eq:009}
    \lambda(\blow{B}{\V a})=\sum_{H\in\C G_\kappa^0}\gamma(H)\kappa!\sum_{\phi\in\Phi_H}\prod_{i=1}^m\frac{a_i^{|\phi^{-1}(i)|}}{|\phi^{-1}(i)|!}.
  \end{equation}
  On the other hand, we have that
  \[\lambda(B',x)=\sum_{H\in\C G_\kappa^0}\gamma(H)(\kappa-1)!\sum_{\phi\in\Phi_H^{i_0}}
  \frac{a_{i_0}^{|\phi^{-1}(i_0)|-1}}{(|\phi^{-1}(i_0)|-1)!}
  \prod_{\substack{i=1\\i\neq i_0}}^m\frac{a_i^{|\phi^{-1}(i)|}}{|\phi^{-1}(i)|!}+O(1/n)\stackrel{\eqref{eq:009}}{=}\frac{1}{\kappa}\frac{\partial}{\partial a_{i_0}} \lambda(\blow{B}{\V a}).\epf\]

Let us illustrate some of the above concepts in the special case of
Example~\ref{ex:TuranFunction} with $\C H=\{K_t\}$ (namely, the Tur\'an function $\ex(n,K_t)$). Here $m=t-1$ and $B=K_{m}$.
Ignoring rounding errors, if we create $B'$ from
the complete $m$-partite graph
$K_{a_1n,\dots,a_mn}=\blow{B}{V_1,\dots,V_m}$ by adding a new vertex $w$ having
$y_1a_1n,\dots,y_ma_mn$ neighbours in $V_1,\dots,V_m$ respectively, then
$\Lambda(B',w)$ is just $\sum_{i=1}^m y_ia_in$, the number of edges at $w$. Thus
$R_{\V a}(\V y)=\lim_{n\to\infty} \Lambda(B',w)/n=\sum_{i=1}^m y_ia_i$. Since
we forbid $K_{m+1}$, a vector $\V y$ is admissible if and only if at least
one $y_i$ is 0. Here $\V v_i$ is the $\V y$-vector corresponding
to $w$ being a twin of the vertices in $V_i$, that is, $\V v_i$ consists of
1s except one 0 at position $i$. Note that $\V a=(1/m,\dots,1/m)$ is the (unique) maximiser
of $\lambda(\blow{B}{\V x})$ for $\V x\in \I S_m$. If we fix this $\V a$
and maximise $R_{\V a}(\V y)$ over admissible $\V y\in[0,1]^m$, then trivially
the set of maximisers is $\{\V v_1,\dots,\V v_m\}$. The following lemma states that one part of this inclusion (namely, that each $\V v_i$ is a maximiser) holds whenever $\V a$ has no zero entries. This makes a perfect combinatorial sense:
in every extremal configuration $\blow{B}{V_1,\dots,V_m}$ all vertices must
make asymptotically the same contribution to $\Lambda$.

\begin{lemma}\label{lm:Strict1} Fix any $\V a\in{}\I S_m$ that maximises $\lambda(\blow{B}{\cdot})$. Suppose that $\V a$ has no zero entries. Then the maximum of $R_{\V a}(\V y)$ over admissible $\V y\in [0,1]^m$
is $\lambda(\C G)$ and, furthermore, $R_{\V a}(\V v_i)=\lambda(\C G)$ for each
$i\in[m]$ (that is, each of the vectors $\V v_i$ is a maximiser).\end{lemma}
\bpf
Since $\V a$ achieves a maximum and lies in the interior of $\I S_m$, we have that $\frac{\partial}{\partial_i} \lambda(\blow{B}{\V a})=\frac{\partial}{\partial_j} \lambda(\blow{B}{\V a})$ for all $i,j\in[m]$. Denote this common value by $R$. By Claim \ref{cl:partial}, we have $R=\frac{1}{k}R _{\V a}(\V v_i)$ for all $i\in{}[m]$.
Since $\lambda(\blow{B}{\V a})$ is a homogeneous polynomial of degree $\kappa$, we have that
 $$\lambda(\blow{B}{\V a})=\frac{1}{\kappa}\sum_{i=1}^ma_i\frac{\partial}{\partial a_i} \lambda(\blow{B}{\V a}).
 $$
  This, the fact that $\V a$ maximises $\lambda(\blow{B}{\V a})$, Claim \ref{cl:partial}, and equality $\sum_{i=1}^ma_i=1$ imply that
\[\lambda(\C G)=\lambda(\blow{B}{\V a})
=\frac{1}{\kappa}\sum_{i=1}^ma_i\frac{\partial}{\partial a_i} \lambda(\blow{B}{\V a})=\sum_{i=1}^ma_i R _{\V a}(\V v_i)=R.\]
Thus $R _{\V a}(\V v_i)=\lambda(\C G)$ for all $i \in{}[m]$.

To prove the first part of the lemma, we derive a contradiction
by assuming that some admissible $\V y\in[0,1]^m$ achieves a strictly greater value. Let $c:=R_{\V a}(\V y)-\lambda(\C G)>0$ and pick some real $\e$ with $0<\e \ll c$.

Here we can start with
$B'=\blow{B}{V_1,\dots,V_m}$ of order $n\to\infty$ with $|V_i|/n\to a_i$
and form $B''$ by adding a set $Y$ of $\e n$ new vertices that span an independent set with
the identical adjacencies to $V_1,\dots,V_m$ governed by $\V y$.
Since the vector $\V y$ is admissible, the obtained
graph is $\C F$-free. Indeed, $\blow{B}{V_1,\dots,V_m}$ plus one vertex $v\in Y$ is
$\C F$-free by the admissibility of $\V y$; by blowing up the vertex $v$ we cannot violate $\C F$-freeness
because each member of $\C F$ is twin-free.

The
contribution of the new vertices to $\lambda$ is $c\e -O(\e^2)$. Indeed,
if we take a random $\kappa$-subset $X$ of $V(B'')$, then with probabilities respectively $1-\kappa \e+O(\e^2)$, $\kappa\e+O(\e^2)$, and $O(\e^2)$, the set $X$ intersects $Y$ in zero, one and at least two vertices;
thus
 $$
 \lambda(B'')=(1-\kappa \e)\lambda(B')+\kappa\e (\lambda(B')+c)+O(\e^2)= \lambda(B')+c\kappa\e+O(\e^2).
$$
 So we see that $\lambda(B'')-\lambda(B')$ can be made strictly positive by choosing small constant $\e\ll c$.
Thus the $(\V a+o(1))$-blow-up $B'$ of $B$ is not asymptotically optimal, contradicting
the optimality of $\V a$.\epf


Let
us say that $B$ is
\emph{$(\lambda,\V a)$-strict} if the set of maximisers of $R_{\V a}(\V y)$ over the admissible $\V y$'s in $[0,1]^m$
is exactly $\{\V v_1,\dots,\V v_m\}$. Call the graph $B$ \emph{$\lambda$-strict}
if $B$ is $\lambda$-minimal and $B$ is $(\lambda,\V a)$-strict for every $\V a\in\I S_m$ that maximises $\lambda(\blow{B}{\V a})$.

Recall that $B$ is \emph{$\lambda$-minimal} if $\lambda(\blow{B'}{})$ is strictly smaller
than $\lambda(\blow{B}{})$ for any proper subgraph $B'$ of $B$.
It trivially follows that such  $B$ is necessarily twin-free and every
maximiser $\V a\in{}\I S_m$ has all coordinates non-zero (and by compactness
at least one maximiser $\V a$ exists). Thus, if $B$ is $\lambda$-strict, then,
for each optimal $\V a$,  each of
$\V v_1,\dots,\V v_m\in[0,1]^m$ is a maximiser of $R_{\V a}$
by Lemma~\ref{lm:Strict1} while
the strictness property requires that there are no other maximisers.

\begin{lemma}\label{lm:Strict2} If $B$ is $\lambda$-strict, then for every $\e>0$ there is $\delta>0$ such that if $\V a\in\I S_m$ and admissible $\V y\in[0,1]^m$ satisfy $\lambda(\blow{B}{\V a})\ge \lambda(\C G)-\delta$ and $R_{\V a}(\V y)\ge \lambda(\C G)-\delta$, then
$\V y$ is $\e$-close to some $\V v_i$.
\end{lemma}
\bpf Suppose there is $\e>0$ that violates the lemma, that is, for every $j\in\I N$ there are $\V a^j$ and $\V y^j$ that violate the conclusion for $\delta=1/j$. By passing to a subsequence, we may assume that these vectors converge
to $\V a$ and $\V y$ respectively. By the continuity of $\lambda(\blow{B}{\V x})$, $\V a$ is a maximiser.
By $\lambda$-minimality of $B$ we have that each $a_i>0$.
By Claim \ref{cl:adm}, we have that $\V y$ is admissible, while by Lemma \ref{lm:Strict1} and Assumption \ref{as:exact} we have that $R_{\V a}(\V y)=\lambda(\C G)$. Since $B$ is $\lambda$-strict, we have that $\V y$ is equal to $\V v_i$ for some $i\in [m]$.
But then $\V y^j$ has to get $\e$-close to $\V v_i$
leading to a contradiction.\epf

Here is another easy consequence of the compactness of $\I S_m$.

\begin{lemma}\label{lm:boundary} If $B$ is $\lambda$-minimal, then there is $\delta>0$ such that for every $\V a\in{}\I S_m$ satisfying $\lambda(\blow{B}{\V a})\ge \lambda(\C G)-\delta$ we have that
each $a_i$ is at least $\delta$.\qed
\end{lemma}

Finally, we call a graph $B$ \emph{$\lambda$-flip-averse} if there is $\delta>0$ such that the following holds.
If we take a blow-up $B'=\blow{B}{V_1,\dots,V_m}$ with $n\ge 1/\delta$ vertices
such that $\lambda(B')\ge \lambda(\blow{B}{})-\delta$
and obtain $B'\oplus xy$ by changing the adjacency between a pair of
distinct nodes $x,y\in{}V(B')$ (possibly from the same part), then either $B'\oplus xy$ contains some $H\in{}\C F$ with
$v(H)\le m+2$  as a subgraph or we have that
 \begin{equation}\label{eq:flip}
 \Lambda(B')-\Lambda(B'\oplus xy)\ge \delta n^{\kappa -2}.
 \end{equation}
 By compactness, the property of being flip-averse can be equivalently re-formulated in terms of the polynomial $\lambda(\blow{B}{\V x})$, where $\delta$ disappears from the definition but then its combinatorial meaning will
be less clear.

\subsection{Main results for perfect stability.}

This section consists of two results, each providing a sufficient condition for perfect stability.
The first one is the following.
\begin{theorem}[Perfect Stability I]\label{th:exact}
Suppose that, in addition to Assumptions~\ref{as:general} and~\ref{as:exact}, the following assumptions hold.
 \begin{enumerate}
 \item\label{it:clstable} The problem is classically $B$-stable.
 \item\label{it:strict} The graph $B$ is $\lambda$-strict.
 \item\label{it:flip} The graph $B$ is $\lambda$-flip-averse.
 \end{enumerate}
Then the problem is perfectly $B$-stable.
 \end{theorem}

\hide{\newcommand{\clambda}{c_{\lambda}}
\newcommand{\ca}{c_{a}}
\newcommand{\cw}{c_{w}}
\newcommand{\cspecial}{c_{s}}
\newcommand{\czero}{\cw}
\newcommand{\ctwo}{c_2}
\newcommand{\cthree}{c_3}
}
\newcommand{\clambda}{c_1}
\newcommand{\cw}{c_2}
\newcommand{\czero}{\cw}
\newcommand{\ctwo}{c_3}
\newcommand{\cspecial}{c_4}
\newcommand{\cthree}{c_5}
\newcommand{\ca}{c_6}
\bpf Given $\lambda,B,\C F$, we fix sufficiently small positive constants $\ca{}\gg \cthree{}\gg \cspecial{}\gg \ctwo{}\gg 
\cw{}\gg\clambda$.
\hide{
First, observe that by the $\lambda$-minimality of $B$ and the compactness of $\I S_m$, we can assume
that for every $\V x\in\I S_m$ with $\lambda(\blow{B}{\V x})>\lambda(\C G)-\clambda$
we have that $\V x$ is at least $\ca$-far in $L_1$-distance from the boundary of $\I S_m$ and $\V x$ is $\ctwo{}$-close to a
$\lambda$-optimal $\V a\in \I S_m$.
Lemma~\ref{lm:Strict2}  implies that if $\lambda(B')\ge \lambda(\C G)-\clambda$
for some blow-up $B'$ of $B$ of sufficiently large order $n$ and we add a new vertex $w$ to $B'$ to obtain
an $\C F$-free graph $B''$ with $\lambda(w,B'')\ge \lambda(\C G)-\clambda$, then by editing
at most $\ctwo{}n$ edges at $w$ we can make $B''$ a blow-up of $B$.
}
To prove the perfect stability we pick some large enough real number $C$ (depending on the previous constants).
In this proof, let asymptotic notation such as $O(1)$ or
$\Omega(1)$ hide constants that depend on $\C F$, $\kappa$, $\g{\cdot}$, and $B$ only (but not on
the constants $c_i$).

\renewcommand{\marginpar}[1]{}
Let $n$ be an integer with $n>C$.
Choosing $C$ large enough, we may assume that $\lambda(\C G)+\clambda/2\ge \lambda(n,\C G)$.
Let $G$ be an arbitrary admissible graph on $[n]$. Assume that $\lambda(G)\ge \lambda(n,\C G)-\clambda/2$
for otherwise the result follows trivially, since $C\clambda/2>1$ and the normalised distance $\dedit$ is always bounded by $1$.\marginpar{$\clambda> 1/C$}
By Condition \ref{it:clstable}, that is, the classical $B$-stability,
there is a partition $[n]=V_1\cup\dots\cup V_m$ such that\marginpar{$\cw\gg \clambda$}
 \beq{c0}
 |W|\le \cw {n\choose 2},
 \eeq
 where $W:=E(G)\bigtriangleup E(B')$ and $B':=\blow{B}{V_1,\dots,V_m}$.
We call pairs in $W$ \emph{wrong}. Assume that the parts $V_1,\dots,V_m$
were chosen so that $|W|$ is minimum. Clearly, this choice of parts $V_i$ implies
that \eqref{eq:c0} still holds. Since the number of $\kappa$-subsets $X$
of $[n]$ such that $G[X]\not\cong B'[X]$ is at most $|W|{n\choose \kappa-2}$,
we conclude that
 \begin{equation}\label{eq:LargeL}
 \Lambda(B')\ge \Lambda(G)-|W|{n\choose \kappa-2}\cdot 2\gmax
  \ge (\lambda(\C G)-O(\cw{})){n\choose \kappa}.
 \end{equation}
 where $\gmax:=\max\{|\g{H}|\mid H\in\C G_\kappa^0\}$.

Let us call a vertex $x$ \emph{special} if $\lambda(G,x)< \lambda(\C G)- \cspecial{}$.
We set
$S$ to be the set of special vertices and $\sigma:=|S|/n$.

For each $i\in{}[m]$, we set $b_i:=|V_i|/n$. By~\eqref{eq:LargeL}, the continuity of $\lambda(\blow{B}{\cdot})$, and the compactness of $\I S_m$, we can assume that the vector $\V b=(b_1,\dots,b_m)$ is $\ctwo{}$-close
to a maximiser $\V a$ of $\lambda(B\blow{}\cdot)$, that is,\marginpar{$\ctwo{}\gg\cw$}
\begin{equation}\label{eq:ab}
\|\V a-\V b\|_1\le \ctwo{}.
\end{equation}
By Lemma~\ref{lm:boundary},
we can assume that each $a_i\ge \ca{}$; thus we conclude that
$b_i\ge \ca{}-\ctwo{}\ge \ca{}/2$ for each $i\in{}[m]$.\marginpar{$\ca\ge \ctwo{}$}

At this point, we can give an informal overview of the rest of the proof. 
First, Claim~\ref{cl:MaxLambda(x)} shows
that, for every vertex $x$ of $G$,  the normalised contribution $\lambda(G,x)$ of a vertex $x$ to $\lambda(G)$ is less than $\lambda(\C G)+\czero{}$ for otherwise the addition of an appropriate number of clones of $x$ to $G$ will bring $\lambda(G)$ well over $\lambda(\C G)$, which is impossible. 
It follows that, in order to avoid $\lambda(G)$ being too small, we have that $\sigma=O(\czero{}/\cspecial{})$. Furthermore, the adjacenty of each vertex $x\in [n]\setminus S$ essentially follows the ideal adjacency of
part-$i$ vertices, for some $i\in [m]$, as this is the only possibility to have $\lambda(G,x)$ close to  $\lambda(\C G)$ by
the assumed $\lambda$-strictness. Since our choice
of the parts $V_i$ minimises the number of wrong adjacencies, this vertex $x$ has to belong to $V_i$ and thus
its wrong degree $|\Gamma_W(x)|$ is necessarily small, see~\eqref{eq:DeltaW}. (Also, somewhat conversely, each vertex $x\in S$ has high wrong degree, just to account for the drop $\lambda(x,G)< \lambda(\C G)-\cspecial{}$.)
This, the near-optimality of $\V b$, the fact that $|S|=\sigma n$ is small and
the $\lambda$-flip-aversion give that every edge-flip inside $[n]\setminus S$ has negative effect on $\lambda$ (Claim~\ref{cl:Gradient}), not only with respect to $G$ but also with respect an arbitrary graph $\widetilde{G}$ obtained from $G$ by changing some adjacencies inside $W$ (Claim~\ref{cl:tildeG}). Thus if we flip $W'$, all wrong pairs outside $S$, and ``fix'' each vertex of $S$, then  $\Lambda$ increases by at least $\Omega(\cspecial{} n^{\kappa-2})$ per one changed edge. (Note that, since all vertices of high $W$-degree are inside the small set $S$, the ``pairwise'' effects can be shown to be negligible.) On the other hand, $|W'|+n|S|$ is clearly an upper bound on the edit distance from $G$ to 
the family $\blow{B}{}$. These two estimates
give the perfect stability. 

Let us provide all the remaining details now.

\Claim{MaxLambda(x)}{For every vertex $x$, we have that $\lambda(G,x)< \lambda(\C G)+ \czero{}$.}

\bcpf We assume on the contrary that there exists a node $x_0$ satisfying $\lambda(G,x_0)\geqslant \lambda(\C G)+ \czero{}$. Set $\e=\czero{}^2$.
Consider $G'$ obtained from $G$ by adding $\e n$ clones of $x_0$. We view $\lambda(G')$ as the expectation of $\g{G'[X]}$
for a random $\kappa$-set $X$. With probability at least $1-\kappa \e$,
the set $X$ is disjoint
from the added clones and its conditional expectation is exactly $\lambda(G)$.
With probability $\kappa\e+O(\e^2)$, the set $X$ has exactly one element from
the added clones and avoids $x_0$. Conditioned on the latter event, $G'[X]$ is the
same as $G[Y]$ where we take a random $\kappa$-subset $Y$ of $V(G)$ conditioned
on $Y\ni x_0$; thus the conditional expectation of $\g{G'[X]}$ is exactly
$\lambda(G,x)$ (which we assumed to be at least $\lambda(\C G)+\czero{}$). Finally, the contribution from the
remaining sets is in the absolute value at most $2\gmax$ times their probability $O(\e^2)$. Also, note our choice of $G$ such that $\lambda(G)\ge \lambda(n,\C G)-\clambda/2\ge \lambda(\C G)-\clambda$. Thus
 \begin{eqnarray*}
  \lambda(G') &\ge&(1-\kappa\e) \lambda(G)+\kappa\e \lambda(G,x)
   +O(\e^2)\\
  &\ge& (1-\kappa\e)(\lambda(\C G)-\clambda/2)+\kappa\e(\lambda(\C G)+\czero{})
 -O(\e^2)\\
  & \ge& \lambda(\C G)+\kappa \czero{}^3-\clambda-O(\czero{}^4).
  \end{eqnarray*}
 This is strictly larger than $\lambda(\C G)$.\marginpar{$\czero{}\gg \clambda$} On the other hand, by Lemma~\ref{lm:ClosedByBlowups}, we have that
$G'$ is admissible and therefore, invoking Lemma \ref{lm:lim}, we get that $\lambda(G')\le \lambda(\C G)+O(1/n)\le \lambda(\C G)+O(1/C)$,
a contradiction.\ecpf

If we pick a uniform random $x\in{}[n]$, then the difference $\lambda(\C G)-\lambda(G,x)$ is never below
$-\czero{}$ by Claim~\ref{cl:MaxLambda(x)}, while with probability $\sigma$ it is at least $\cspecial{}$. On the other hand, the average of $\lambda(\C G)-\lambda(G,x)$ over
$x\in V(G)$ is $\lambda(\C G)-\lambda(G)\le \clambda$.
Thus $-(1-\sigma)\czero{}+\sigma \cspecial\le \clambda$\marginpar{$\czero{}\ge \clambda$} and, roughly,
$\sigma\le 2\czero{}/\cspecial$.

Take any $x\in{}[n]$. Let $B_x'$ be obtained from $B'$ by changing adjacencies at $x$ so that
$\Gamma_{B_x'}(x)=\Gamma_G(x)$. We have that $\lambda(B_x',x)=R_{\V b}(\V y_x)+O(1/n)$, where $\V y_x=(y_{x,1},...,y_{x,m})$ is an element of $[0,1]^m$ defined by $y_{x,i}:=|\Gamma_G(x)\cap V_i|/|V_i|$ for all $i\in{}[m]$.
We also define another element $\V y'_x=(y'_{x,1},...,y'_{x,m})$ of $[0,1]^m$ by setting $y_{x,i}':=y_{x,i}$ unless if $y_{x,i}\le \ctwo{}/m$ (resp.\ $y_{x,i}\ge 1-\ctwo{}/m$), then we set $y_{x,i}':=0$ (resp.\ $y_{x,i}':=1)$. Clearly,
 \begin{equation}\label{eq:yy'}
  \|\V y_x-\V y_x'\|_1\le \ctwo{}.
  \end{equation}

\Claim{YAdmissible}{The vector $\V y_x'$ is admissible.}

\bcpf Suppose that the claim does not hold. Let this be
witnessed by a vector $\V b\in\{0,1\}^v$ and a map $\phi:[v]\to [m]$. Then $y_{x,\phi(i)}'\in\{0,1\}$ implies $b_i=y_{x,\phi(i)}'$, while the graph $B(\phi,\V b)$ is
of order $v+1$ and not $\C F$-free. As we observed after the definition of
an admissible vector, one can assume that $v\le 2m$. If $y_{x,i}'$ does not belong to $\{0,1\}$, then $y_{x,i}$ is
$\ctwo{}/m$-far from $0$ and $1$. Also, we know that each $V_i$ has at least $\ca{}n/2$ vertices.

Let us show that $B'_x$ has at least $\Omega((\ctwo{}\ca{}n)^{v})$ copies of $B(\phi,\V b)$ via $x$.
In fact, it is enough to consider only the copies where the vertex $v+1$ of $B(\phi,\V b)$ is mapped into $x$.
For $i\in [v]$, let $T_i$ be $V_i\setminus \Gamma_G(x)$ if $b_i=0$ and  $V_i\cap \Gamma_G(x)$ if $b_i=1$; note that $T_i$ always has at least $|V_i|\times \ctwo{}/m-O(1)$ vertices. Now, if we map each $i\in [v]$ arbitrarily into $T_i$, then these vertices together with $x$ form a copy of $B(\phi,\V b)$ in $B'_x$, giving at least the stated number of copies.

Each of the above copies contains a wrong
pair which is not adjacent to $x$. (Recall that $G$ is $\C F$-free but $B(\phi,\V b)$ is not and that the vertex $x$ has the same neighbourhoods in $G$ and $B_x'$.)  On the other
hand, each wrong pair disjoint from $x$ can be counted at most $n^{v-2}$ times.
This gives at least $\Omega((\ctwo{}\ca{}n)^{v})/n^{v-2}$ wrong pairs, contradicting~\eqref{eq:c0} since $\cw$ is sufficiently small with respect to $\ctwo{}$ and $\ca$
(and $v\le 2m$).\marginpar{$\ca\gg \cw$}\ecpf\\
By~\eqref{eq:ab} we also have that $|R_{\V b}(\V y_x)-R_{\V a}(\V y_x)|\le O(\ctwo{})$. 
On the other hand, there are at most $|W|{n-3\choose \kappa-3}$
$\kappa$-subsets $X$ of $[n]$ satisfying $x\in X$ and $G[X]\not\cong B_x'[X]$, because each such set must contain
a wrong pair disjoint from $x$. Thus by~\eqref{eq:c0}, we have that
$|\lambda(G,x)-\lambda(B_x',x)|\le O(\cw)$. Also, observe that
$\lambda(B_x',x)=R_{\V b}(\V y_x)+O(1/n)$. By~\eqref{eq:yy'}, we get that $|R_{\V a}(\V y_x)-R_{\V a}(\V y_x')|\le O(\ctwo{})$. By the Triangle Inequality, we derive
that
 $$
 |\lambda(G,x)-R_{\V a}(\V y_x')|\le O(\ctwo{}).
 $$

Suppose furthermore that $x\in{}[n]\setminus S$. By the definition of $S$, we have that $\lambda(G,x)\ge \lambda(\C G)-\cspecial$ and therefore
  \begin{equation}\label{eq:Ray'}
   R_{\V a}(\V y_x')\ge  \lambda(G,x)-O(\ctwo{})\ge \lambda(\C G)-O(\cspecial).
  \end{equation}
By Assumption~\ref{it:strict}, Lemma~\ref{lm:Strict2}
and Inequality~\eqref{eq:Ray'}, we conclude that $\V y_x'$ is $(\cthree{}/2)$-close
(in the $L_1$-norm) to the ``adjacency vector'' $\V v_i$ of some $i\in{}[m]$.\marginpar{$\cthree{}\gg\cspecial$}
By~\eqref{eq:yy'},\marginpar{$\cthree{}\ge 2\czero{}$}
 \begin{equation}\label{eq:yvi}
  \|\V y_x-\V v_i\|_1\le \cthree{}.
   \end{equation}

Next, let us show that $x$ belongs to $V_i$.
Suppose on the contrary that $x$ belongs to $V_j$ for some $j\not=i$. By the twin-freeness of $B$
(which trivially follows from the $\lambda$-minimality of $B$),
there is some
$h\in{}[m]$ which is adjacent to exactly one of $i$ and $j$, say $ih\in E(B)$
but $jh\not\in E(B)$. The vertex $x$ is adjacent in $G$ to $y_{x,h}|V_h|$
vertices of $V_h$. But we know that $|V_h|\ge \ca n/2$ and, by~\eqref{eq:yvi}, $y_{x,h}\ge v_{i,h}-\cthree{}=1-\cthree{}$. On the other hand, $B'$ has no edges between $V_i\ni x$
and~$V_h$. Thus $x$ belongs to at least $(1-\cthree{})\ca n/2$ wrong pairs having an endpoint in~$V_h$. Let as denote this set of edges by $A$. Consider
changing the partition $V_1\cup\dots\cup V_m$ by
moving $x$ to~$V_i$. Observe that the new set of wrong pairs will differ from the old one only on edges containing~$x$.
By~\eqref{eq:yvi}, at most $\cthree{}n$ edges can be introduced into the set of wrong pairs, while every edge in $A$ will not be contained, anymore, in the new set of wrong pairs. Thus the number of wrong pairs will strictly decrease.
This contradicts the choice of the partition $V_1\cup\dots\cup V_m$ and, in particular, the minimality of~$|W|$.
Thus indeed $x\in V_i$, as claimed.

Thus, again by~\eqref{eq:yvi}, we have that
 \beq{DeltaW}
 |\Gamma_{W}(x)|\le \cthree{}n,\quad \forall\, x\in [n]\setminus S.
 \eeq

\Claim{Gradient}{For every pair $xy$ in $W\cap {[n]\setminus S\choose 2}$, the graph $B'\oplus xy$ (which is obtained from $B'$ by changing the adjacency of $xy$) is $\C F$-free and satisfies
$\Lambda(B') - \Lambda(B'\oplus xy)
\ge \ca n^{\kappa-2}$.}

\bcpf Suppose on the contrary that $B'\oplus xy$
contains a forbidden subgraph $H\in\C F$.
Since $\C F$ consists of twin-free graphs
and $B'\oplus xy$ has at most $m+2$ pairwise non-twin vertices,
we can assume that $H$ has $v\le m+2$ vertices.
In fact, we must have
at least ${(\ca{}/2-\sigma)n\choose v-2}$ copies of $H$ on $[n]\setminus S$ via $xy$ in $B'\oplus xy$, since $|V_i|\ge \ca{}n/2$ for
each $i\in{}[m]$. Notice that the vertex set of each such copy
contains a pair from $W$ different from $xy$. By~\eqref{eq:DeltaW}, we have at most $2\cthree{}n$ wrong pairs
adjacent to $xy$, each in at most $n^{v-3}$ copies of $H$; while every other wrong pair appears in at most $n^{v-4}$ copies of $H$. This gives that the total number of $H$-subgraphs on $[n]\setminus S$ via $xy$ is at most\marginpar{$\cthree{}\ge \cw$}
 $$
 2\cthree{}n \cdot n^{v-3}+ |W|\cdot n^{v-4} \le 4\cthree{}n^{v-2},
 $$
 where we used~\eqref{eq:c0}. This is strictly less than ${\ca{}n/2\choose v-2}$, a contradiction.\marginpar{$\ca\ge \cthree{}$}
This contradiction shows that no such $H$ exists, proving the first part of the claim.

 The second part follow from Assumption~\ref{it:flip} of the theorem.\ecpf

\begin{claim}\label{cl:tildeG}
  Let $\widetilde{G}$ be an arbitrary (not necessarily $\C F$-free) graph having $[n]$ as a vertex set and such that $\widetilde{W}\subseteq W$, where $\widetilde{W}:=E(\widetilde{G})\bigtriangleup E(B')$. Then for every pair $xy$ in $\widetilde{W}\cap{[n]\setminus S\choose2}$ we have
  \beq{incr}
  \Sigma':=\sum_{X\in {[n]\setminus S\choose \kappa}}\left(\gamma((\widetilde{G}\oplus xy)[X])-\gamma(\widetilde{G}[X])\right)>\ca n^{\kappa-2}/2.\eeq
\end{claim}
\bcpf
Let us estimate $\Sigma'-\Sigma''$, where we define
 \begin{eqnarray*}
 \Sigma''&:=& \sum_{X\in {[n]\setminus S\choose \kappa}}\left(\gamma(B'[X]) - \gamma((B'\oplus xy)[X])\right).
 \end{eqnarray*}
 Let $X$ be a $\kappa$-subset of $[n]\setminus S$ that contributes different amounts to $\Sigma'$ and $\Sigma''$. Clearly, both $x$ and $y$ belong to $X$; also $X$ has to
contain at least one further pair $ab\in{}\widetilde{W}$. The number of the $\kappa$-subsets $X$ containing a pair $ab\in{}\widetilde{W}$ satisfying $\{a,b\}\cap \{x,y\}=\emptyset$ is at most $|\widetilde{W}|\leq|W|$ (the number of choices of $ab$) times ${n-4\choose \kappa-4}$ (the number of choices of $X\setminus\{a,b,x,y\}$).
Likewise, the number of the $\kappa$-subsets $X$ containing a pair $ab\in{}\widetilde{W}$ satisfying $\{a,b\}\cap \{x,y\}\neq\emptyset$ is at most the number of wrong pairs
adjacent to $x$ or $y$, which by \eqref{eq:DeltaW} satisfies $$
 |\Gamma_{\widetilde{W}}(x)|+|\Gamma_{\widetilde{W}}(y)|\leqslant|\Gamma_W(x)|+|\Gamma_W(y)|
 \le 2\cthree n,
 $$
 times ${n-3\choose \kappa-3}$.  Thus,~\eqref{eq:c0} gives that $|\Sigma'-\Sigma''|\le O(\cthree{} n^{\kappa-2})$. On the other hand, the sum
 $$
 \Sigma''':=\sum_{X\in {[n]\choose \kappa}\setminus {[n]\setminus S\choose \kappa}}\left(\gamma(B'[X]) - \gamma((B'\oplus xy)[X])\right)
 $$
 has at most $|S|n^{\kappa-3}$ non-zero terms (all such $X$ have to contain the pair $xy$ as well as intersect~$S$). Observe that $\Sigma''+\Sigma'''=\Lambda(B') - \Lambda(B'\oplus xy)$ is at least $c_6n^{\kappa-2}$ by Claim~\ref{cl:Gradient}.
 Thus $\Sigma'\ge c_6n^{\kappa-2}/2$, as desired. \ecpf

Enumerate $W':=W\cap {[n]\setminus S\choose 2}$ as $\{e_1,\dots,e_w\}$. Let $G_0:=G$ and
for $i=1,\dots,w$, let $G_i=G_{i-1}\oplus e_i$; that is, we flip the wrong pairs on $[n]\setminus S$ in some order. The final graph $G_w$ coincides with $B'$ on $[n]\setminus S$. By using Claim~\ref{cl:tildeG} to estimate the effect of each of the $w$ flips, we conclude that
\hide{
Thus, for every $\widetilde{G}$ as in the claim above, we have that
\[  \sum_{X\in{[n]\setminus S\choose\kappa}}\gamma((\widetilde{G}\oplus xy)[X])-\gamma(\widetilde{G}[X])>(\ca/2-O(\sigma))n^{\kappa-2}.
\]
}%
\begin{equation}
  \label{eq:020}
   \sum_{X\in{[n]\setminus S\choose\kappa}}\left(\gamma(B'[X])-\gamma(G[X])\right)\geq w\ca n^{\kappa-2}/2.
\end{equation}
On the other hand, we have that
\begin{equation}
  \label{eq:021}
  \sum_{\substack{X\in{[n]\choose\kappa}\\X\cap S\neq\emptyset}}\left(\gamma(B'[X])-\gamma(G[X])\right)\geq \sum_{x\in S}\left(\Lambda(B',x)-\Lambda(G,x)\right)-O(|S|^2n^{\kappa-2}).
\end{equation}
For each vertex $x\in{}S$, the value $\lambda(G,x)$ is at most
$\lambda(\C G)-\cspecial$ by the definition of $S$. By Claim~\ref{cl:partial}, the value $\lambda(B',x)$
is equal to $\frac{1}{\kappa}\frac{\partial}{\partial_i}\lambda(\blow{B}{\V b})+O(1/n)$ where $i\in[m]$ is the
index of the part $V_i$ that contains $x$. Since $\V b$ is $\ctwo{}$-close to an optimal vector
(namely, a vector
$\V a\in\I S_m$ that satisfies $\lambda(\blow{B}{\V a})=\lambda(\C G)$), we have that
 $$
  \frac{1}{\kappa}\frac{\partial}{\partial_i}\lambda(\blow{B}{\V b})
  \ge \frac{1}{\kappa}\frac{\partial}{\partial_i}\lambda(\blow{B}{\V a})-O(\ctwo{})
  =\lambda(\blow{B}{})-O(\ctwo{}).
   $$
Thus $\lambda(B',x)-\lambda(G,x)\ge \cspecial-O(\ctwo{})\ge \cspecial/2$ for each $x\in S$ and invoking \eqref{eq:020} and \eqref{eq:021}, we get
\begin{equation}\label{eq:B'-G}
 \Lambda(B')-\Lambda(G)\ge |W'|\ca n^{\kappa -2}/2 +|S|\frac{\cspecial{}}2 {n-1\choose
 \kappa-1}- O(|W'|\sigma n^{\kappa-2}+|S|^2n^{\kappa-2}).
 \end{equation}
\hide{
Here is the explanation for the third term. Before modification, each
vertex of $x\in S$ had $\lambda(G,x)$ by $\cspecial$ smaller than $\lambda(\C G)$,
which is by at least $\cspecial-\clambda/2<0.9\cspecial$ smaller than $\lambda(G)$,
the average over all $x$. Each new attachment of $x\in S$ follows $B'=\blow{B}{V_1,\dots,V_m}$
and we know that $\V b=(|V_1|/n,\dots,|V_m|/n)$ is $\ctwo{}$-close to an optimal $\V a$.
At the exact optimality, each vertex $x$ would satisfy $\lambda(\blow{B}{\V a},x)=\lambda(\C G)$, so for $B'$ we have to add error term of order $O(\ctwo{})$. Thus
when we replace the neighbourhood of each vertex $x\in S$ from that in $G$ by that
in $B'$, we win at least $0.9\cspecial-O(\ctwo{}+\cw+\sigma^2)>\cspecial/2$, where the $\cw$-term comes from wrong pairs disjoint from $S$ while $\sigma^2$ counts the contribution to $\lambda$
from $\kappa$-sets having at least two vertices in $S$.
}
By~\eqref{eq:B'-G} (and our bounds on $|W'|\le |W|\le \cw{n\choose 2}$ and $|S|/n=\sigma\le 2\czero{}/\cspecial\ll \min(\ca,\cspecial)$), we have that, for example,
\marginpar{$\ca,\cspecial\gg \czero{}$}
 $$
  \Lambda(n,\C G)-\Lambda(G)\ge \Lambda(B')-\Lambda(G)\ge |W'|\ca n^{\kappa-2}/4+|S|\frac{\cspecial}4{n-1\choose \kappa-1}
  \ge c_3(
  |W'|+|S|n)\frac{{n\choose k}}{{n\choose2}}.
   $$
Observing that $|W'|+|S|n\ge |W|$ and $\dedit(G,\blow{B}{})=|W|/{n\choose 2}$, we derive the perfect stability.\epf

\hide{
{\color{blue}
Note that the condition that the graphs in $\C F$ are twin-free is needed in Theorem~\ref{th:exact}. For example, suppose that
$\C F$ consists only of the 2-uniform blow-up $K_3(2,2,2)$ of $K_3$, then it is easy to see every $\C F$-free
graph of order $n\to\infty$ is $o(n^2)$-close to a $K_3$-free graph. (This directly follows from the removal lemma but also a direct proof is possible, using the supersaturation method
of Erd\H os and Simonovich.) Thus the Tur\'an problem for $K_3(2,2,2)$ is classically $(K_2,(1/2,1/2))$-stable.
(If needed, one can probably show that it is even robustly stable.) However, the conclusion of Theorem~\ref{th:exact}
fails: take the Tur\'an graph and connect one vertex to everything in its part.
In fact, we can take a Tur\'an graph and add a $C_4$-free graph into one part, so
the Tur\'an problem for $K_{2,2,2}$ is not even robustly stable.
}
}

\begin{theorem}[Perfect Stability II]\label{th:exact2}
 Suppose that Assumptions~\ref{as:general} and \ref{as:exact} are satisfied,
the problem is robustly $B$-stable and $B$ is $\lambda$-minimal.
Then the problem is perfectly $B$-stable.
\end{theorem}

\bpf
Clearly, the perfect stability will follow by Theorem~\ref{th:exact} if we show that its Assumptions~\ref{it:clstable},~\ref{it:strict} and~\ref{it:flip} are satisfied. Assuming that the problem is robustly $B$-stable, we trivially have that the problem is classically $B$-stable, that is, Assumptions~\ref{it:clstable} of Theorem \ref{th:exact} is satisfied. 
Thus it is enough to verify Assumptions~\ref{it:strict} and~\ref{it:flip} of Theorem~\ref{th:exact}. 

Roughly speaking, our proof is based on the following idea. For example, suppose that Assumption~\ref{it:strict} (the strictness of $\lambda$) fails. Let this be witnessed by a vector $\V y\in [0,1]^m$. Then we take a blow-up $G=\blow{B}{V_1,\dots,V_m}$ of order $n$ with optimal part ratios and add a set $Z$ of $\e n$ twin vertices, each attached to $G$ according to $\V y$. Since $\V y$ is $\Omega(1)$-far from each canonical attachment $\V v_i$, the new graph $G'$ has normalised edit distance $\Omega(\e)$ to the family~$\blow{B}{}$. On the other hand, if we take a random $\kappa$-subset $X\subseteq V(G')$ then it either is disjoint from $Z$ (and the conditional expectation of $\gamma(G'[X])$ is exactly $\lambda(G)$), or contains exactly one vertex of $Z$ (and the conditional expectation of $\gamma(G'[X])$ is $\lambda(\C G)+o(1)$ by the choice of $\V y$), or contains at least two vertices of $Z$ (which has probability $O(\e^2)$). We conclude that $|\lambda(G')-\lambda(G)|=O(\e^2)$, a contradiction to the robust stability. 
Likewise, if some edge flip violates Assumption~\ref{it:flip} (the flip aversion of $\lambda$), then one ``magnifies'' this by flipping all pairs between two appropriately placed sets of size~$\e n$.

Let us continue with the formal proof. Let the robust stability of the problem be satisfied with constant~$C$. Given $\lambda$, $B$ and
$C$, we choose a small enough quantity $c>0$.

In order to prove that the problem is strict, we assume on the contrary that there exist a maximiser $\V a$ in $\I S_m$ of $\lambda(B\blow{}\cdot)$ and an admissible $\V y$ in $[0,1]^m$ violating  $\lambda$-strictness.
Since $B$ is $\lambda$-minimal, we have that each $a_i\ge c$. We set
 $$\delta:=\min_{i\in[m]}\|\V y - \V v_i\|_1>0,$$
 and we pick some positive real $\e$ satisfying $\e\ll \min(c,\delta)$.

Let $G$ be a blow-up $\blow{B}{V_1,\dots,V_m}$ on $n\to\infty$ vertices with $|V_i|/n\to a_i$.
Let $G'$ be obtained from $G$ by adding a set $Z$ of $\e n$ twins whose attachment to $V(G)$ is given
by the vector $\V y+o(1)$, where we insist that if $y_i=0$ (resp.\ $y_i=1$),
then each $z\in{}Z$ is adjacent to no vertex in $V_i$ (resp. every vertex in $V_i$).
Since $\V y$ is admissible and the graphs in $\C F$ are twin-free,
$G'$ is $\C F$-free. Since $R_{\V a}(\V y)=\lambda(\C G)$, we have that the average of $\g{G'[X]}$ over the $\kappa$-subsets $X$ of $V(G')$
with $|X\cap Z|=1$ is $\lambda(\C G)+o(1)$. Thus it follows that $\lambda(\C G)-\lambda(G')$ is
at most $O(\e^2)$. By robust stability, the normalised distance from $G'$ to some blow-up $B'=\blow{B}{U_1,\dots,U_m}$ of $B$
is $O(\e^2)$. Clearly, $\lambda(B')\ge \lambda(G')-O(\e^2)\ge \lambda(\C G)-O(\e^2)$.

Recall that we have partitions $V_1\cup\dots \cup V_m\cup Z=U_1\cup\dots\cup U_m$.
We have that each $|U_i|\ge cn$ for otherwise we obtain the contradiction that
 $$\lambda(B')\le \lambda'+O(c)<\lambda(\blow{B}{})-O(\e^2),
  $$
 where
$\lambda'<\lambda(\blow{B}{})$ is the maximum of $\lambda$ over all blow-ups of proper subgraphs of~$B$. Similarly, each $V_i$ has at least $cn$ elements.

\Claim{AlmostSame}{There is an automorphism $\sigma:[m]\to [m]$ of $B$ such that for each $i$
 \begin{equation}\label{eq:AlmostSame}
  |U_{\sigma(i)}\bigtriangleup V_i|\le 2\e n/c.
 \end{equation}
 }

\bcpf We show first that for each $i\in{}[m]$ there exists $\sigma(i)\in{}[m]$ satisfying \eqref{eq:AlmostSame} (and then
observe that the map $\sigma:[m]\to [m]$ is an automorphism of $B$). Take any $i\in{}[m]$. Suppose that there is no choice of $\sigma(i)$ satisfying~\eqref{eq:AlmostSame}.
We pick $x\in{}[m]$ such that $|U_x\cap V_i|\ge|V_i|/m\ge cn/m>\e n/cm$. We
distinguish the following two cases.\medskip

\noindent\textbf{Case I:} There exists $y\in{}[m]$ such that $y\neq x$ and $|U_y\cap V_i|>\e n/cm$.\medskip

\noindent Since $B$ is twin-free (which follows by the $\lambda$-minimality of $B$),
pick $h\in[m]$ such that exactly one of $x,y$ is a $B$-neighbour of $h$. Then,
every $v\in{}U_h\setminus Z$ is incident to at least $\e n/cm$ pairs on which the graphs $G'$ and $B'$ differ,
because each $v\in{}U_h\setminus Z$ has different $B'$-adjacencies to $V_i\cap U_x$ and $V_i\cap U_y$ but the same
$G'$-adjacency to all vertices of $V_i\supseteq (V_i\cap U_{x})\cup (V_i\cap U_{y})$. Thus $\Dedit(G',B')\ge
(1/2)\cdot |U_h\setminus Z|\cdot (\e n/cm)$ which is not $O(\e^2n^2)$, a
contradiction.\medskip

\noindent\textbf{Case II:} For every $y\in{}[m]$ such that $y\neq x$ we have that $|U_y\cap V_i|\le\e n/cm$.\medskip

\noindent It holds that $|V_i\setminus U_x|\le\e n/c$. Since we work under the assumption that there is no appropriate choice of $\sigma(i)$, we have, in particular, that $|U_{x}\bigtriangleup V_i|> 2\e n/c$ and therefore $|U_x\setminus V_i|>\e n/c$. We pick $j\in{}[m]$ with $j\neq i$ such that $U_x\cap V_j>\e n/cm$. Arguments similar to the ones used in Case I lead to a contradiction.

To complete the proof we show that $\sigma$ is an automorphism of $B$.
Let us observe that $\sigma$ is an injection. Indeed, suppose on the contrary that there exist $i,j$ and $x$ in $[m]$ such that $i\neq j$ and $\sigma(i)=\sigma(j)=x$.
Then we have that
\[|U_x\bigtriangleup V_j|\geqslant|U_x\setminus V_j|\geqslant|U_x\cap V_i|
\stackrel{\eqref{eq:AlmostSame}}{\geqslant}|V_i|-2\e n/c\geqslant cn-2\e n/ c\]
contradicting \eqref{eq:AlmostSame}. To prove that $\sigma$ is edge and non-edge preserving, we assume on the contrary that there exists a pair of nodes $ij$ such that $\sigma$ does not preserve adjacency. Then the graphs $G'$ and $B'$
differ on every pair $uv$ with $u\in V_i\cap U_{\sigma(i)}$ and $v\in V_j\cap U_{\sigma(j)}$ generating at least
$((c-2\e/c)n)^2\gg \e^2n^2$
such pairs. The latter is a contradiction to $\Dedit(G',B')=O(\e^2n^2)$. The claim is proved.
\ecpf

By relabelling $U_1,\dots,U_m$, we can assume that the bijection $\sigma$ of Claim~\ref{cl:AlmostSame}
is the identity map. We expand $(V_i)_{i=1}^m$ to a partition $(V'_i)_{i=1}^m$ of the vertex set of $G'$ setting $V'_i=V_i\cup (U_i\cap Z)$ for each $i\in{}[m]$. Clearly
\begin{equation}\label{eq:AlmostSameNew}
  |V'_i\bigtriangleup U_i|\leqslant 2\e n/c
\end{equation}
for all $i\in{}[m]$. Finally, we set
\[\Delta_1:=E(G')\bigtriangleup E(\blow{B}{V'_1,\dots,V'_m})\text{
and } \Delta_2:=   E(B') \bigtriangleup E(\blow{B}{V'_1,\dots,V'_m}) .\]
Each vertex $v\in{}Z$
is adjacent to at least $\delta n/2$ pairs in $\Delta_1$,
because $\V y$ is $\delta$-far from $\V v_1,\dots,\V v_m$, and at most $2\e m n/c$ pairs in $\Delta_2$.
Thus the symmetric difference between $G'$ and $B'$
is at least $\e n\times (\delta/2-2\e m/c)n\gg \e^2v(G')^2$, a contradiction which shows that
the graph $B$ is $\lambda$-strict.

Next, let us prove the $\lambda$-flip-aversion of $B$. We pick some positive real $\e\ll c$ and
towards a contradiction we assume that there exists some integer $n$ with $n>1/\e^3$,
an almost optimal blow-up $B'=\blow{B}{V_1,\dots,V_m}$ on $[n]$ and some pair $x,y$ of distinct nodes
such that the graph $B'\oplus xy$ contains no forbidden graph of order at most $m+2$ and
\begin{equation}
  \label{eq032}
  \Lambda(B')-\Lambda(B'\oplus xy)<\e^3 n^{\kappa-2}.
\end{equation}

Let $i,j\in{}[m]$ be such that $x\in V_i$ and $y\in V_j$. We pick subsets $X$ and $Y$ of $V_i$ and $V_j$ respectively with cardinality $\e n$ each. If $i=j$ then we choose $X$ and $Y$ to be disjoint.
Let $\mathcal{B}$ be the set of all pairs of nodes with one node in $X$ and one in $Y$. Also let $G$ be the graph obtained by flipping the adjacency between each pair in $\mathcal{B}$.
Since each of $X$ and $Y$ consists of twins, $G$ does not contain any forbidden subgraph.

Let us show that
\begin{equation}\label{eq:034}
  \lambda(B')-\lambda(G)\leqslant O(\e^3).
\end{equation}
Indeed, let $\mathcal{A}$ be the set of all $\kappa$-element subsets of $V=V_1\cup...\cup V_m$.
We partition $\mathcal{A}$ into $\mathcal{A}_0$, $\mathcal{A}_1$ and $\mathcal{A}_{\geqslant 2}$, the set of all
$Z\in{}\mathcal{A}$ containing  respectively zero, one and at least two pairs of $\mathcal{B}$. Finally, for each $e\in{}\mathcal{B}$, we set $\mathcal{A}^e$, $\mathcal{A}_{1}^e$ and  $\mathcal{A}^e_{\geqslant2}$ to be the set of all $Z\in{}\mathcal{A}$, $\mathcal{A}_{1}$ and $\mathcal{A}_{\geqslant 2}$ respectively, containing $e$. Note that if $Z\in\mathcal{A}_{\geqslant 2}$, then $|Z\cap (X\cup Y)|\ge 3$ and thus $|\mathcal{A}_{\geqslant 2}|=O(\e^3n^{\kappa})$. We are going to use this fact a couple of times in the following chain of equalities.
 \hide{
 \begin{equation}
   \label{eq:035}
   \begin{split}
     \Lambda(B')-\Lambda(G) &= {n\choose \kappa}^{-1}\sum_{Z\in\mathcal{A}}\left(\lambda(B'[Z])-\lambda(G[Z])\right) \\
     &=  {n\choose \kappa}^{-1}\sum_{e\in\mathcal{B}}\sum_{Z\in\mathcal{A}_1^e}
     \left(\lambda(B'[Z])-\lambda(G[Z])\right)+O(\e^3)\\
     &= {n\choose \kappa}^{-1}\sum_{e\in\mathcal{B}}\sum_{Z\in\mathcal{A}_1^e}\lambda(B'[Z])-\lambda(B'\oplus e[Z])+O(\e^3) \\
     &= {n\choose \kappa}^{-1}\Big(\sum_{e\in\mathcal{B}}\sum_{Z\in\mathcal{A}^e}\lambda(B'[Z])-\lambda(B'\oplus e[Z])\\
     &\;\;\;\;\;\;\;\;\;\;-\sum_{e\in\mathcal{B}}\sum_{Z\in\mathcal{A}_{\geqslant 2}^e}\lambda(B')[Z]-\lambda(B'\oplus e)[Z]\Big)+O(\e^3)\\
     &={n\choose \kappa}^{-1}\sum_{e\in\mathcal{B}}\Lambda(B')-\Lambda(B'\oplus e)+O(\e^3)\\
     &={n\choose \kappa}^{-1}\sum_{e\in\mathcal{B}}\Lambda(B')-\Lambda(B'\oplus xy)+O(\e^3)
     \stackrel{\eqref{eq032}}{\leqslant}O(\e^3).
   \end{split}
 \end{equation}}
\begin{equation*}
  \begin{split}
    \Lambda(B')-\Lambda(G) &= \sum_{Z\in\mathcal{A}}\left(\lambda(B'[Z])-\lambda(G[Z])\right) \\
    &=  \sum_{e\in\mathcal{B}}\sum_{Z\in\mathcal{A}_1^e}
    \left(\lambda(B'[Z])-\lambda(G[Z])\right)+O(\e^3n^{\kappa})\\
    &= \sum_{e\in\mathcal{B}}\sum_{Z\in\mathcal{A}_1^e}\left(\lambda(B'[Z])-\lambda((B'\oplus e)[Z])\right)+O(\e^3n^{\kappa}) \\
    &=\sum_{e\in\mathcal{B}}\left(\Lambda(B')-\Lambda(B'\oplus e)\right)+O(\e^3n^{\kappa})\\
    &=\sum_{e\in\mathcal{B}}\left(\Lambda(B')-\Lambda(B'\oplus xy)\right)+O(\e^3n^{\kappa})
    \stackrel{\eqref{eq032}}{\leqslant}O(\e^3n^{\kappa}).
  \end{split}
\end{equation*}
Therefore, by the almost-optimality of $B'$ we have that $|\lambda(G)-\lambda(\mathcal{G})|\leqslant O(\e^3)$. By the assumed
robust stability, there exists some blow-up $B''=\blow{B}{U_1,...,U_m}$ of $B$ such that $\dedit(B'',G)=O(\e^3)$.
Following arguments as in the proof of Claim~\ref{cl:AlmostSame}, we can assume
 that
$|V_h\bigtriangleup U_h|\le \e n/c$ for every $h\in[m]$.

Then we distinguish the following three (non-exclusive) cases.
\begin{enumerate}
  \item[(i)] $|X\setminus U_i|\ge \e n/2$.
  \item[(ii)] $|Y\setminus U_j|\ge \e n/2$.
  \item[(iii)] $|X \cap U_i|> \e n/2$ and $|Y \cap U_j|> \e n/2$
\end{enumerate}
We complete the proof by showing that each case leads to a contradiction and, in particular, we show that each case yields that
$\dedit(G,B'')=\Omega(\e^2)$. Indeed, let us assume (i). Then there is $i'\not=i$ such that $|X\cap U_{i'}|\ge
\e n/2m$. Pick $h\in [m]$ such that the $B$-adjacencies of $\{h,i\}$ and $\{h,i'\}$ differ. We have at least
$(cn-\e/c)n$ vertices in $U_h\cap V_h$.
Thus the symmetric difference
between $G$ and $B''$ is at least $(cn-\e/c)n\times \e n/2m\gg \e^2 n$. Likewise, case (ii) leads to a contradiction. Finally assuming case (iii) we have that $G$ and $B''$ differ on every pair with one node in
$X \cap U_i$ and one in $Y \cap U_j$. Thus the symmetric
difference between $G$ and $B''$ is at least $(\e n/2)^2$. \epf

\section{Finding optimal asymptotic part ratios}\label{se:asympt_part_rations}

In this section, we provide some analysis related to the values of $\V a$ in $\I S_m$ that maximise the
function $\lambda(\blow{B}{\cdot})$.

While in all examples from
Section~\ref{se:examples} the optimal vector $\V a$ was uniform, this is not
always the case.
For example, it was conjectured in~\cite{pikhurko+vaughan:13} (based on the
numerical evidence
from Flagmatic) that the asymptotically extremal value for Erd\H os' $f(n,4,4)$-problem is attained
by a blow-up of a specific 8-part graph $B$. If the conjecture is true, then
the optimal blow-up of
$B$ that minimises the number of $\overline K_4$-subgraphs is not uniform (in fact,
the optimal part ratios
are some irrational numbers). Alternatively, here is a simple although rather
artificial example that illustrates the point.

\begin{example}[Simple problem with a non-uniform optimal vector]
 Let $\C F$ consist of all odd cycles plus the graph with 3 vertices and one
edge. Then $\C F$-free graphs
on $[n]$ are exactly complete bipartite graphs, that is, blow-ups of $B=K_2$.
 Let $\kappa=6$. Let $\g{H}=0$ except one defines $\g{H}$ for
 $H\in\{K_{0,6},K_{1,5},K_{2,4},K_{3,3}\}$ so that
 $\lambda(K_{xn,(1-x)n})=p(x)+o(1)$, where, e.g.\
 $$
 p(x)=12(x-1/2)^6-217(x-1/2)^4+24(x-1/2)^2.
 $$
 This polynomial $p$ is symmetric around $1/2$ and its maximum on $[0,1]$ is
 attained at
 $x_0=(3-\sqrt{2})/6=0.264...$ and $1-x_0$. Finding
the maximum of $\lambda(\Forb{\C F})=\lambda(\blow{B}{})$ over $\I
S_2=\{(x,1-x)\mid x\in[0,1]\}$
amounts to optimising $p(x)$ over $x\in[0,1]$ which is not attained for
$(1/2,1/2)$.\qed\end{example}


Let us prove a sufficient condition that implies the uniqueness of the maximiser and happens to apply to many concrete problems.


\begin{lemma}\label{lm:unique} Let all assumptions of Theorem~\ref{th:stab1} apply.
View the graph $\tau$ from Assumption~2 also as a type and assume additionally that the flag algebra certificate $\C C$ includes a matrix $Q^\tau$
of co-rank 1 associated to~$\tau$. 
Then the vector $\V a$ is the unique  maximiser of $\lambda(B\blow{}{\cdot})$ in $\mathbb{S}_m$.
%
\end{lemma}

\bpf Let $\V b\in\I S_m$ be a maximiser of $\lambda(\blow{B}\cdot)$. By Assumption~2(b), we have $\lambda(\Forb{\C F}\cup\{\tau\})<\lambda(\Forb{\C F})=\lambda(\blow{B}{\V b}))$. Thus there is a strong homomorphism $f$ from $\tau$ into $B[\{i\in[m]:b_i>0\}]$.
Fix one such $f$. 

For large $n$, let $G=\blow{B}{V_1,\ldots,V_m}$ with $|V_i|=b_in+O(1)$ and take an (injective) embedding $\psi:V(\tau)\to V(G)$ such that $\psi(x)\in V_{f(x)}$ for every $x\in V(\tau)$. Define $\V x_{\V b}$ to be the limit as $n\to\infty$ of the vector $\V x$ from~\req{x} normalised so that the sum of entries is~$1$.  Clearly, the limit does not depend on the choice of $\psi$. Arguing as in the proof of Lemma~\ref{lm:forced}, we conclude that $\V x_{\V b}$ is a zero eigenvector of~$Q^\tau$. Of course, the same applies to the vector~$\V x_{\V a}$. Since $Q^\tau$ is of co-rank~1,
we have that $\V x_{\V b}=\V x_{\V a}$. However, $\V b$ is uniquely determined from $\V x_{\V b}$. Namely,  by Assumption~2(c), the $i$-th entry $b_i$ is the $\ell$-th root, $\ell:=(N-v(\tau))/2$, of the entry of $\V x_{\V b}$ that corresponds to the $\tau$-flag obtained by adding some $\ell$ new vertices from $V_i$ to the $\psi$-image of $\tau$; this follows by recalling that the vector $\V x_{\V b}$ encodes the limiting distribution of the $\tau$-subflag of $G$ induced by a random $\kappa$-subset containing $\psi(V(\tau))$.
Thus $\V b=\V a$ and $\V a$ is indeed the unique maximiser of $\lambda(\blow{B}\cdot)$ in $\I S_m$.
\epf


\section{Computer implementation}

Combining Theorems \ref{th:stab1} and \ref{th:exact2} we obtain the following result, which provides sufficient conditions for perfect stability.
The verification of these conditions can be carried out by a computer. In the next section we include such applications.

\begin{theorem}
\label{thm:pc}
 Let Assumption \ref{as:general} and Part~\ref{as:exact.i} of Assumption~\ref{as:exact}
apply. Also, we assume all of the following.
\begin{enumerate}
 \item\label{it:FAProof} We have a vector $\V a\in\I S_m$ with no zero entries and a certificate
 $\C C$ on $N$ vertices that proves $\lambda(\C G)\le \lambda(\blow{B}{\V a})$. (Thus, by Assumption~\ref{as:general}.3,
 we know that $\lambda(\C G)=\lambda(\blow{B}{\V a})$.)
   \item\label{it:tau} There is a graph $\tau$ of order at most $N-2$ satisfying the following.
   \begin{enumerate}
     \item[(a)] $\lambda(\Forb{\C F})>\lambda(\Forb{\C F\cup\{\tau\}})$.
     \item[(b)] There exists the unique
(up to automorphisms of $\tau$ and $B$) strong homomorphism $f$ from $\tau$ into $B$.
     \item[(c)] For every distinct
$x_1$ and $x_2$ in $V(B)$ we have $\Gamma_B(x_1)\cap f(V(\tau))\not= \Gamma_B(x_2)\cap f(V(\tau))$.
   \end{enumerate}
 \item\label{it:sharp} Every $\C C$-sharp graph of order $N$ admits a
strong homomorphism
into $B$.
\end{enumerate}
Additionally, suppose that at least one of the following two statements holds:
\begin{enumerate}
  \item[(i)] \label{it:rk}the certificate $\C C$ contains (as a type) the graph $\tau$ from Assumption~\ref{it:tau} above and the corresponding matrix $Q^\tau$ in $\C C$ is of co-rank 1, or
  \item[(ii)]\label{it:banB} 
  $\lambda(\Forb{\C F\cup\{B\}})<\lambda(\Forb{\C F})$.
\end{enumerate}
Then the problem is
  perfectly $B$-stable.
\end{theorem}

\bpf Clearly, all assumptions of Theorem~\ref{th:stab1} are satisfied, so the problem is robustly $B$-stable.
By
Theorem~\ref{th:exact2}, it is enough to check only that $B$ is $\lambda$-minimal.

If Condition (i) holds, then the $\lambda$-minimality of $B$ follows from Lemma~\ref{lm:unique} (and
the assumption that $\V a$ has no zero entries).
So assume that Condition (ii) holds. Let $B'$ be an arbitrary proper subgraph of $B$ and let $B''$ be any blow-up of $B'$ on $n\to\infty$ vertices. Since $B$ is twin-free by Condition 2(c), we have that $B''$ is $B$-free and thus $B''$ belongs to
$\Forb{\C F\cup\{B\}}$. Thus
 $$
 \lambda(B'')\leqslant\lambda(\Forb{\C F\cup\{B\}})+O(1/n)<\lambda(\Forb{\C F})+O(1/n).
 $$
 Again, we conclude that $B$ is $\lambda$-minimal, as desired.
 \epf

\brm If Assumptions~\ref{it:FAProof}, \ref{it:tau}, \ref{it:sharp} and (i)\hide{\ref{it:rk}} of Theorem \ref{thm:pc} are satisfied, then we have that $\lambda(\blow{B}\cdot)$ admits $\V a$ as a unique maximiser (see Lemma~\ref{lm:unique}). This is not the case if Assumptions~\ref{it:FAProof}, \ref{it:tau}, \ref{it:sharp} and (ii)\hide{\ref{it:banB}} of Theorem \ref{thm:pc} are satisfied, when the uniqueness of $\V a$ as a maximiser of $\lambda(\blow{B}\cdot)$ is not guaranteed. In this case, one has to investigate the uniquess of $\V a$ by other means.
\hide{In particular, the uniqueness of $\V a$ as a maximiser of $\lambda(\blow{B}\cdot)$ is not guaranteed. However, Condition (ii)\hide{\ref{it:banB}} of Theorem \ref{thm:pc} implies that $B$ is $\lambda$-minimal. If we have in addition that the equation $M \V z=\V e_1$, where $M$ is as defined in Section \ref{se:asympt_part_rations}, has unique solution, then $\V a$ is the unique maximiser of $\lambda(\blow{B}\cdot)$ (see Lemma \ref{lm:stab2}, as well as, the remark after it).}

\section{Applications of the general theorems}\label{Applications}

Below is a list of results that directly follow by Theorem \ref{thm:pc} by
running our computer code. The ancillary folder of the arxiv version of this
paper contains, for each problem, the flagmatic script \texttt{*.sage} which
was used to generate the certificate and the transcript of the session
\texttt{*.txt} when the code is run. Due to arxiv's file size limitations, ancillary folder only contains certificates \texttt{*.js} with $N \leq 6$. All certificate files are in Flagmatic's Github directory at:

\begin{center}\url{https://github.com/jsliacan/flagmatic/tree/master/certificates}.\end{center}

For example, for the $f(n,4,3)$-problem discussed in
Section~\ref{se:fkl}, these are the files \texttt{f43.sage}, \texttt{f43.txt}
and \texttt{f43.js} respectively.

The reader who would like to verify these results has the following options.
\begin{description}

\item[Generate certificates from scratch using flagmatic:] For this the reader would need to install our version of \emph{Flagmatic} (which is built upon version 2.0 of Emil Vaughan), the \emph{Sage} environment, and an SDP solver such as \emph{CSDP} or \emph{SDPA/SDPA-DD}. The required version of \emph{Flagmatic} can be downloaded from this URL:\\
  $$
  \mbox{\url{https://github.com/jsliacan/flagmatic}}
  $$
which in particular contains a \texttt{README.md} file with directions on how
to install it and run our scripts.

\item[Run our verifier script \texttt{inspect\_certificate.py}:] This stand-alone script (which is written in \emph{Python/Sage} and uses exact arithmetic) can be used to verify the bound given by each certificate. It is available at the above URL. Its source code is relatively short and well-documented. (Also, the Appendix to the arXiv version of this paper~\cite[Appendix]{OSTarxiv} contains some further notes on our implementation.) For example, the complete verification of the certificates \texttt{f43.js}, \texttt{f43\_stab.js} can be invoked with the following shell command:
  $$
  \mbox{\scriptsize \texttt{sage -python inspect\_certificate.py f43.js --stability 3/25 "4:121324" "5:1223344551" f43\_stab.js}}
  $$
The full details on how to use the \texttt{inspect\_certificate.py} verification script can be found at the end of the \texttt{README.md} file at \url{https://github.com/jsliacan/flagmatic}.
\item[Write an independent verifier:] The information on how the data inside our certificate files are organised can be found in~\cite[Appendix]{OSTarxiv}.
 \end{description}

\hide{

{\color{blue} We need to fix the introduction. Also we should just mention what Flagmatic is, since we refer to it in the examples.

The examples below should be thought of as a proof of concept for Flagmatic stability verification. In all these cases, the conditions in robust and perfect stability can be verified readily by hand (save for the flag algebra computations proving that forbidding $\tau$ and $B$ (one by one) reduces objective value).}
}

In the following, we describe some of the input values (such as $N$ and $\V a$)
that determine $\lambda(\C G)$ and prove perfect stability in
Theorems~\ref{th:ErdosPr3}--\ref{th:Y}.

\begin{description}
  \hide{\color{blue}\item \emph{Tur\'an problem for cliques}: For each integer
                $k$ with $k\ge3$, we have $\C F=\{K_k\}$, $\kappa=2$ and
                $\gamma$ is equal to zero except $\gamma(K_2)=1$. The result
                follows by a straightforward application of Theorem \ref{thm:pc}
                for $N=k$, $B=K_{k-1}$, $\V a$ the vector in $S_{k-1}$ having
                each entry equal to $k-1$ and $\tau=K_{k-2}$.}
  \item[\emph{Minimising the number of independent
  sets in triangle-free graphs (Theorem~\ref{th:ErdosPr3}):}]
        Recall that $k\in\{4,...,7\}$,    $\C F = \{K_3\}$, $\kappa = k$ and $\gamma$ is equal to zero
         except $\gamma(\overline{K}_k)=-1$. Theorem \ref{th:ErdosPr3} for
         $k=4,5$ follows
         by Theorem \ref{thm:pc} for $N=5$, $B=C_5$, $\V
         a=(1/5,\ldots,1/5)\in\I S_{5}$
         and $\tau=K_2\cup K_1$, that is, the disjoint union of an edge and
         a single node (see scripts \texttt{f43.sage} and \texttt{f53.sage}).

         Unfortunately, our code could not generate certificates
         when $k \in \{6,7\}$. This computationally demanding task (with $N=8$) seems to be very sensitive to the obtained numerical SDP solution and the version of \textit{Sage}. However, the corresponding certificates
         have already been produced by Pikhurko and Vaughan~\cite{pikhurko+vaughan:13}
          and
         we include them in the arxiv version of this paper. By running our script \texttt{inspect\_certificate.py} on them, one can confirm that the problem is perfectly stable in these two cases, where we let $B$ be the Clebsch graph,
         $\V a=(1/16,\dots,1/16)\in\I S_{16}$, and $\tau$ be the 5-cycle $C_5$
         with one isolated vertex added. (Interestingly, the correct asymptotic
         of $f(n,6,3)$ can be obtained already for $N=7$ but we could not
         satisfy  Condition~2 of Theorem~\ref{thm:pc} with this $N$.) In all the cases above, the uniqueness of $\V a$ follows from Lemma~\ref{lm:unique}, since the corank of $Q_\tau$ is 1 for each $k\in \{4,\ldots,7\}$.

  \item[\emph{Maximising the number of pentagons in triangle-free graphs
  (Theorem~\ref{th:pentagons}):}]
  Recall that the problem
        is defined by $\C F=\{K_3\}$, $\kappa=5$, and $\gamma(H)$ equals zero,
except $\gamma(C_5)=1$. Theorem \ref{th:pentagons} follows by Theorem \ref{thm:pc}
        for $N=5$, $B=C_5$, $\V a$ the vector in $\I S_5$ having each entry
        equal to $1/5$ and $\tau=K_2\cup K_1$, that is, the disjoint union of
        an edge and a single node. The uniqueness of $\V a$ follows from Lemma~\ref{lm:unique}, since the co-rank of $Q_\tau$ is 1.

  \item[\emph{Inducibility of the cycle on four vertices
  (Theorem~\ref{th:c4}):}] Recall that the
  problem
        is defined by $\C F=\emptyset$, $\kappa=4$, and $\gamma(H)$ equals zero,
        except $\gamma(C_4)=1$. Theorem \ref{th:c4} follows by Theorem \ref{thm:pc}
        for $N = 5$, $B = K_2$, $\V a=(1/2,1/2)$ and $\tau = K_1$. The uniqueness of $\V a$ follows from Lemma~\ref{lm:unique}, since the co-rank of $Q_\tau$ is 1.

  \item[\emph{Inducibility of $K_4$ minus an edge (Theorem~\ref{th:k4-}):}]
  Recall that the problem is
  defined by
        $\C F=\emptyset$, $\kappa=4$, and $\gamma(H)$ equals zero, except $\gamma(K_4^-)=1$.
        Theorem \ref{th:k4-} follows by Theorem \ref{thm:pc} for $N=7$, $B=K_5$, $\V a$ the vector in $\I S_5$ having each entry equal to $1/5$ and $\tau=K_5$. The uniqueness of $\V a$ follows from Lemma~\ref{lm:unique}, since the co-rank of $Q_\tau$ is 1.

  \item[\emph{Inducibility of $K_{3,2}$ (Theorem~\ref{th:k32}):}] Recall that
  the problem is defined by
        $\C F=\emptyset$, $\kappa=5$, and $\gamma(H)$ equals zero,
        except $\gamma(K_{3,2})=1$. Theorem \ref{th:k32} follows by Theorem \ref{thm:pc}
        for $N=6$, $B = K_2$, $\V a=(1/2,1/2)$ and $\tau=K_2$. The uniqueness of $\V a$ follows from Lemma~\ref{lm:unique}, since the co-rank of $Q_\tau$ is 1.

  \item[\emph{Inducibility of $K_{2,2,1}$ (Theorem~\ref{th:k221}):}] Recall
  that the problem is defined
  by
        $\C F=\emptyset$, $\kappa=5$, and $\gamma(H)$ equals zero, except $\gamma(K_{2,2,1})=1$.
        Theorem \ref{th:k221} follows by Theorem \ref{thm:pc}
        for $N=6$, $B = K_3$, $\V a=(1/3,1/3,1/3)$ and $\tau=K_2$. The uniqueness of $\V a$ follows from Lemma~\ref{lm:unique}, since the co-rank of $Q_\tau$ is 1.

  \item[\emph{Inducibility of $P_3\cup K_2$ (Theorem~\ref{th:p3k2}):}] Recall
  that the problem is
  defined by
        $\C F=\emptyset$, $\kappa=5$, and $\gamma(H)$ equals zero, except $\gamma(P_3\cup K_2)=1$.
        Theorem \ref{th:p3k2} follows by Theorem \ref{thm:pc}
        for $N = 6$, $B = K_3 \cup K_3$, that is, the disjoint union of two triangles, $\V a$ the
        vector in $\I S_6$ having each entry equal to $1/6$ and $\tau = K_2 \cup K_2$, that is, the
        disjoint union of two edges. The uniqueness of $\V a$ follows from Lemma~\ref{lm:unique}, since the co-rank of $Q_\tau$ is 1.

  \item[\emph{Inducibility of the ``Y'' graph (Theorem~\ref{th:Y}):}] Recall
  that the problem is
  defined by
        $\C F=\emptyset$, $\kappa=5$, and $\gamma(H)$ equals zero, except
        $\gamma(\mathrm{Y})=1$.
        Theorem \ref{th:p3k2} follows by Theorem \ref{thm:pc}
        for $N = 6$, $B = C_5$, that is, the cycle on 5 vertices, $\V a$ the
        vector in $\I S_5$ having each entry equal to $1/5$ and $\tau = P_4$, that is, the
        path on 4 vertices. The uniqueness of $\V a$ follows from Lemma~\ref{lm:unique}, since the co-rank of $Q_\tau$ is 1.

\end{description}

\subsection{Inducibility of the Paw graph}\label{paw}

The value $i(F_{\text{paw}})$ has been calculated
by Hirst~\cite{hirst-inducibility}, where $F_{\text{paw}}$ is the paw graph,
that is, the graph obtained by adding a pendant edge to a triangle. We work on the complementary problem.
We set $\C F=\emptyset$ and $\gamma$ the map taking the value 0
on every graph of order 4 except the disjoint union of $P_3$ and a single
vertex, that we denote by $F$, where it takes the value 1. From Hirst's work it
follows that $i(F)=3/8$ and an asymptotically extremal construction is a balanced blow-up of
the graph consisting of two disjoint edges.
In this section, we show that the problem is $K_2\cup K_2$-perfectly stable.

\begin{figure}[htb]
\centering \includegraphics[]{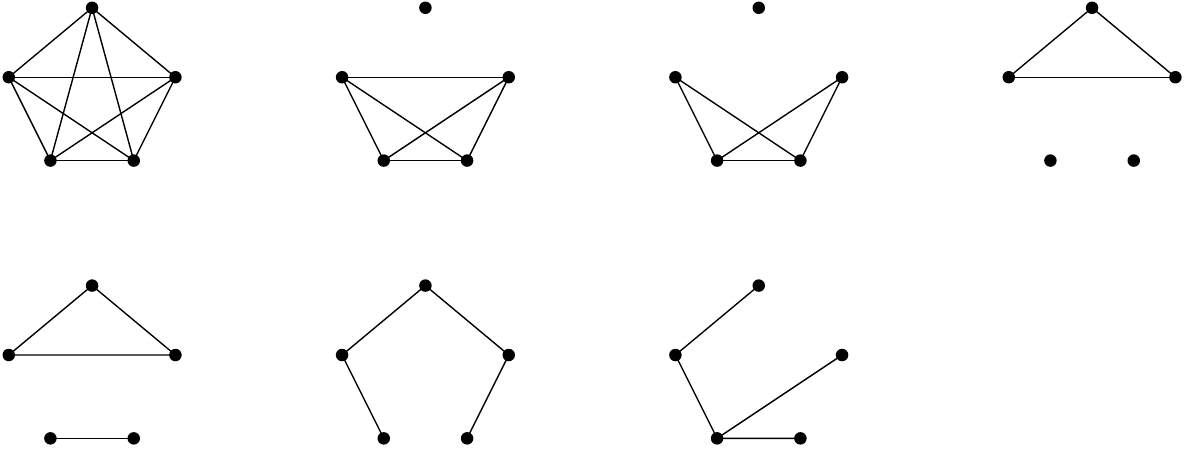}
\caption{Sharp graphs that are not a blow of $K_2\cup K_2$.}\label{fg:2}
\end{figure}

Unfortunately, our result does not follow directly by Theorem \ref{thm:pc},
since Condition \ref{it:sharp} does not hold for our flag algebra
certificate. In particular, according to our certificate the sharp graphs
consist of the blow-ups of $K_2\cup K_2$ on $5$ vertices and the graphs listed
in Figure~\ref{fg:2}.
Let us denote by $\C S$ the set of sharp graph on 5 vertices and by $\C{NS}$
the set of the non-sharp ones.
However, letting $B= K_2\cup K_2$, $N=5$ and $\tau$ be the disjoint union of an
edge and a single vertex, we have that Assumptions~\ref{it:FAProof} and
\ref{it:tau} of Theorem \ref{thm:pc} are satisfied. We refer to them as P1 and
P2 respectively.
The perfect stability of the problem follows by a sequence of lemmas.

\begin{lemma}
  \label{lm:PawMinimal}
  The graph $K_2\cup K_2$ is $\lambda$-minimal.
\end{lemma}
\bpf
Let $\mathbf{a}=(a_1,a_2,a_3,a_4)$ in $\I S_{4}$ and set $B=K_2\cup K_2$. It is
easy to see that
\[\begin{split}\lambda(\blow{B}{\mathbf{a}})
&=12 ( a_1^2a_2(a_3+a_4) + a_1a_2^2(a_3+a_4) + (a_1+a_2)a_3^2a_4 +
(a_1+a_2)a_3a_4^2 )\\
&=12(a_1a_2+a_3a_4)(a_1+a_2)(a_3+a_4).
\end{split}\]
\sloppy To prove the $\lambda$-minimality of $B$, it suffices by symmetry to show that the
maximum
value achieved by $\lambda(\blow{B}{\mathbf{a}})$ for
$\mathbf{a}=(a_1,a_2,a_3,a_4)$ in $\I S$ with $a_4=0$ is strictly less than
$\lambda(\blow{B}{(1/4,1/4,1/4,1/4)})=3/8$. Equivalently, it suffices to show
that the maximum of the map
\[f(x_1,x_2)=12x_1x_2(x_1+x_2)(1-x_1-x_2)\]
for $(x_1,x_2)$ in
$D=\{(x_1,x_2):x_1,x_2\geqslant0\;\text{and}\;x_1+x_2\leqslant1\}$, is strictly
less than $3/8$.

Indeed, observe that $f$ vanishes on the boundary of $D$ and therefore the
maximum is achieved in the interior of $D$. The partial derivatives of $f$ at
the maximum satisfy the following:
\begin{equation}
  \label{eq:partialdev1}
  \frac{\partial f}{\partial x_1}=12(x_2(x_1+x_2)(1-x_1-x_2) +
  x_1x_2(1-2x_1-2x_2))=0
\end{equation}
and
\begin{equation}
  \label{eq:partialdev2}
  \frac{\partial f}{\partial x_2}=12(x_1(x_1+x_2)(1-x_1-x_2) +
  x_1x_2(1-2x_1-2x_2))=0.
\end{equation}
Subtracting \eqref{eq:partialdev1} and \eqref{eq:partialdev2}, we obtain that
$x_1=x_2$. Plugging it into \eqref{eq:partialdev1}, we get that $f$ achieves a
maximum at $(3/8,3/8)$. Thus the maximum of $f$ is $81/256$, which is strictly
smaller that $3/8$ and the proof of the lemma is complete.
\epf

\begin{lemma}
  \label{lm:PawStab}
  The problem is classically $K_2\cup K_2$-stable.
\end{lemma}
\bpf
  Let $\e$ be a positive real. By the Induced Removal Lemma of Alon et
  al.~\cite{alon+fischer+krivelevich+szegedy:00} there exists a
  positive real
  $\eta$ such that for every graph $G$ of order at least $1/\eta$ satisfying
  $p(H,G)\leqslant\eta$ for all $H\in{}\C{NS}$, we have that there exists a
  graph $G'$ of the same order as $G$ such that $\dedit(G',G)\le \e$ and each induced subgraph of $G'$ belongs to $\C S$.
  By \eqref{eq:ai+bi}, there exists a positive real $\delta$ such that for each
  graph $G$ of order at least $1/\delta$ satisfying $\lambda(\C
  G)-\lambda(G)\leqslant\delta$, we have that $p(H,G)\leqslant \eta$ for all
  $H\in{}\mathcal{S}$ and therefore there exists some graph $G'$ of the same
  order as $G$ such that $\dedit(G,G')\le \e$
  and each induced subgraph of $G'$ belongs to $\C S$.

Since $G'$ is close to $G$ and $\lambda(G)$ is close to
  $\lambda(\C G)$, we get that $\lambda(G')$ is close to $\lambda(G)$ and
  therefore, by P2(a), we have that $\tau$ embeds into $G'$. Recall that $\tau$
  is
  the disjoint union of an edge and a single vertex. Without loss of
  generality, we may assume that $V(\tau)=[3]$ and $\{1,2\}$ forms an edge in
  $\tau$. Since $G'$ admits an induced copy of $\tau$, there exists an
  injective strong homomorphism
  $\psi:[3]\to V(G')$ between $\tau$ and $G'$. For every $s\in{}2^{[3]}$,
  where we view $2^{[3]}$ as the set of maps from $[3]$ to $\{0,1\}$, we define
   \begin{equation}
     V'_s=\{x\in V(G')\setminus\mathrm{Im}(\psi):\{x,\psi(j)\}\in
     E(G')\;\text{iff}\;s(j)=1,\;\text{for all}\;j\in[3]\}.
   \end{equation}
   Clearly, $(V'_s)_{s\in2^{[3]}}$ forms a partition of
   $V(G')\setminus\mathrm{Im}(\psi)$. Let $G''$ be the graph obtained by deleting
   the nodes of $G'$ that belong to some $V'_s$ of cardinality at most 3.
   Finally, set $V_s = V'_s \cap V(G'')$ for all $s\in{}2^{[3]}$.
   Thus we have that $(V_s)_{s\in2^{[3]}}$ forms a partition of
   $V(G'')\setminus\mathrm{Im}(\psi)$, each $V_s$ is either empty
    or contains at least four elements and every induced subgraph of $G'$ on five vertices
   belongs to $\C S$.
   \hide{Moreover, by changing the adjacency on at
   most $O(n)$ pairs of nodes, we may assume that for all $s\in{}2^{[3]}$ we
   have that either $V_s=\emptyset$ or $|V_s|\geqslant4$.
   Finally, again by changing the adjacency on at
   most $O(n)$ pairs of nodes, we may assume that 
   \begin{enumerate}
     \item[(Q.i)] Either there is no edge between $V_{(1,0,0)}$ and
     $V_{(0,0,1)}$ or there is some $z\in{}V_{(0,0,1)}$ having at least two
     neighbours in $V_{(1,0,0)}$.
     \item[(Q.ii)] Either there is no edge between $V_{(0,1,0)}$ and
     $V_{(0,0,1)}$ or there is some $z\in{}V_{(0,0,1)}$ having at least two
     neighbours in $V_{(0,1,0)}$.
     \item[(Q.iii)] Either there is no edge between $V_{(1,0,0)}$ and
     $V_{(1,1,0)}$ or there is some $z\in{}V_{(1,1,0)}$ having at least two
     neighbours in $V_{(1,0,0)}$.
     \item[(Q.iv)] Either there is no edge between $V_{(0,1,0)}$ and
     $V_{(1,1,0)}$ or there is some $z\in{}V_{(1,1,0)}$ having at least two
     neighbours in $V_{(0,1,0)}$.
   \end{enumerate}}
\begin{claim}\label{cl:G''}
  The graph $G''$ is a blow-up of the disjoint union of two edges, or 
  a disjoint union of a complete graph and a blow-up of an edge, or
  the disjoint union of a complete graph and an empty graph, or
  the disjoint union of two complete graphs.
\end{claim}

Before we give the proof of Claim~\ref{cl:G''}, let us show how it implies
Lemma~\ref{lm:PawStab}. Observe that an isolated clique can contain at most one
vertex of an $F$-subgraph. Thus if we remove all edges inside such cliques in
$G''$, then we do not decrease the number of $F$-subgraphs.
By Claim~\ref{cl:G''}, the resulting graph $G'''$ is a blow-up of $K_2\cup K_2$.
By Lemma~\ref{lm:PawMinimal}, $G'''$ cannot be a blow-up of $K_2\cup K_1$. This
easily implies that $G''$ itself is a blow-up of $K_2\cup K_2$. Since $G$
and $G''$ are close to each other, Lemma~\ref{lm:PawStab} follows.

\smallskip\noindent{\it Proof of  Claim~\ref{cl:G''}.}
To prove the claim, it suffices to show that
\begin{enumerate}
  \item\label{it:paw.01}  $V_{(1,1,1)}$ is empty,
  \item\label{it:paw.02} both $V_{(1,0,1)}$ and $V_{(0,1,1)}$ are empty,
  \item\label{it:paw.03} $G''[V_{(1,1,0)}]$ is complete,
  \item\label{it:paw.04} both $G''[V_{(1,0,0)}]$ and $G''[V_{(0,1,0)}]$ are empty
  graphs,
  \item\label{it:paw.05} $G''[V_{(0,0,1)}]$ is either complete or empty graph,
  \item\label{it:paw.06} $G''[V_{(0,0,0)}]$ is an empty graph,
  \item\label{it:paw.07} for every $(z_1,z_2)\in V_{(1,0,0)}\times
  V_{(0,1,0)}$, we have that $z_1,z_2$ form an edge in $G''$,
  \item\label{it:paw.08} for every $(z_1,z_2)\in V_{(1,0,0)}\times
  V_{(0,0,0)}$, we have that $z_1,z_2$ do not form an edge in $G''$,
  \item\label{it:paw.09} for every $(z_1,z_2)\in V_{(0,0,0)}\times
  V_{(0,0,1)}$, we have that $z_1,z_2$ form an edge in $G''$,
  \item\label{it:paw.10} for every $(z_1,z_2)\in V_{(0,1,0)}\times
  V_{(0,0,0)}$, we have that $z_1,z_2$ do not form an edge in $G''$,
  \item\label{it:paw.11} there is no edge between $V_{(1,0,0)}$ and
  $V_{(0,0,1)}$, as well as, between $V_{(0,1,0)}$ and $V_{(0,0,1)}$,
  \item\label{it:paw.12} for every $(z_1,z_2)\in V_{(0,0,0)}\times
  V_{(1,1,0)}$, we have that $z_1,z_2$ do not form an edge in $G''$,
  \item\label{it:paw.13} for every $(z_1,z_2)\in V_{(0,0,1)}\times
  V_{(1,1,0)}$, we have that $z_1,z_2$ do not form an edge in $G''$,
  \item\label{it:paw.14} if $V_{(0,0,0)}\neq\emptyset$, then $G''[V_{(0,0,1)}]$
  is an empty graph and
  \item\label{it:paw.15} if $V_{(0,1,0)}\neq\emptyset$ or
  $V_{(1,0,0)}\neq\emptyset$ then $V_{(1,1,0)}=\emptyset$.
\end{enumerate}
\begin{figure}[htb]
\centering \includegraphics[]{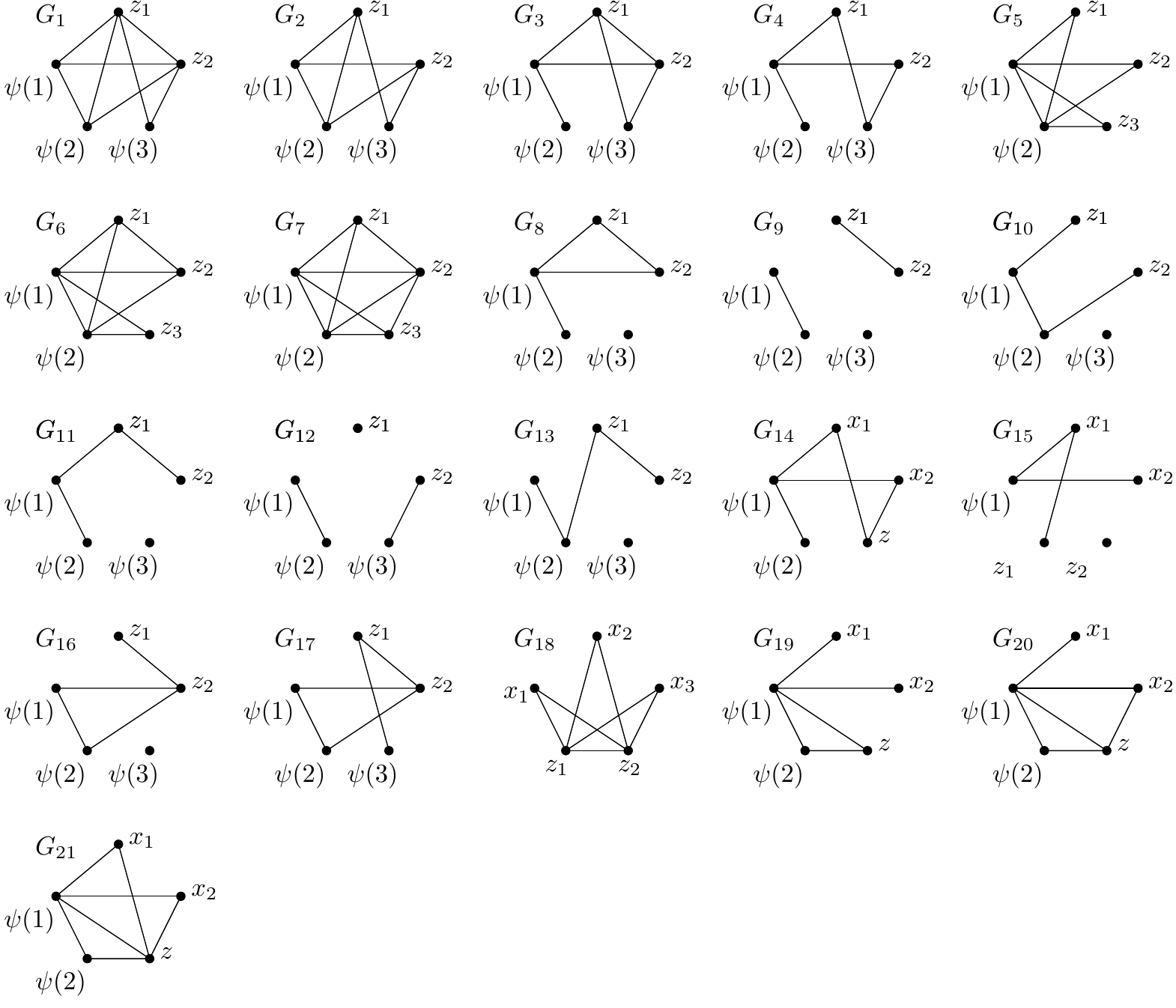}
\caption{Non-sharp graphs that are used in the proof of Claim~\ref{cl:G''}}
\label{fg:3}
\end{figure}
To prove \ref{it:paw.01}, we assume, on the contrary, that $V_{(1,1,1)}$ is
non-empty. Thus $V_{(1,1,1)}$
contains at least four elements. We pick distinct $z_1$ and $z_2$ in
$V_{(1,1,1)}$.
There are two cases: either $z_1$ and $z_2$ form an edge in $G''$ or not. The induced
subgraphs $G''[\mathrm{Im}(\psi)\cup\{z_1,z_2\}]$ of $G''$ are
the graphs $G_1$ and $G_2$ in Figure~\ref{fg:3} respectively. Neither of them belongs
to~$\C S$.

Concerning \ref{it:paw.02}, the arguments justifying that $V_{(1,0,1)}$ and
$V_{(0,1,1)}$ are empty are identical. So we only show that $V_{(1,0,1)}$ is
empty. Again we assume the contrary and pick two distinct elements $z_1$
and $z_2$ in $V_{(1,0,1)}$. Then the induced subgraph of $G''$ on
$\mathrm{Im}(\psi)\cup\{z_1,z_2\}$ is either $G_3$ or $G_4$ in
Figure~\ref{fg:3},
depending on whether $z_1z_2$ forms an edge or not in $G''$. Neither $G_3$ nor $G_4$
belongs to~$\C S$.

Towards \ref{it:paw.03}, assuming the contrary, we pick distinct $z_1,z_2$
and $z_3$ in $V_{(1,1,0)}$ that do not form a triangle. Then the induced
subgraph of $G''$ on $\{\psi(1),\psi(2),z_1,z_2,z_3\}$ is either the graph $G_5$
or $G_6$ or $G_7$ from Figure~\ref{fg:3}, depending on whether $z_1,z_2,z_3$
span zero,
one or two edges respectively. None of these graphs belongs to $\C S$.

Concerning \ref{it:paw.04}, the arguments justifying that $G''[V_{(1,0,0)}]$ and
$G''[V_{(0,1,0)}]$ are empty graphs are identical.
So we only show that $G''[V_{(1,0,0)}]$ is an empty graph.
Assume on the contrary that there exist
  distinct $z_1$ and $z_2$ in $V_{(1,0,0)}$ that form an edge in $G''$.
  Then $G''[\mathrm{Im}(\psi)\cup\{z_1,z_2\}]$ is the graph $G_8$ in
  Figure~\ref{fg:3}
  and does not belong to $\C S$.

To see \ref{it:paw.05}, recall that $V_{(0,0,1)}$ is either
empty or contains at least four elements. If $V_{(0,0,1)}$ is empty then our
claim holds trivially. So let us assume that  $V_{(0,0,1)}$ is of cardinality
at least 4 and pick $z_1,z_2,z_3,z_4\in{}V_{(0,0,1)}$.
Let $H=G''[\{\psi(3),z_1,z_2,z_3,z_4\}]$.  Observe that $\psi(3)$ is of degree 4
in $H$. The only graphs in $\C S$ that contain a node of degree 4 are the star
and the complete graph. Thus $H$ has to be isomorphic to one of these two, yielding
that $G''[\{z_1,z_2,z_3,z_4\}]$ is either an empty or a complete graph,
respectively, on 4 vertices, and $G''[V_{(0,0,1)}]$ is either an empty or a
complete graph, respectively.

To prove \ref{it:paw.06}, we assume the contrary and pick $z_1,z_2$ in
$V_{(0,0,0)}$
that form an edge in $G''$.  Then the induced subgraph of $G''$ on
$\mathrm{Im}(\psi)\cup\{z_1,z_2\}$ is the graph $G_9$ in Figure~\ref{fg:3},
which does not belong to $\C S$.

To prove \ref{it:paw.07}, assuming the contrary, we pick $z_1$ in
$V_{(1,0,0)}$ and $z_2$ in  $V_{(0,1,0)}$ that do not form an edge. Then
graph $G_{10}$ in Figure~\ref{fg:3} is the induced subgraph of $G''$ on
$\mathrm{Im}(\psi)\cup\{z_1,z_2\}$ and does not belong to $\C S$.

To prove \ref{it:paw.08}, assuming the contrary, we pick $z_1$ in
$V_{(1,0,0)}$ and $z_2$ in  $V_{(0,0,0)}$ that form an edge. Then
graph $G_{11}$ in Figure~\ref{fg:3} is the induced subgraph of $G''$ on
$\mathrm{Im}(\psi)\cup\{z_1,z_2\}$ and does not belong to $\C S$.

Similarly, to prove \ref{it:paw.09}, assuming the contrary, we pick $z_1$ in
$V_{(0,0,0)}$ and $z_2$ in  $V_{(0,0,1)}$ that do not form an edge. Then
graph $G_{12}$ in Figure~\ref{fg:3} is the induced subgraph of $G''$ on
$\mathrm{Im}(\psi)\cup\{z_1,z_2\}$ and does not belong to $\C S$.

To prove \ref{it:paw.10}, assuming the contrary, we pick $z_1$ in
$V_{(0,1,0)}$ and $z_2$ in  $V_{(0,0,0)}$ that form an edge. Then
graph $G_{13}$ in Figure~\ref{fg:3} is the induced subgraph of $G''$ on
$\mathrm{Im}(\psi)\cup\{z_1,z_2\}$ and does not belong to $\C S$.

Both assertions in \ref{it:paw.11} follow by identical arguments. So let us
show that there is no edge between
$V_{(1,0,0)}$ and $V_{(0,0,1)}$. First, we show that there is no vertex in one of these sets
having more than one neighbour in the other. There are two cases and the arguments are similar.
So we show that there is no vertex in $V_{(0,0,1)}$ having at least two
neighbours in $V_{(1,0,0)}$.
Indeed, assuming the contrary, we have that there exists $z\in{}V_{(0,0,1)}$ having at least two
neighbours in $V_{(1,0,0)}$, say $x_1,x_2$. Then $G''$ induces on
$\{\psi(1),\psi(2),x_1,x_2,z\}$ the graph $G_{14}$, which does not belong to
$\C S$. Finally, assuming that there is an edge between $V_{(1,0,0)}$ and $V_{(0,0,1)}$,
we can find $x_1 , x_2\in{}V_{(1,0,0)}$ and $z_1 , z_2\in{}V_{(0,0,1)}$ such that
$x_1 , z_1$ form an edge, while $x_2 , z_2$ do not.
Then $G''$ induces on
$\{\psi(1) , x_1 , x_2 , z_1 , z_2\}$ the graph $G_{15}$, which does not belong to
$\C S$.

To prove \ref{it:paw.12}, assuming the contrary, we pick $z_1$ in
$V_{(0,0,0)}$ and $z_2$ in  $V_{(1,1,0)}$ that form an edge. Then the
graph $G_{16}$ in Figure~\ref{fg:3} is the induced subgraph of $G''$ on
$\mathrm{Im}(\psi)\cup\{z_1,z_2\}$ and does not belong to $\C S$.

Similarly, to prove \ref{it:paw.13}, assuming the contrary, we pick $z_1$ in
$V_{(0,0,1)}$ and $z_2$ in  $V_{(1,1,0)}$ that form an edge. Then the
graph $G_{17}$ in Figure~\ref{fg:3} is the induced subgraph of $G''$ on
$\mathrm{Im}(\psi)\cup\{z_1,z_2\}$ and does not belong to $\C S$.

Towards \ref{it:paw.14}, we assume on the contrary that
$V_{(0,0,0)}$ is non-empty and $G''[V_{(0,0,1)}]$ is not an empty
graph.
Thus there exist
 $z_1,z_2\in{}V_{(0,0,1)}$ that form an edge in $G''$. Pick distinct
$x_1,x_2,x_3\in{}V_{(0,0,0)}$. Invoking Items \ref{it:paw.06}
and \ref{it:paw.09}, we have that induced subgraph of $G''$ on
$\{x_1,x_2,x_3,z_1,z_2\}$ is the graph $G_{18}$ in Figure~\ref{fg:3}, which
does not
belong to $\C S$.

Concerning \ref{it:paw.15}, the arguments yielding that $V_{(1,1,0)}=\emptyset$
assuming $V_{(1,0,0)}\neq\emptyset$ are identical to
the ones yielding that $V_{(1,1,0)}=\emptyset$ assuming
$V_{(0,1,0)}\neq\emptyset$. So let us show the first implication.
Assume on the contrary that both $V_{(1,1,0)}$ and $V_{(1,0,0)}$ are non-empty.
We pick $x_1 , x_2\in{}V_{(1,0,0)}$ and $z\in{}V_{(1,1,0)}$. By item \ref{it:paw.04},
we have that $x_1 , x_2$ do not form an edge in $G''$.
Thus the induced
subgraph of $G''$ on $\{\psi(1),\psi(2),x_1,x_2,z\}$ is either the graph $G_{19}$
or $G_{20}$ or $G_{21}$ from Figure~\ref{fg:3}, depending on whether
$z,x_1,x_2$ span zero,
one or two edges respectively. None of these graphs belongs to $\C S$.
This finishes the proof of Claim~\ref{cl:G''} (and thus of
Lemma~\ref{lm:PawStab}).\epf

We have the following strengthening of Lemma \ref{lm:PawMinimal}.

\begin{lemma}
  \label{lm:PawLambdaMaximazer}
  The only maximiser of $\lambda(\blow{K_2\cup K_2}{\cdot})$ is the vector
  $(1/4,1/4,1/4,1/4)$.
\end{lemma}
\bpf
Let $\mathbf{a}=(a_1,a_2,a_3,a_4)$ in $\I S_{4}$ and set $B=K_2\cup K_2$. As we
have already mentioned
\[\lambda(\blow{B}{\mathbf{a}})
=12(a_1a_2+a_3a_4)(a_1+a_2)(a_3+a_4).
\]
To prove that the only maximiser of $\lambda(\blow{B}{\cdot})$ is the vector
$(1/4,1/4,1/4,1/4)$, we show, equivalently, that the map
\[f(x_1,x_2,x_3)=12(x_1x_2+x_3(1-x_1-x_2-x_3))(x_1+x_2)(1-x_1-x_2)\]
with $(x_1,x_2,x_3)$ in
$D=\{(x_1,x_2,x_3):x_1,x_2,x_3\geqslant0\;\text{and}\;x_1+x_2+x_3\leqslant1\}$
admits the vector $(1/4,1/4,1/4)$ as the unique maximiser. By Lemma
\ref{lm:PawMinimal}, no maximiser of $f$ is on the boundary of $D$. Thus, we are
interested in the points belonging to the interior of $D$, where all the
partial derivatives of $f$ vanish. Hence, the following equations should be
satisfied.
\begin{equation}
  \label{eq:PawExtrPartialDev1}
  \frac{\partial f}{\partial
  x_1}=12\big((x_2-x_3)(x_1+x_2)(1-x_1-x_2)+(x_1x_2+x_3(1-x_1-x_2-x_3))(1-2x_1-2x_2)\big)=0,
\end{equation}
\begin{equation}
  \label{eq:PawExtrPartialDev2}
  \frac{\partial f}{\partial
  x_2}=12\big((x_1-x_3)(x_1+x_2)(1-x_1-x_2)+(x_1x_2+x_3(1-x_1-x_2-x_3))(1-2x_1-2x_2)\big)=0
\end{equation}
and
\begin{equation}
  \label{eq:PawExtrPartialDev3}
  \frac{\partial f}{\partial x_3}=12(1-x_1-x_2-2x_3)(x_1+x_2)(1-x_1-x_2) = 0.
\end{equation}
By \eqref{eq:PawExtrPartialDev3} and recalling that we are only interested in
points that belong to the interior of $D$, we get that $x_3=1-x_1-x_2-x_3$,
while by subtracting \eqref{eq:PawExtrPartialDev1} and
\eqref{eq:PawExtrPartialDev2}, we get that $x_1=x_2$. Combining these two we
get, in particular, that
$x_1+x_3=1/2$. Plugging the last three equalities into
\eqref{eq:PawExtrPartialDev1}, we have that
\[\begin{split}
  0
  &=12\big((2x_1-1/2)2x_1(1-2x_1)+(x_1^2+(1/2-x_1)^2)(1-4x_1)\big) = 3(1-4x_1)^3.
\end{split}\]
Thus $x_1=1/4$ and the result follows readily.
\epf

\begin{lemma}
  \label{lm:PawFlipAverse}
  The graph $K_2\cup K_2$ is $\lambda$-flip-averse.
\end{lemma}
\bpf
Set $B=K_2\cup K_2$ and let $B'=\blow{B}{V_1,V_2,V_3,V_4}$ be a blow-up of $B$
on $n$ vertices, with $|V_i|=n/4+O(1)$ for all $i=1,2,3,4$.
Let $i,j\in{}[4]$, $x\in{}V_i$ and $y\in{}V_j$, with $x\neq y$.
It suffices to distinguish the following three cases.

If $ij$ forms an edge in $B$, then the number of $F$-subgraphs in $B'$ (resp.\ $B'\oplus xy$) that use the pair $xy$ is $n^2/4+O(n)$ (resp.\ $n^2/16+O(1)$) and we have that
\begin{equation}
  \label{eq:045}
  \Lambda(B')-\Lambda(B'\oplus xy)=n^2/4-n^2/16+O(n)=3n^2/16+O(n).
\end{equation}

If $i\neq j$ and $ij$ do not form an edge, then
\begin{equation}
  \label{eq:046}
  \Lambda(B')-\Lambda(B'\oplus xy)=3n^2/16-n^2/8+O(n)=n^2/16+O(n).
\end{equation}

If $i=j$, then $B'\oplus xy$ has no copies via $xy$ and
\begin{equation}
  \label{eq:044}
  \Lambda(B')-\Lambda(B'\oplus xy)=n^2/8+O(n).
\end{equation}
By \eqref{eq:044}, \eqref{eq:045} and \eqref{eq:046}, the result follows.
\epf

\begin{lemma}
  \label{lm:PawStrict}
  The graph $K_2\cup K_2$ is $\lambda$-strict.
\end{lemma}
\bpf
We set $\mathbf{a}=(1/4,1/4,1/4,1/4)$ and $B=K_2\cup K_2$. By Lemma
\ref{lm:PawLambdaMaximazer}, we have that $\mathbf{a}$ is the unique maximiser
of $\lambda(\blow{B}{\cdot})$. Thus it suffices to check that $B$ is
$(\lambda,\mathbf{a})$-strict. Indeed, let us fix some $\mathbf{y}$ in
$[0,1]^4$ that maximises $R_\mathbf{a}(\cdot)$. We will show that $\mathbf{y}$
has exactly one non-zero entry which is equal to 1.

Let $B'$ be a balanced blow-up of $B$ of order $n$. Let us denote by $G$ the
graph obtained by attaching to $B'$ a new node $w$ not belonging to $V(B')$
with adjacencies governed by $\mathbf{y}$. In particular, if $y_i=0$ (resp.
$y_i=1$) for some $i\in{}[4]$, then $w$ is attached to no vertex (resp. to all
vertices) in $V_i$.
We also define $\mathcal{H}$ to be the set of all graphs $H$ on $5$ vertices satisfying
$H\sim G[X]$ for some $X\in{}{V(G)\choose 5}$ with $w\in X$.
We have the following claim.
\begin{claim}
  \label{cl:PawStrictGraphs}
  $\mathcal{H}\subseteq\mathcal{S}$.
\end{claim}
\bpf[Proof of Claim \ref{cl:PawStrictGraphs}.]
Let $\e$ be a positive real. We denote by $G'$ the graph obtained by adding
$\e n$ twins of $w$ in $G$. Set $V=V(B')$ and $V'=V(G')\setminus V$. Let
$\mathcal{A}_0={V\choose5}$, $\mathcal{A}_1$ to be the set of all $X$ in
${V(G'\choose5}$ having exactly one element in $V'$ and $\mathcal{A}_2$ the set
of all $X\in{}{V(G'\choose5}$ having at least two elements in $V'$.
Since $\mathbf{y}$ maximises $R_\mathbf{a}(\cdot)$, by Lemma \ref{lm:Strict1},
we have that $R_\mathbf{a}(\mathbf{y})=\lambda(\mathcal{G})$.
Therefore,
\begin{equation}
  \label{eq:047}
  \begin{split}
    \lambda(G')
    &={n+\e n\choose k}^{-1}\sum_{X\in {V(G')\choose k}}\lambda(G'[X])\\
    &=\frac{|\mathcal{A}_0|}{{n+\e n\choose k}}\sum_{X\in
    \mathcal{A}_0}\lambda(G'[X])
    +\frac{|\mathcal{A}_1|}{{n+\e n\choose k}}\sum_{X\in
    \mathcal{A}_1}\lambda(G'[X])
    +\frac{|\mathcal{A}_2|}{{n+\e n\choose k}}\sum_{X\in
    \mathcal{A}_2}\lambda(G'[X])\\
    &\geqslant (1-k\e)\lambda(\C G)+k\e(1-k\e)\lambda(\C G)+O(\e^2)+O(1/n)
  \end{split}
\end{equation}
and hence we get
\begin{equation}
  \label{eq:048}
  \lambda(\C G)-\lambda(G')\leqslant O(\e^2) + O(1/n).
\end{equation}
By \eqref{eq:ai+bi}, there exists some positive real $\eta$
independent from $\e$ and $n$ satisfying
\begin{equation}
  \label{eq:049}
  \lambda(\C G)-\lambda(G')\geqslant \eta p(H,G') +O(1/n)
\end{equation}
for all $H\in{}\mathcal{NS}$. Finally, observe that for every $H\in{}\C H$ we
have that $p(H,G')=\Omega(\e)$.
Thus by \eqref{eq:048}, \eqref{eq:049} and a choice of a sufficiently small
$\e$, $\C H$ and $\mathcal{NS}$ are disjoint.
\epf

\begin{figure}[htb]
    \centering \includegraphics[width=.4\textwidth]{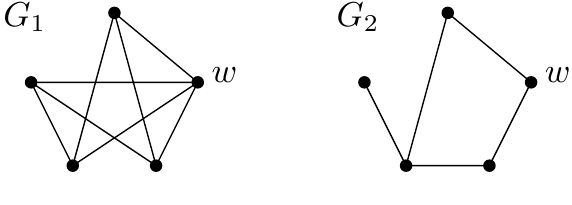}
    \caption{The graphs $G_1$ and $G_2$ used in the proof of
    Lemma~\ref{lm:PawStrict}.}
    \label{fg:4}
  \end{figure}

Since $R_\mathbf{a}(0,0,0,0)=3/16$, $\mathbf{y}$ has at least one
non-zero coordinate. Next, let us observe that $\mathbf{y}$ cannot have two
positive coordinates corresponding to adjacent nodes of $B$. Indeed, assuming
the contrary, we pick four nodes in $V(B')$ adjacent to $w$ and inducing a
balanced complete bipartite graph. Together with $w$, they induce the graph
$G_1$ in Figure~\ref{fg:4}, which does not belong to $\C S$, though it belongs
to $\C
H$ by definition of $\C H$, contradicting Claim \ref{cl:PawStrictGraphs}.
Hence, $\mathbf{y}$ has either one non-positive coordinate, or two positive
coordinates that correspond to non-adjacent nodes $i,j$ in $B$.
Assuming that the second case occurs, picking two nodes adjacent to $w$ from $V_i$, one node
adjacent to $w$ from $V_j$ and one node non-adjacent to $w$ from $V_{i'}$, where $i'$ is the element of $V(B)$
adjacent to $i$ in $B$, we have that these nodes together with $w$ induce the graph
$G_2$ from Figure~\ref{fg:4} that once again does not belong to $\C S$ though
it
belongs to $\C H$, contradicting Claim \ref{cl:PawStrictGraphs}. Therefore, $\mathbf{y}$ has exactly on positive coordinate. Finally, we observe
that the non-zero coordinate of $\mathbf{y}$ is equal to 1. Indeed, let us
assume the contrary and let $y_i$ be the non-zero coordinate of $\mathbf{y}$.
Also let $i'$ be the adjacent node to $i$ in $B$. Picking two nodes in $V_i$
adjacent to $w$, a non-adjacent node to $w$ in $V_i$ and a node in $V_{i'}$,
together with $w$ we induce the graph $G_2$ in Figure~\ref{fg:4} that belongs
to $\C H$
and not to $\C S$ contradicting Claim \ref{cl:PawStrictGraphs}.
\epf

By Lemmas \ref{lm:PawStab}, \ref{lm:PawFlipAverse} and~\ref{lm:PawStrict}, the assumptions of Theorem \ref{th:exact} are satisfied and therefore the
problem is perfectly $K_2\cup K_2$-stable.

\section{Tur\'an problem}\label{se:turan}

This section is devoted to the proof of Theorem \ref{th:turan}. As we have already mentioned
Part 1 of Theorem \ref{th:turan} is known (see \cite[Lemma~2.3]{roberts+scott}), we focus on the proof of the second part of Theorem \ref{th:turan}.
Let
 \begin{equation}\label{eq:v0}
 v_0:=\max\{v(H'):H'\in\C H\},
 \end{equation}
 which is finite as $\C H$ was assumed to be finite. Recall that we defined
 $\C H^\uparrow$ to be the
 collection of graphs obtained by adding missing edges
 to the graphs in $\C H$.
Before we start the proof, we provide an equivalent reformulation of the
property in Part 2 of Theorem \ref{th:turan}. As in the theorem,
let  $m=\min\{\chi(H)\mid H\in\C H\}-1$. Then the following two statements are equivalent.
 \begin{enumerate}
   \item[(i)] There is a constant $D$ such that for every $q$ if we add at least $Dq$
 edges into a part of $K_m^q$ then the obtained graph
 is not $\C H^\uparrow$-free.
   \item[(ii)] There is
   a forest $W$ such that the graph obtained from  $K_m^{q_0}$, $q_0:=v(W)$,
by adding $W$ into one part is not $\C H^\uparrow$-free. \end{enumerate}

    Let us first assume (i) and prove (ii).
 Let $Z$ be a graph with minimum degree at least $2 D$ and girth strictly greater than $v_0$. Set $q:=v(Z)$.
 Also, let $V_1,...,V_m$ be disjoint sets, each of cardinality $q$, and let $G$ to be the graph obtained from $K_m(V_1,...,V_m)$ by adding a copy of $Z$ in $V_1$.
 Since $Z$ is of minimum degree $2D$, we have that $Z$ contains at least $Dq$ edges. By (i),
 $G$ has a (not necessarily induced) subgraph  $H\in{}\C H$. Since $v(H)\le v_0$, we conclude that $H[V_1]$ contains no cycle and thus $H$ is as desired.

 Assuming (ii), we claim that (i) holds with $D:=q_0$. This is a consequence of the well-known fact that if $G$ is a graph with $|E(G)|>(q_0-1)v(G)$, then $G$
 contains
 a copy of the forest $W$ (not necessarily as an induced subgraph).  Indeed, by e.g.~\cite[Theorem~2.5]{BondyMurty:gt}, $G$ contains a non-empty subgraph $G'$ of minimum
 degree at least $q_0$ where the required copy of $W$ can be easily found.

Moreover, we will need the following result, which follows by the Ramsey Theorem \cite{Ramsey} and elementary probabilistic estimates
(see e.g.\ \cite[Lemma 2.7]{Dodos+Kanellopoulos+Tyros} for a proof).
\begin{lemma}
  Let $\varepsilon , \theta$ be reals with $0 < \theta < \varepsilon$ and $\ell_1 , \ell_2$ be positive integers with $\ell_1 < \ell_2$.
  Then there exists a positive integer $\ell_3$ with the following property. For every probability space $(\Omega, \Sigma, \mu)$ and every sequence
  $(A_j)_{j=1}^{\ell_3}$ such that $\mu(A_j)\ge\varepsilon$ for all $j\in{}[\ell_3]$, we have that there exists a subset $L$ of $[\ell_3]$ of cardinality $\ell_2$ such that for every subset $K$ of $L$ of cardinality $\ell_1$ we have that $$\mu\Big(\bigcap_{j\in K} A_j\Big) \ge \theta^{\ell_1}.\qed$$
\end{lemma}
An iterated use of the above lemma yields the following, which we will use in the proof of the second part of Theorem \ref{th:turan}.
\begin{lemma}\label{lm:turan}
  Let $\varepsilon , \theta$ be reals with $0 < \theta < \varepsilon$ and $q,\ell$ be positive integers.
  Then there exists a positive integer $k=k(q,\ell,\theta,\varepsilon)$ with the following property. Let $(\Omega_1, \Sigma_1, \mu_1),..., (\Omega_q, \Sigma_q, \mu_q)$ be probability spaces and for each $i\in{}[q]$ let
  $(A^i_j)_{j=1}^{k}$ be a sequence in
  $\Sigma_i$ such that $\mu_i(A_j^i)\ge\varepsilon$ for all $j\in{}[k]$. Then there exists a subset $L$ of $[k]$ of cardinality $\ell$ such that for every $i\in{}[q]$ we have that $$\mu_i\Big(\bigcap_{j\in L} A^i_j\Big) \ge \theta^{\ell}.\qed$$
\end{lemma}

\bpf [Proof of Part 2 of Theorem \ref{th:turan}.]
Recall that the theorem of Erd\H os~\cite{erdos:67} and
Simonovits~\cite{simonovits:68} states that the problem is classically
stable with $B=K_m$. Thus the only twin-free graph $B$, which can have  the property
that the problem is robustly $B$-stable, is $K_m$.
Let $t_m(n)$ be the maximum size of a $K_m$-blow-up of order $n$; it is easy to see that the maximum is attained if and only if any two part sizes differ at most by $1$.

According to the discussion in the beginning of this section, it suffices to prove equivalence between robust $K_m$-stability and Condition (i) stated above. If Condition~(i) fails, then for each $D$ we can construct an $\C H^\uparrow$-free graph $G_D$  by adding
$Dq$ edges to $K_m^q$ for some~$q=q(D)$. This graph $G_D$ of order $n:=mq$ exceeds $t_m(n)$, the maximum size of a $K_m$-blow-up on $n$ vertices, by $Dn/m$. Thus problem is not robustly $K_m$-stable.


 Let us show the converse direction. Let $D$ satisfy
 Condition (i) and define $v_0$ by~\eqref{eq:v0}. Given $\C H$ and $D$,
 we choose positive constants in this order
 $$
 c\gg c_3\gg c_2\gg c_1\gg c_0,
 $$
 each being sufficiently small depending on the previous ones.
 Assume on the contrary that the problem is
 not robustly $K_m$-stable. Hence,
  there exists an $\C H^\uparrow$-free graph
  $G$ with $n\ge 1/c_0$ vertices satisfying
  \begin{equation}
    \label{eq:turan08}
    t_m(n) - e(G) + n < c_1 \Dedit( G , \blow{K_m}{} ).
  \end{equation}
  Let $V_1,...,V_m$ be a max-cut partition of $V(G)$ and set $T:=\blow{K_m}{V_1,...,V_m}$. Since $e(G)\ge t_m(n)-2c_1{n\choose 2}$,
  we have by the Erd\H os-Simonovits Stability Theorem~\cite{erdos:67,simonovits:68} that
  \begin{equation}
    \label{eq:turan08.1}
    |E(G) \bigtriangleup E(T)|\le c_2 {n \choose 2}.
  \end{equation}
  It routinely follows that
  \begin{equation}
    \label{eq:turan08.2}
    (1/m-c_3)n\le|V_i|\le(1/m+c_3)n,\quad \mbox{for all $i\in{}[m]$.}
  \end{equation}

   Next we observe that in each $G[V_i]$ there are only a few vertices of high degree. More precisely, we have the following claim.
  \begin{claim}\label{cl:turan_02.1.1}
  For every $i\in{}[m]$, the induced subgraph $G[V_i]$ has at most $k(m, 2D, 2mc/3, 3mc/4)$ vertices of degree at least $c n$, where $k()$ satisfies Lemma \ref{lm:turan}.
\end{claim}
\bpf
We set $k=k(m, 2D, 2mc/3, 3mc/4)$ and assume on the contrary that there exist $i_0\in{}[m]$ and $x_1, ..., x_k\in{}V_{i_0}$ such that the degree of each $x_j$ in $G[V_{i_0}]$ is at least $c n$. By the max-cut property of $V_1,...,V_m$ we have for each
$i\in{}[m]$ and $j\in{}[k]$ that the set of all neighbours of $x_j$ in $V_i$, which we denote by $A_j^i$,
is of cardinality at least $c n$ and therefore, by \eqref{eq:turan08.2}, of uniform density at least $\frac{c}{1/m + c_3} \ge 3m c / 4$. Applying Lemma \ref{lm:turan},
we obtain a subset $L$ of $[k]$ of cardinality $2D$ such that for each $i\in{}[m]$, setting $B_i:=\bigcap_{j \in L}A_j^i$ (which is the set of vertices in $V_i$ that are $G$-adjacent to $x_j$ for all $j\in{}L$),
we have that $$
|B_i|\ge (2mc/3)^{2D}\,|V_i|\ge (c^{2D} / 2^{2D}) n.$$

We pick arbitrary subsets $Y_1,...,Y_m$ of $B_1,...,B_m$ respectively, of cardinality $(c^{2D} / 2^{2D}) n $ each.
We set $Y:=Y_1 \cup\dots\cup Y_m$ and $Z:=\{x_j : j \in L\}$.
Observe that $G[Y]$ cannot contain a copy (not necessarily induced) of $K_m^{4D}$.
Indeed, assume on the contrary that there exist pairwise disjoint $4D$-subsets $W_1,...,W_m$ of $Y_1,\dots,Y_m$, respectively, such that $E(G)\supseteq E(\blow{K_m}{W_1,...,W_m})$.
Let $W_1'$ be the set obtained by deleting $2D$ vertices from $W_1$ and adding the set $Z$.
Then $E(G)$ is a superset of $E(G[W_1'])\cup E(\blow{K_m}{W_1',W_2,...,W_m})$. Observing that $G[W_1']$ contains at least $2D \cdot 2D =D\cdot |W_1'|$ edges, we get that $G$ is not $\mathcal{H}^\uparrow$-free, a contradiction.

Thus, for every choice of a $4D$-subset $W_j$ of $Y_j$, for $j=1,...,m$, there should be
at least one missing edge (that is, an edge of $T$ but not of $G$). Notice that there are ${c^{2D} n / 2^{2D}  \choose 4D} ^ m$ choices of $(W_1,\dots,W_m)$.
On the other hand,
a missing edge can be overcounted at most
\begin{equation}
  {c^{2D} n / 2^{2D}-1 \choose 4D-1}^2{c^{2D} n / 2^{2D}\choose 4D}^{m-2}=\frac{(4D)^2}{(c^{2D} n / 2^{2D})^2}{c^{2D} n / 2^{2D} \choose 4D}^m
\end{equation}
times. Thus $ E(T) \setminus E(G) $ is of cardinality at least $ ( c^{4D} / 2^{4+4D} D^2 ) n^2  $ contradicting \eqref{eq:turan08.1}.
\epf

We set $K:=m\cdot k(m, 2D, 2mc/3, 3mc/4)$. Let $U'$ be the set all vertices having at least $c n$ neighbours  within their part.
By Claim \ref{cl:turan_02.1.1} we have that
\begin{equation}\label{eq:turan09}
  |U'|\le K
  \le c_2 n.
\end{equation}
  We also set $U''$ to be the set of all vertices $x$ in $V(G)\setminus U'$ such that $d_T(x) - d_G(x) \ge c n$. By \eqref{eq:turan08.1}, we get that
  \begin{equation}
    \label{eq:turan10}
    |U''|\le \frac{c_2}{c}\, n
    \le (c_3-c_2) n.
  \end{equation}
  Thus, setting $\mathcal{E}''$ to be the set of all pairs $e$ of vertices in $V(G)$ satisfying $e \cap U'' \neq \emptyset$, we have that
  \begin{equation}\label{eq:turan12}
    | E(T) \cap \mathcal{E}'' | - | E(G) \cap \mathcal{E}'' |
    \ge |U''|\, c n - {|U''| \choose 2} \stackrel{ \eqref{eq:turan10} }{ \ge } |U''| (c - c_3) n 
    \ge  \frac{c}{2}\, |U''|\, n.
  \end{equation}
  Moreover, setting $U:=U' \cup U''$, by \eqref{eq:turan09} and \eqref{eq:turan10}, we have that $|U| \le c_3 n$.
  Also,  set $V':=V\setminus U$ and $V'_i:=V_i\setminus U$ for each $i\in{}[m]$.
  We have the following claim.
  \begin{claim}\label{cl:turan_02.1}
   Let $i\in{}[m]$ and $X$ be a subset of $V' \setminus V_i$ with at most $v_0$ elements. Then $V'_i$ has at least $(1-3 c m v_0 ) |V_i'|$ vertices $G$-adjacent to every node in~$X$.
 \end{claim}
 \bpf
 For every $x\in{}V' \setminus V_i$ we have the following. Let $j$ be the unique element of $[m]$ satisfying $x\in V_j$. Since $x\not\in U''$, we have that
 $d_G(x)\ge d_T(x) -  c n$. Invoking the fact that $x$ has at most $c n$ $G$-neighbours in $V_j$, since $c \not \in U'$, and $x$ is $T$-adjacent to all vertices in $V'_i$, it follows that $x$ is $G$-adjacent to all  but at most $2 c n$ vertices in $V'_i$. By \eqref{eq:turan08.2}, \eqref{eq:turan09} and \eqref{eq:turan10}, we get that $n / m \le |V'_i| / (1 - 2 c_3 m)$
 and therefore $x$ is $G$-adjacent to at least
 \[|V_i'| - 2cm \frac{n}{m} \le \Big( 1 - \frac{2 c m}{1 - 2 c_3 m} \Big)|V_i'|
 \le ( 1 - 3 c m )|V_i'| \]
 vertices in $V_i'$. Since $X$ has at most $v_0$ elements, the claim follows.
 \epf

 Next, we observe that in each $V_i'$ we have a few edges. In particular, we have the following.

  \begin{claim}
    \label{cl:turan_03}
    For every $j\in{}[m]$, we have that $G[V_j']$ contains less than $(5D/2m)n$ edges.
  \end{claim}
  \bpf
    We assume on the contrary that there is some $j\in{}[m]$ such that
    $G[V_j']$ contains at least $5D(n/2m)$ edges. Without loss of generality,
    let $j=1$.  Let $X_1$ be a random subset
    of $V'_1$ of size $n/2m$.
    Then the expected number of edges in $G[X_1]$ is
    at least
    $$5 D ( n / 2 m ){ |V'_1| - 2 \choose n / 2 m - 2 }{ |V'_1| \choose n / 2 m }^{-1}\stackrel{\eqref{eq:turan08.2}}{>}D \cdot n / 2 m
     .$$
    We pick
    $X_1$ so that $G[X_1]$ has at least $D \cdot n / 2 m$ edges and, for each $i\in{}\{1,\dots,m\}$,  we pick an arbitrary $(n/2m)$-subset $X_i$ of $V'_i$.
    We set $F$ to be the graph obtained by adding to $\blow{K_m}{X_1,...,X_m}$ the edges of $G[X_1]$.
    By the choice of $D$
    there is an injective homomorphism $f$, that is, an injective map sending edges to edges, from
    some $H\in{}\C H$ into $F$.
    We will arrive to a contradiction by constructing an injective homomorphism $f'$ from $H$ into $G$.
    To this end, we set $Y_i:=\{h\in H: f(h) \in X_i\}$ for each $i\in{}[m]$. We inductively define $f'$ on each $Y_i$.  For every $h\in{}Y_1$ we set $f'(h)=f(h)$. Then, for each $i=2,...,m$, assuming that $f'$ has been defined on $\bigcup_{j=1}^{i-1}Y_j$, we extend $f'$ on $Y_i$ by using arbitrary elements of $V_i'$
    which are adjacent in $G$ to every element of  $f'(\cup_{j=1}^{i-1}Y_j)$. Claim~\ref{cl:turan_02.1} guarantees that such a selection is feasible. It follows easily that $f'$ is indeed an injective homomorphism. Thus its image contains a (not necessarily induced) copy of $H$, which contradicts that $G$ is $\mathcal{H}^\uparrow$-free.
  \epf

We have
\begin{eqnarray*}
     t_m(n)-e(G) & \ge& e(T)-e(G) \\
       & \ge& e( T[V'] ) - e( G[V'] ) + | E(T) \cap \mathcal{E}''  |
            - | E(G) \cap \mathcal{E}''  | - |U'|\,n \\
       & \stackrel{\eqref{eq:turan09},\eqref{eq:turan12}}{\ge}& e( T[V'] ) - e( G[V'] ) + \frac{c}{2}\, |U''|\,n - Kn\\
       & \stackrel{\text{Claim}\;\ref{cl:turan_03}}{\ge}& | E( T[V'] ) \bigtriangleup E( G[V'] ) |
            + \frac{c}{2}\, |U''|\,n - (5 D + K) n \\
       & \stackrel{\eqref{eq:turan09}}{\ge}& | E( T[V'] ) \bigtriangleup E( G[V'] ) | + |U'|\,n
            + \frac{c}{2}\, |U''|\,n - (5 D + 2 K) n \\
       & \ge& \frac{c}{2}\, \big( | E( T[V'] ) \bigtriangleup E( G[V'] ) | + |U|\,n \big) - (5 D + 2 K) n\\
       & \ge& \frac{c}{2}\, | E( T ) \bigtriangleup E( G ) |  - (5 D + 2 K) n
       \ \ge\  \frac{c}{2}\,\Dedit(G,\blow{K_m}{})- (5 D + 2 K) n.
  \end{eqnarray*}
 By combining this with \eqref{eq:turan08},  we get
\[\frac{1}{c_1}\Big(t_m(n)-e(G) + \frac{n-1}{2}\Big)
\le \Dedit(G,\blow{K_m}{})\le \frac{2}{c}\big( t_m(n) - e(G) + (5 D + 2 K) n \big),\]
which is a contradiction since we have assumed that $c_1 \ll c$.
\epf

\section{Concluding Remarks}

Theorem \ref{th:turan} implies that the three notions of stability introduced in Section~\ref{Intro} are non-equivalent. Indeed, let $K_{a,b,c}$ denote the complete 3-partite graph with part sizes $a,b,c$. Then the Tur\'an problem  $\ex(n,K_{2,2,2})$ is classically $K_2$-stable (by~\cite{erdos:67,simonovits:68}) but not robustly $K_2$-stable by Part~2 of Theorem \ref{th:turan}. (Namely, one can add a  $C_4$-free bipartite graph of size $\Omega(q^{3/2})$ into one part of $K_2^q$, which will not violate the property of being $\{K_{2,2,2}\}^\uparrow$-free.) Also, $\ex(n,K_{2,2,1})$ is robustly but not perfectly $K_2$-stable by Theorem \ref{th:turan}.



Theoretically, one should be able to write a computer code that takes as input only a family $\C F$ of twin-free graphs and $\Lambda$ and then tries to figure out everything else (namely $B$, $\V a$, $N$, and $\C C$) automatically. For lower bounds, computer can enumerate all small $B$
such that $\blow{B}{}$ is $\C F$-free  and then
use Gr\"obner bases calculations to
calculate $\lambda(\blow{B}{})$, thus identifying best possible $B$.
For upper bounds, computer may start with
largest feasible $N$ (which is $8$ for graphs unless $\C G$ is rather structured),
outputting some floating-point number $c$ as an upper bound. Furthermore, if $c$ seems to coincide with $\lambda(\blow{B}{})$, then the steps of finding smallest $N$ that works and rounding the solution (using
$B$ as conjectured extremal configuration) could be also automated.
However, the
human intuition (based on various heuristics, symmetries, structure of admissible graphs, etc) is usually superior to the brute force search for plausible extremal
configurations. Of course, the more powerful combination would be when computer search is
restricted to a narrow set of plausible examples suggested by the user. It would be
interesting to write a computer code that has this wider functionality and yet requires
little coding from the user.

If the maximiser $\V a$ of $\lambda(\blow{B}{})$ is unique (up to symmetry), one may be tempted
to define another version of stability where one wishes to relate
$\Lambda(G)-\Lambda(n,\C G)$ to the distance from $G$ to $\blow{B}{V_1,\dots,V_m}$
with $|V_i|=a_i/n+O(1)$. However, here the dependence is in general worse.
For example, consider the Tur\'an problem for triangle which is perfectly $K_2$-stable.
Here the optimal $\V a$ is unique: $(1/2,1/2)$. However, for
$G:=K_{(1/2-\e)n,(1/2+\e)n}$ we have that $\lambda(n,\C G)-\lambda(G)=O(\e^2)$, which is much smaller than $\dedit(G,K_{n/2,n/2})=\Omega(\e)$ when $\e\to 0$.

\renewcommand{\baselinestretch}{1.0}

\end{document}